\documentclass{articleChaosCP}
\usepackage{style}
\usepackage[normalem]{ulem}
\usepackage{xfrac}
\DeclareMathAlphabet{\mathpzc}{OT1}{pzc}{m}{it}

\newcommand*{\Scale}[2][4]{\scalebox{#1}{$#2$}}%

\hypersetup{pdfauthor={Name}}
     
\begin{document}
\title{Breakdown of homoclinic orbits to $L_1$ of the hydrogen atom in a circularly
       polarized microwave field}
\author{%
\begin{minipage}[c]{0.4\linewidth}
\centering
Amadeu Delshams$^{1}$\\
\href{mailto:Amadeu.Delshams@upc.edu}{\texttt{Amadeu.Delshams@upc.edu}}
\end{minipage}
\begin{minipage}[c]{0.4\linewidth}
\centering
Merc\`e Oll\'e$^{2}$\\
\href{mailto:Merce.Olle@upc.edu}{\texttt{Merce.Olle@upc.edu}}
\end{minipage}\\[0,2in]
\begin{minipage}[c]{0.4\linewidth}
\centering
Juan Ramon Pacha$^{3}$\\
\href{mailto:Juan.Ramon.Pacha@upc.edu}{\texttt{Juan.Ramon.Pacha@upc.edu}}
\end{minipage}
\begin{minipage}[c]{0.4\linewidth}
\centering
\'Oscar Rodr\'{\i}guez$^{4}$\\
\href{mailto:Oscar.Rodriguez@upc.edu}{\texttt{Oscar.Rodriguez@upc.edu}}
\end{minipage} \\ 
\begin{minipage}[c]{0.85\linewidth}
\small
\centering
\vspace{2\baselineskip}\noindent
$^{1}$Lab of Geometry and Dynamical Systems and IMTech,\\
Universitat Polit\`ecnica de Catalunya (UPC) and Centre de Recerca Matem\`atica (CRM),\\ Barcelona, Spain\\
\vspace{0.5\baselineskip}
$^{2}$Dept.~of Mathematics and IMTech, UPC and CRM,\\
 Barcelona, Spain
\\
\vspace{0.5\baselineskip}
$^{3}$Dept.~of Mathematics, UPC,\\
 Barcelona, Spain
 \\
\vspace{0.5\baselineskip}
$^{4}$Dept.~of Mathematics and IMTech, UPC,\\
Barcelona, Spain
\end{minipage}
\date{\today}
}

\maketitle

\begin{abstract}
We consider the Rydberg electron in a circularly polarized microwave field, whose dynamics is described by a 2 d.o.f. Hamiltonian, which is a perturbation of size $K>0$ of the standard rotating Kepler problem. In a rotating frame, the largest chaotic region of this system lies around a center-saddle equilibrium point $L_1$ and its associated invariant manifolds. 
We compute the distance between stable and unstable manifolds of $L_1$ by means of a semi-analytical method, which consists of combining normal form, Melnikov, and averaging methods with numerical methods. Also, we introduce a new family of Hamiltonians, which we call \emph{Toy CP systems}, to be able to compare our numerical results with the existing theoretical results in the literature.
It should be noted that the distance between
these stable and unstable manifolds is exponentially small in the perturbation parameter $K$ (in analogy with the $L_3$ libration point of the R3BP).
\end{abstract}


\section{Introduction}
\subsection{State of the art of the CP Problem}
The \emph{CP Problem} consists of the motion in an hydrogen atom placed in an
external circularly polarized microwave field.  This system has been intensively
studied, both
experimentally~\cite{PhysRevLett.64.511,nature-569-75-77-2019,PhysRevA.90.063402},
as theoretically, from
classical~\cite{PhysRevA.47.R1612,PhysRevA.50.5077,Ganesan,PhysRevLett.84.610,OP2018}
and quantum~\cite{Ganesan,ZakrzewskiDGR1993,Bialynicki-BirulaKE1994} points of
view, as well as numerically
\cite{PhysRevA.45.R2678,BuchleitnerDG1995,PhysRevA.51.1508,
PhysRevLett.74.1720,BrunelloUF1996,BarrabesOBFM2012,O2018,BOR2024}.
In a rotating frame (see Section~\ref{sec:CPmodel}), the CP problem reads as a 2
degree-of-freedom (d.o.f.) Hamiltonian system. For a small enough field strength
or a large enough angular frequency of the microwave field, the CP problem turns
out to be just a perturbation of some size, say $K>0$, of the rotating Kepler
  problem, and by KAM theory the phase space is pretty full of 2-dimensional
  invariant tori with non-conmensurable frequencies. These 2-dimensional KAM
  tori preclude the existence of unstable trajectories, and in particular of
  ionization, that is, the loss of the electron, which turns out to be the main
  topic dealt with in all these studies. Consequently, ionization only takes
  place when all KAM invariant tori break, which happens when the perturbation
  parameter $K$ is not so small.

However, for small $K$, the motion is not completely regular, and erratic
trajectories can appear. Such kind of chaotic motion takes place outside the
KAM tori, and is typically associated with the existence of saddle invariant objects.
There are, at least, three zones where saddle objects may appear in the CP
problem: (i) near collision with the nucleus, inside the first KAM torus,
(ii)~very far of the nucleus, say at infinity, outside the last KAM torus, and
(iii)~in the resonant zones, between the KAM invariant tori.

The aim of this work is to study the largest chaotic resonance zone, which
happens to be associated with the lowest-order commensurability of the
frequencies of the unperturbed invariant tori, and takes place around a
saddle-center equilibrium point in a rotating frame, called $L_1$, together with
its associated invariant unstable and stable manifolds. Saddle Lyapunov periodic
orbits emerge from $L_1$, which are very important, indeed crucial, for the
global behavior of the system~\cite{BarrabesOBFM2012}. Indeed, it is a standard
principle in dynamical systems that the saddle invariant objects and their
associated invariant manifolds are the main landmarks for understanding the global
dynamical behavior of a dynamical system.

Using an adequate cross-section, we measure the distance between the two
branches of unstable and stable asymptotic trajectories to $L_1$, which turns out to be positive and exponentially small with respect to the perturbation
parameter $K$. This is a consequence of the fact that $L_1$ is a \emph{weak}
saddle equilibrium point, that is, its real characteristic exponents are
$\pm\Ord{\sqrt{K}}$, so they  tend to zero when the size of the perturbation $K$
tends to zero. The measure of this distance is performed both
analytically~\eqref{eq:homoL1Anal} and numerically~(\ref{eq:NumerDistance}--\ref{eq:ResonantAsympt}) with
a fairly good agreement, although its exponentially smallness makes difficult a
very precise fitting.

The numerical measure requires using high-precision computations combined with
the parameterization method (described in Appendix \ref{appendix:parmeth}) to get a local computation of the invariant manifolds,
as well as Taylor method 
 (implemented on a robust,
fast and accurate package by Jorba and Zou
\cite{JorbaZ2005}) 
to integrate the system without a significant loss of
precision. Multiple precision computations have been carried out using the \texttt{mpfr} library (see~\cite{10.1145/1236463.1236468}). 

Analytic approximations are also obtained for the invariant manifolds of $L_1$. For
this, we simply introduce action-angle variables, which are nothing more than
the Delaunay variables for the Kepler problem, combined with a subsequent change
to the Poincaré variables to avoid the degeneracy of the circular solutions of
Kepler's problem. In these variables, the CP problem consists of a perturbation
of an integrable system formed by an isochronous rotor plus an \emph{amended}
pendulum. This \emph{amended} pendulum depends on the perturbation parameter and
gives rise to two \emph{separatrices}, which also depend on the perturbation
parameter.

When the perturbation is taken into account, we introduce a
direct approximation of the splitting of these separatrices, based on the study of the variational equations along them (this is the so-called Melnikov method) and we compare the distance between these approximate invariant manifolds with the one obtained with the  
aforementioned numerical calculations, performed both in synodic coordinates~\eqref{eq:NumerDistance} for the rotating frame (in the Poincar\'e section $y=0$), and in Poincar\'e variables~\eqref {eq:ResonantAsympt}.

We have also verified that to have a good analytical approximation, it is
crucial to choose the right amended pendulum. For this reason we also present an
additional approximation equivalent to Melnikov method, based on the application
of the averaging method. We check that the choice of an adequate degree of the
amended pendulum as a function of the perturbation  parameter is totally
necessary to obtain a good fit with the numerical calculations.

Unfortunately, there are no theoretical results in the literature that can be directly applied to prove analytically the numerical results found for the CP problem, due to the fact that the perturbation presents a high-order singularity, and that the integrable Hamiltonian where perturbative theoretical methods should be applied  depends essentially on the perturbative parameter $K$. To discuss whether Melnikov method, when applied to a Hamiltonian given by an integrable part plus a perturbation, correctly predicts the found fitting formulas, it is necessary to consider both the integrable part and the perturbation.

This discussion has led us to the analysis of a more general Hamiltonian that we have called \emph{Toy CP problem}~\eqref{eq:ToyCP}, which, apart from $K$, also depends on two other parameters $a$ and $m$. More details are provided in subsection \ref{subsec14}. The observed phenomenology opens the door to a new
world of theoretical studies.

It has to be noticed that the distance studied is much smaller in one of the
separatrix, the \emph{internal} one, than in the other one, the \emph{external}
separatrix. For this reason, in this paper we focus on measuring the distance between the
external asymptotic trajectories to $L_1$.

This is a first step to measuring the distance between the two branches of
asymptotic unstable and stable surfaces to the Lyapunov periodic orbits close to
$L_1$, and to show that they intersect along only two transverse homoclinic
orbits. This computation will be performed in a future paper, and will provide a measure
of the chaotic region close to $L_1$, which is proportional to the
transversality of these branches.

It is worth remarking that the CP problem is similar to the Planar Circular
Restricted 3-Body Problem (R3BP), since both are perturbations of the rotational
Kepler problem. Since the equilibrium point $L_1$ of the CP problem is weakly
hyperbolic, it is also similar to the libration point $L_3$ of the R3BP, where
similar features happen and theoretical results are
available~\cite{BaldomaGG2021Arxiv-I, BaldomaGG2021Arxiv-II}. Therefore, we have
applied a number of tools coming from Celestial Mechanics or, more generically,
Hamiltonian systems, like invariant manifolds, Melnikov method, averaging method,
etc.

\subsection{The model for the CP problem}\label{sec:CPmodel}
In the simplest case (assuming \emph{planar} motion for the electron) the
\emph{classical} motion is governed by a system of two 2nd-order ODE
\[
\begin{split}
 \ddot{X} &= -\dfrac{X}{R^{3}}
 		-F\cos \left(\oomega s\right), 
			\qquad R^2=X^2+Y^2,\\
 \ddot{Y} &= -\dfrac{Y}{R^{3}}
	 	-F\sin \left(\oomega s\right), 
			\qquad \dot{} =\frac {\df{}}{\df{s}},
\end{split}
\]
where $\oomega > 0$ is the \emph{angular frequency
of the microwave field} and $F>0$ is the \emph{field strength}.

This system can be written as a periodic in time 2 d.o.f Hamiltonian

\begin{displaymath}
 H (X,Y,P_X,P_Y) =
    \frac{1}{2}\left(P_X^2 + P_Y^2 \right)-\frac{1}{R} +
    F\left(X \cos\left(\oomega  s\right) +
      Y\sin\left(\oomega  s\right)\right).
\end{displaymath}
As in the R3BP (Restricted 3-Body Problem), we can get rid of the time
dependence introducing rotating coordinates
$\left(\widetilde{X},\widetilde{Y},P_{\widetilde{X}},
P_{\widetilde{Y}}, s\right)$
\[
\begin{pmatrix}
	\widetilde{X}\\\widetilde{Y}\end{pmatrix}=R(\oomega s)
	\begin{pmatrix}
		X\\
	     	Y
	\end{pmatrix},\quad
\begin{pmatrix}
		P_{\widetilde{X}}\\
	     	P_{\widetilde{Y}}
\end{pmatrix} = R(\oomega s)
	\begin{pmatrix}
		P_X\\P_Y
        \end{pmatrix},\quad 
		R(\oomega s) = \begin{pmatrix}
					\cos \oomega s&\sin \oomega s\\
				       -\sin \oomega s&\cos \oomega s
			       \end{pmatrix}
\]
to get an \emph{autonomous} Hamiltonian
\[
 H \left(\widetilde{X},\widetilde{Y},
 P_{\widetilde{X}},P_{\widetilde{Y}}\right) =
    \frac{1}{2}\left(P_{\widetilde{X}}^2 + P_{\widetilde{Y}}^2 \right)
   	 -\frac{1}{\widetilde{R}}
         -\oomega \left(\widetilde{X}P_{\widetilde{Y}}
	 -\widetilde{Y}P_{\widetilde{X}}\right) 
    	 +F \widetilde{X},
\]
where $\widetilde{R}^2=\widetilde{X}^2+\widetilde{Y}^2$.  To get rid of the
angular frequency $\oomega$ we perform the change
$(\widetilde{X},\widetilde{Y})=a(\xs,\ys)$,
$\left(P_{\widetilde{X}},P_{\widetilde{Y}}\right)=a\oomega(\pxs,\pys)$, which is a
symplectic transformation with multiplier $a^2 \oomega$ with $a=1/\oomega
^{2/3}$, plus a scaling in time $t=\oomega s$ to get a new Hamiltonian in the
rotating and scaled coordinates $(\xs,\ys,\pxs,\pys)$, usually also called synodic
coordinates in the literature.
\begin{equation}
\label{eq:NumericalH}
H(\xs,\ys,\pxs,\pys)=\frac12(\pxs^2+\pys^2)-\frac1{r}-(\xs\pys-\ys\pxs)+K \xs,\qquad 
		r=\sqrt{\xs^2+\ys^2},
\end{equation}
where $K=F/\oomega ^{4/3}>0$, with associated Hamiltonian equations
\begin{equation}
\label{eq:NumericalEquations}
\begin{aligned}
\dot \xs &= \pxs +\ys, &  \dot{\pxs}&=\pys-\frac{\xs}{r^3}-K, \\
\dot \ys &= \pys -\xs,&  \dot{\pys}&=-\pxs-\frac{\ys}{r^3}.
\end{aligned}
\end{equation}
It is very important and useful to realize that this Hamiltonian system is
invariant under the reversibility
\begin{equation}
\label{eq:reversibility}
(t,\xs,\ys,\pxs,\pys)\mapsto (-t,\xs,-\ys,-\pxs,\pys).
\end{equation}
Note that for $K=0$ we get just the \emph{rotating Kepler Hamiltonian}. Although
the parameter $K$ only appears, additively, in the equation of the momentum
$\pxs$ of equations~\eqref{eq:NumericalEquations}, for $K>0$ the CP system is
\emph{much more complicated} than the rotating Kepler problem, as we will
describe below.

In the literature, it is also typical to find the  system associated with the Hamiltonian \eqref{eq:NumericalH} written as a system of second-order ODE
\begin{equation}\label{eq:nousist}
    \left\{
    \begin{aligned}
        &\ddot{\xs} - 2\dot{\ys} = \Omega_\xs(\xs,\ys),\\
        &\ddot{\ys} + 2\dot{\xs} = \Omega_\ys(\xs,\ys),
    \end{aligned}
    \right.
\end{equation}
where 
\begin{equation*}
    \Omega(\xs,\ys) = \frac{1}{2}r^2 + \frac{1}{r} - K\xs,
\end{equation*}
and $r=\sqrt{\xs^2+\ys^2}$. In this way, the system \eqref{eq:nousist} has a first integral, known as the Jacobi first integral, defined by
\begin{equation}\label{eq:jacobi}
    C = 2\Omega(\xs,\ys) - (\dot{\xs}^2+\dot{\ys}^2),
\end{equation}
and related to the Hamiltonian \eqref{eq:NumericalH} by $C = -2H$. For a fixed value of $K$ and the first integral $C$, the admissible regions of motion, known as Hill regions, are defined by
\[
\mathcal{R} = \mathcal{R}(K,C) = \left\{(\xs,\ys)\in\R^2 \,|\, 2\Omega(\xs,\ys) \geq C\right\}.
\]

\subsection{Invariant objects}
The equations for equilibrium points $\dot \xs= \dot \ys =\dot \pxs=\dot \pys=0$ in
system~\eqref{eq:NumericalEquations} are just
\[
\pxs=-\ys,\ \pys =\xs, \ \ \xs-\frac{\xs}{r^3}=K,\ \ys-\frac{\ys}{r^3}=0,
\]
which give a whole circle $r=1$ for $K=0$, that is, for the rotating Kepler 
Hamiltonian. For the CP problem $K>0$ there are just two equilibrium points
\[
L_1=(\xs_1,0,0,\xs_1),\qquad L_2=(\xs_2,0,0,\xs_2),
\]
with $\xs_1^3-K\xs_1^2-1 =0$ and $\xs_2^3-K\xs_2^2+1 =0$ so that $\xs_1=\xs_1(K)$,
$\xs_2=\xs_2(K)$ are of the form 
\[
\xs_1(K)=-1+\frac{K}{3}+\Ord{K^2},\qquad \xs_2(K)=\xs_1(-K)=-1+\frac{K}{3}+\Ord{K^2}.
\]

\begin{figure}[!ht]
    \centering
    \includegraphics[trim={0mm 0mm 0mm 0mm},clip,width=0.32\textwidth]{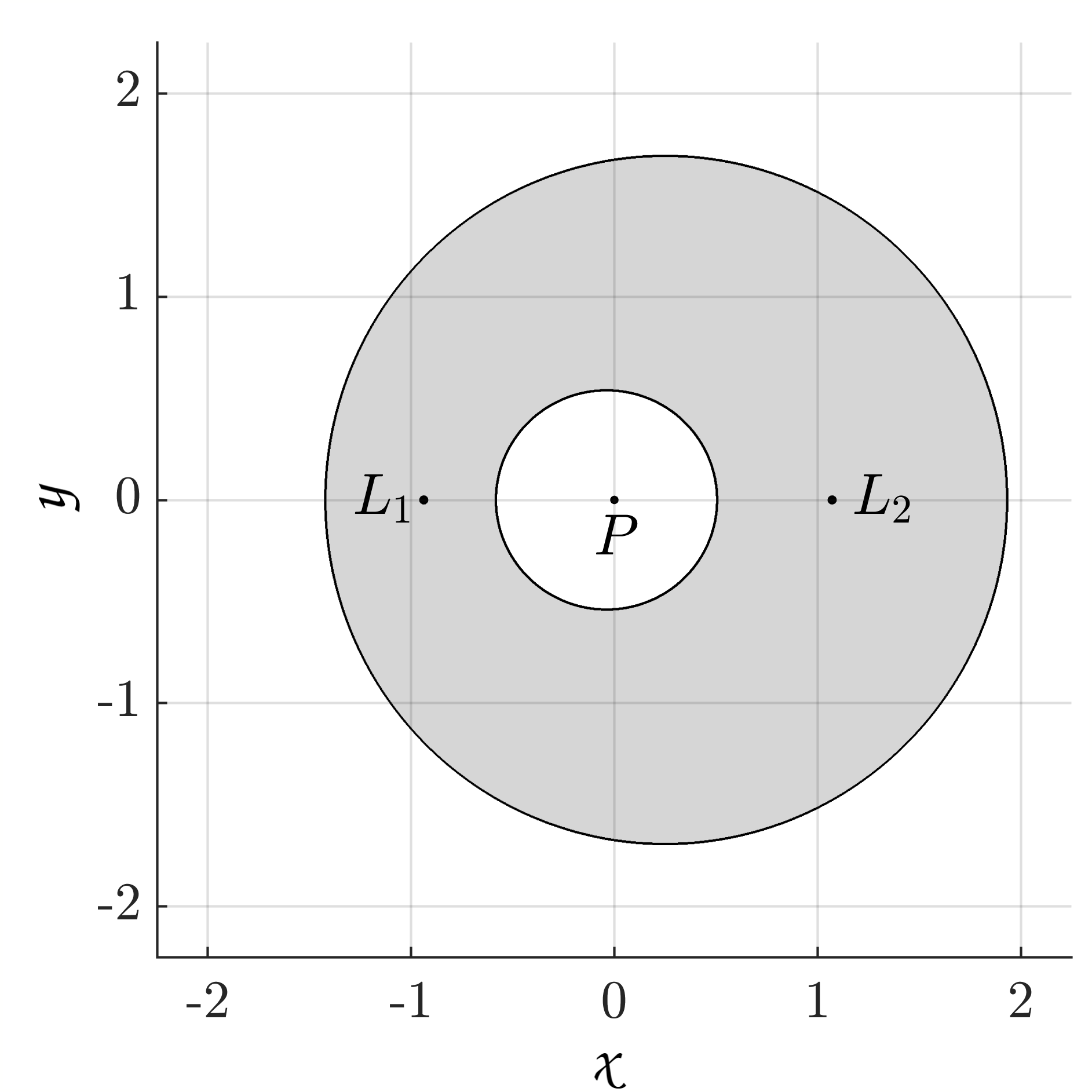} 
    \includegraphics[trim={0mm 0mm 0mm 0mm},clip,width=0.32\textwidth]{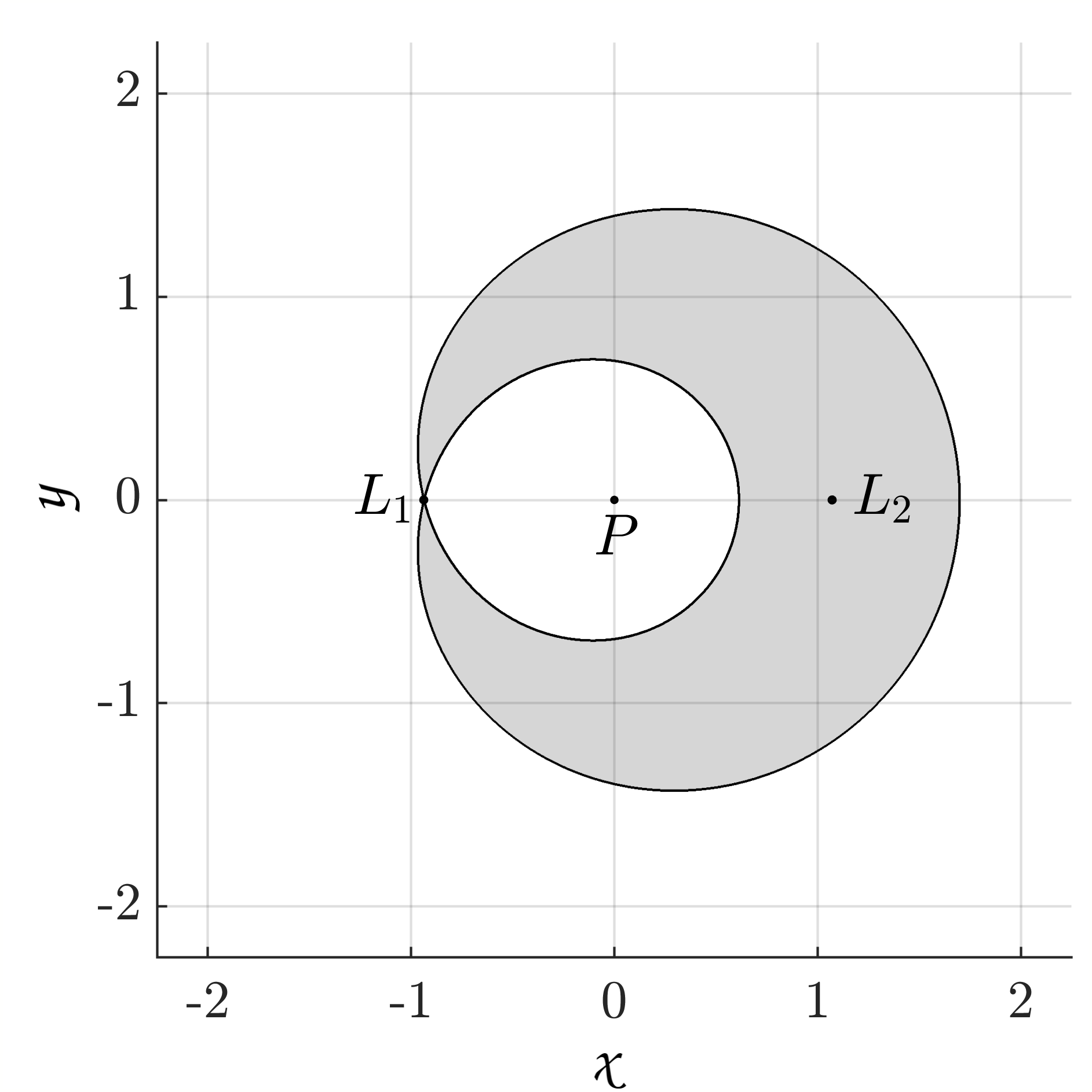} 
    \includegraphics[trim={0mm 0mm 0mm 0mm},clip,width=0.32\textwidth]{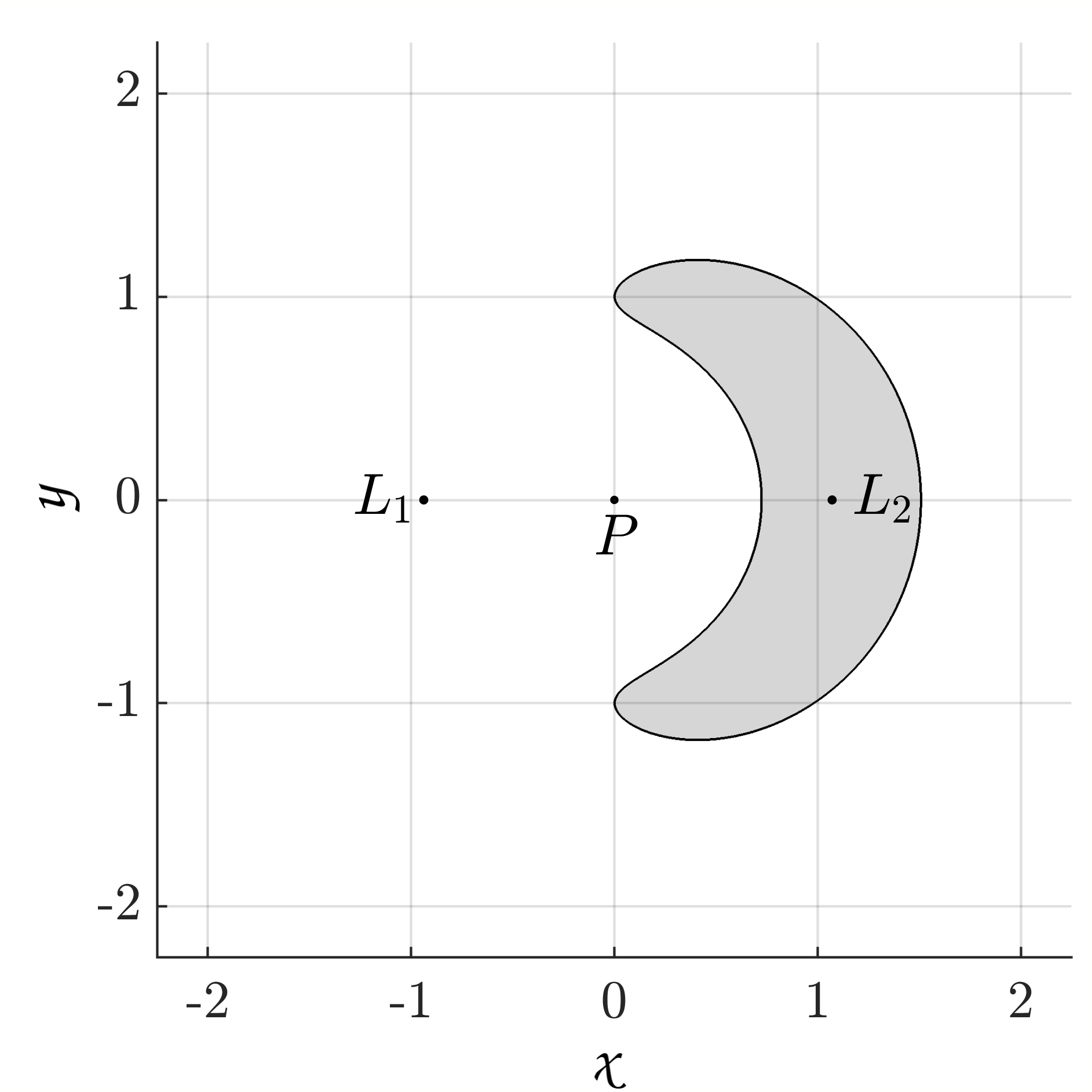} 
    \caption{In gray the forbidden region of motion for $K=0.2$ and $H=-2<h_1$, $H=h_1$ and $H=-1.5\in(h_1,h_2)$ from left to right. $P$ denotes the position of the nucleus.
}
    \label{figHillReg}
\end{figure}

Since $L_i$ are critical points of the potential, the topology of the Hill regions changes at the values of the Jacobi constant corresponding to these points, or equivalently, at the values $h_i = H(L_i)$ of the Hamiltonian \eqref{eq:NumericalH}, for $i=1,2$.

Note that for values $H\leq h_1$, ionization cannot occur, and for values $H\geq h_2$, the entire configuration space becomes a valid region of motion (see Figure~\ref{figHillReg}).
This is of significant importance, as the point $L_1$ and the dynamics associated with the invariant objects emanating from it act as the channel for ionization at the lowest energy levels where it is possible.

In particular, 
$L_1$ is a center-saddle equilibrium point for all $K>0$, with characteristic
exponents
\[
\pm \iu (1+K)(1+\Ord{K}),  \quad \pm \sqrt{3K}(1+\Ord{K}),
\]
and associated 1D invariant manifolds, $\Cu(L_1), \Cs(L_1)$ for
$h_1$. Emerging form $L_1$ there is a family of saddle Lyapunov
periodic orbits $OL$ around $L_1$ for $h>h_1$, with associated
2D-invariant
manifolds $W^{\un}(OL), W^{\st}(OL)$, also called \emph{whiskers}.

Notice that the saddle character of $L_1$ is \emph{weak} with respect to the
parameter $K$, because its saddle characteristic exponents $\pm
\sqrt{3K}(1+\Ord{K})$ tend to $0$ as $K$ tends to $0$. This will imply a very
small distance between $\Cu( L_1)$ and $\Cs( L_1)$
and a very 
small transversality between $W^{\un}(OL)$ and $W^{\st}(OL)$, indeed
exponentially small in $K$. This weak saddle
character of $L_1$ is totally analogous to the weak saddle character of the libration point $L_3$ of the R3BP. 

\begin{figure}[ht!]
    \centering
    \includegraphics[trim={0mm 0mm 0mm 0mm},clip,width=0.32\textwidth]{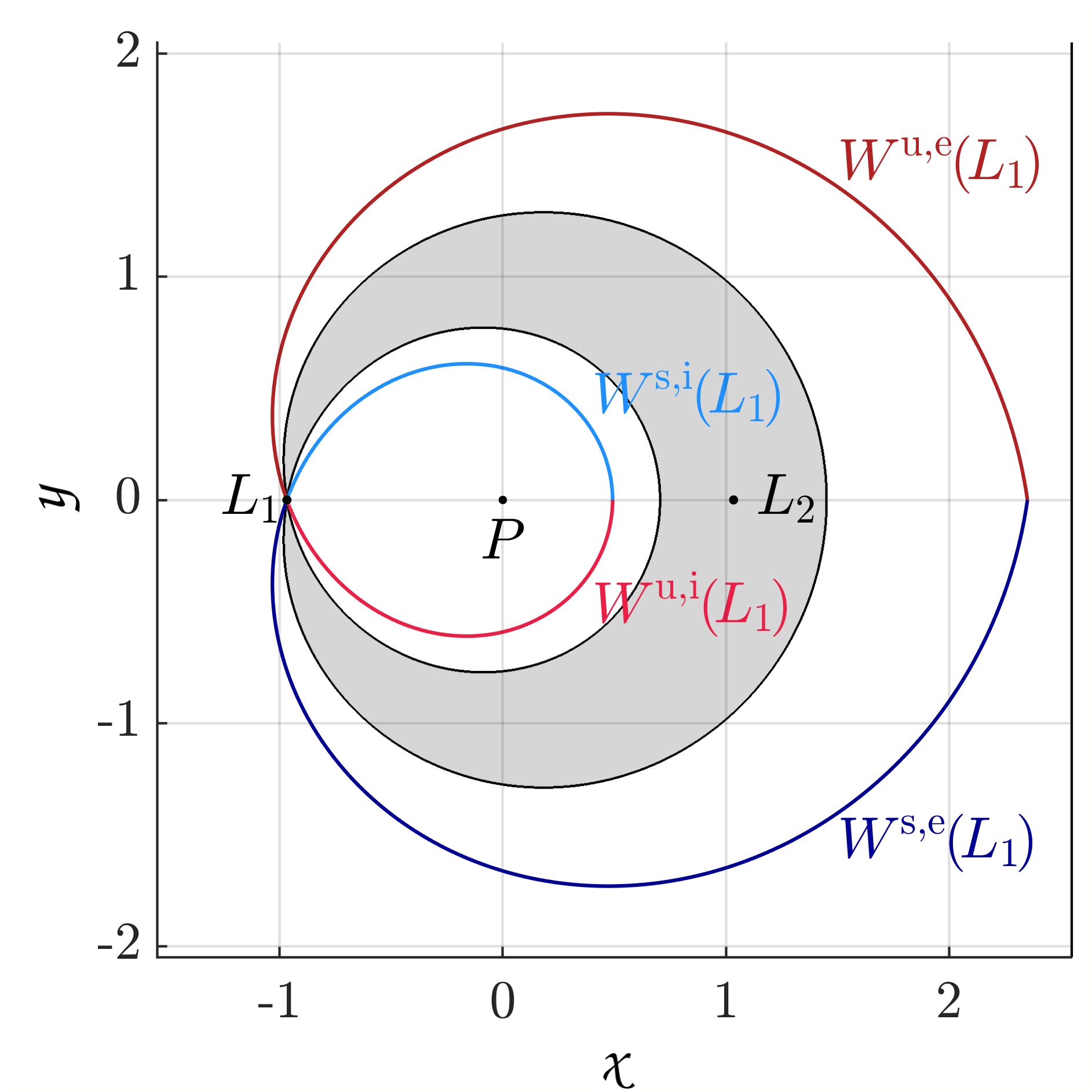}
    \includegraphics[trim={0mm 0mm 0mm 0mm},clip,width=0.32\textwidth]{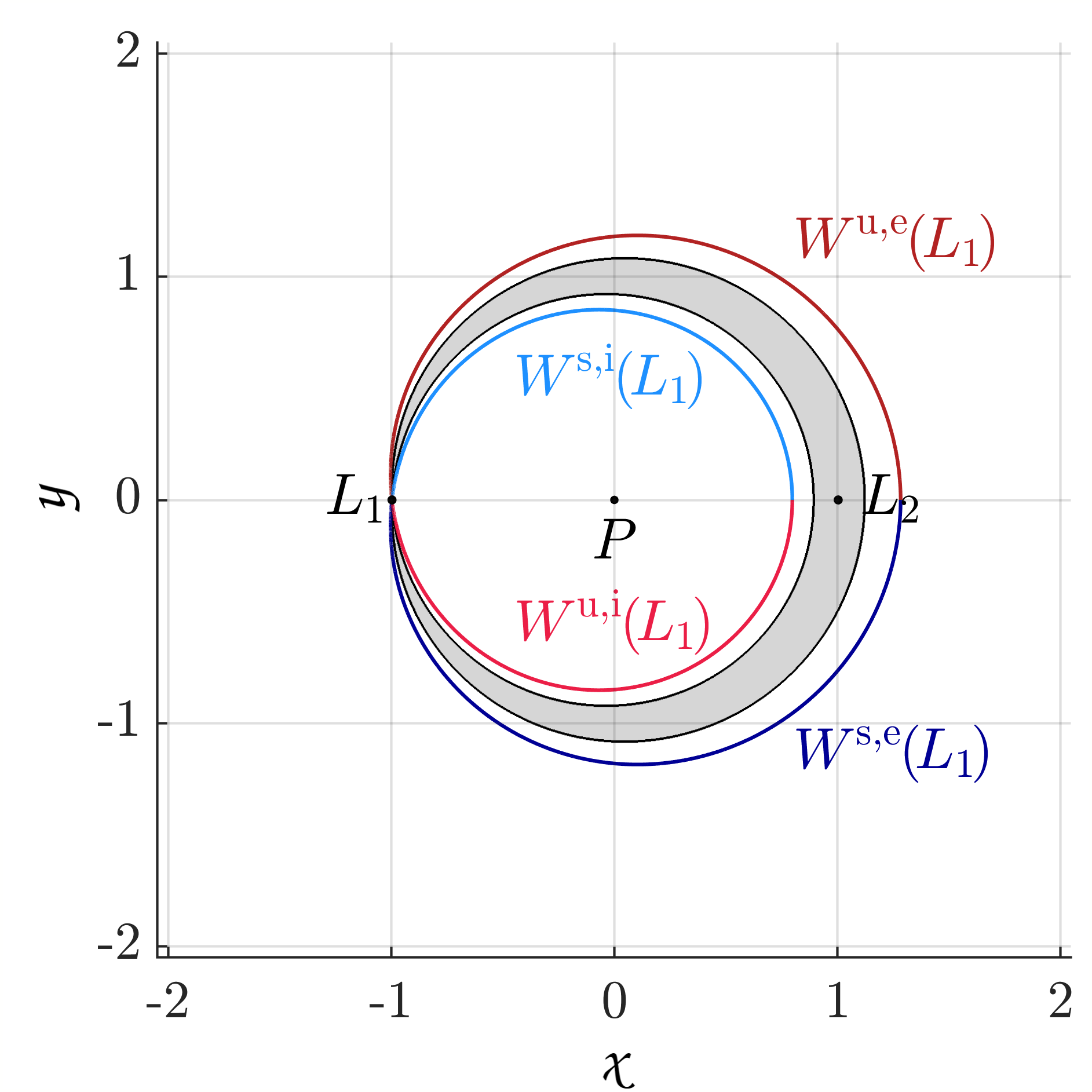}
    \includegraphics[trim={0mm 0mm 0mm 0mm},clip,width=0.32\textwidth]{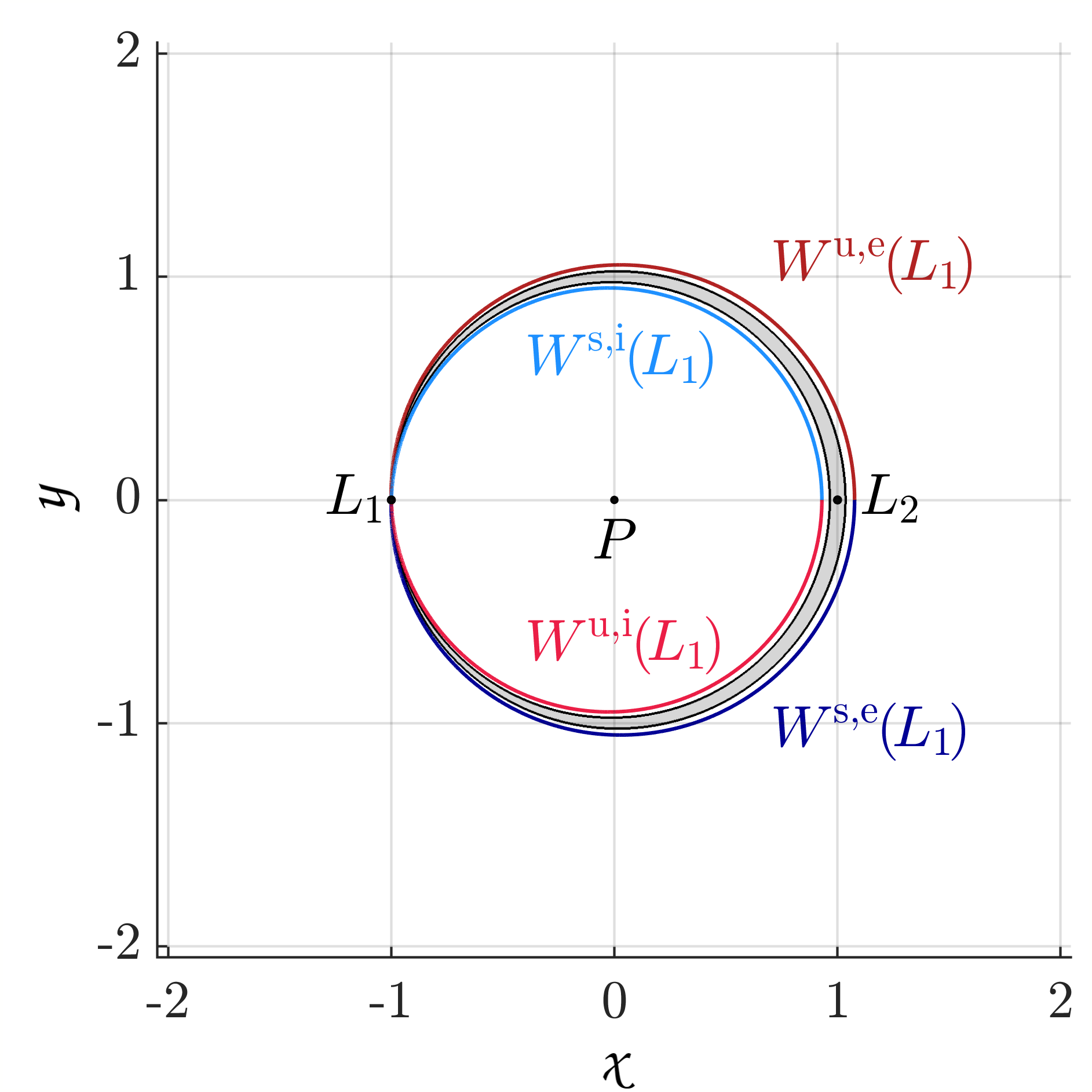}
    \caption{$(\xs,\ys)$ projection of the external $W^{\un,\mathrm{e}}(L_1)$    and internal $W^{\un,\mathrm{i}}(L_1)$ unstable manifolds of $L_1$ (similarly for the stable ones) for $K=0.1$, $0.01$ and $0.001$ from left to right. In gray the forbidden region of motion for $h_1$.
}
    \label{figWUSL1}
\end{figure}

In this paper, we will focus only on the 
invariant manifolds $\Cu(L_1), \Cs(L_1)$.
We distinguish between the external branches of the manifolds,
 $W^{\un,\mathrm{e}}(L_1)$, $W^{\st,\mathrm{e}}(L_1)$
and the internal  branches  $W^{\un,\mathrm{i}}(L_1)$, $W^{\st,\mathrm{i}}(L_1)$. 
The four different branches for $K=0.1$, $0.01$ and $0.001$ are displayed in Figure \ref{figWUSL1}. Notice that the external branches seem to coincide
(similarly the internal ones) but they do not. Precisely, the
purpose of this paper is to analyze the small
distance---splitting---between them. 

Actually we will concentrate on the external branches
that will be simply denoted, from now on, by $\Cu( L_1), \Cs( L_1)$ (unless there are possible misunderstandings).

\subsection{Main Results}\label{subsec14}
\begin{enumerate}
\item
We have numerically found an asymptotic formula for the splitting between the invariant
unstable and stable manifolds (external branches) 
of $L_1$, $W^{\un,\st}(L_1)$, at (the first crossing with)  $\ys=0$ in synodic (rotating)
coordinates $(\xs,\ys,\dot{\xs},\dot{\ys})$~\eqref{eq:NumericalEquations}, which has the expression
\begin{equation}
\label{eq:NumerDistance}
\Delta \dot \xs\sim \varepsilon A\omega^r \exp \left(\frac{\omega\pi}{2}\right),
\qquad \varepsilon=\left(\frac{K}{3}\right)^{1/4}, 
\quad \omega=-\frac{1}{3\varepsilon ^2}=-\frac{1}{\sqrt{3K}},
\end{equation}
with $r=1.61111 \ldots =29/18$ and the leading constant $A=16.055307843\dots$ See \eqref{ajustxp}  and \eqref{Ar} in subsection \ref{subsec3.2}.

\item
We obtain a similar asymptotic formula for the same splitting 
but computed in Poincar\'e (or resonant) variables $(x,y,q,p)$, \eqref{resvar1}
and \eqref{eq:xyqp},
\begin{equation}
\label{eq:ResonantAsympt}
\Delta p\sim \varepsilon \bar A\omega^{\bar r} \exp \left(\frac{\omega\pi}{2}\right),
\end{equation}
see \eqref{deltapu}, with $\bar r=2.11111\ldots=19/9=r+1/2$ and the leading constant $\bar A=27.80860891\dots$ (see  \eqref{Arbar} in subsection \ref{subsec33}), 
where the addition of $1/2$ to $r$ to get $\bar r$ is due to the non-canonical change~\eqref{eq:change} from synodic to resonant variables. See \eqref{relarrbar} in subsection \ref{subsec33}.

\item
We emphasize the ``key point'' of the resonant variables: the analysis of the behavior of the invariant manifolds
$W^{\un,\st}(L_1)$ for the CP problem can be  regarded 
as the behavior of a perturbation of a homoclinic separatrix. Indeed, the Hamiltonian in the resonant variables is
\begin{displaymath}
H(x,y,q,p)=\omega\frac{q^2+p^2}{2}+\frac{y^2}{2}+\cos x-1 
+\varepsilon H_1(x,y,q,p;\varepsilon), 
\end{displaymath}
and the homoclinic separatrix is given by $q=p=0$ and \eqref{eq:StandardSeparatrices}.

\item
We show numerically that the Melnikov method predicts correctly the numerical asymptotic results~\eqref{eq:ResonantAsympt}, as long as we take a \emph{suitable} integrable Hamiltonian system, formed by a rotor and an \emph{amended} pendulum, which is not the natural pendulum (as happens in a lot of other problems, see, for instance,~\cite{Simo1994,GuardiaOS2010,BaldomaM2012,Guardia13}).
 To do so, we consider a \emph{truncated} Hamiltonian once we have applied
a change of variables (from synodic to Poincar\' e ones). More precisely,
\begin{equation}
\label{eq:hamresontruncat}
H(x,y,q,p)=\omega\frac{q^2+p^2}{2}+\frac{y^2}{2}+\cos x-1 - \frac{2}{3}\varepsilon ^2y ^3
+\varepsilon\left( \frac{3q}{2}-\frac{q}{2}\cos 2x-\frac{p}{2}\sin 2x\right), 
\end{equation} 

\item
We introduce a new family of perturbations of a rotor and a pendulum, which we call \emph{Toy CP problem}
\begin{equation}
\label{eq:ToyCP}
H(x,y,q,p)=\omega\frac{q^2+p^2}{2}+\frac{y^2}{2}+\cos x-1 - \frac{2a}{3}\varepsilon ^2y ^3
+\varepsilon^m\left( \frac{3q}{2}-\frac{q}{2}\cos 2x-\frac{p}{2}\sin 2x\right), 
\end{equation}
in order to check if  our results could be provided by previous analytical results in the literature. For $m=1$ Hamiltonian~\eqref{eq:ToyCP} is just the truncated Hamiltonian
~\eqref{eq:hamresontruncat} of the CP problem after a normal form expansion or, equivalently, after an averaging process. It has to be noticed that in the standard application of normal form or averaging, one typically takes $a=0$ because the term $\displaystyle\frac{2a}{3}\varepsilon^2 y^3$ is of higher order than the $O(\varepsilon)$ perturbation. 

Nevertheless, we found out that the numerically computed values of the splitting depend \emph{essentially} on the parameter $a$. Particularly the exponent $r$ in formula~\eqref{eq:NumerDistance}. For the purposes of this paper there are, at least, two different scenarios: the ``standard'' case $a=0$ and the \emph{new} case $a=1$. Concerning the integrable part of the Hamiltonian \eqref{eq:ToyCP}, for $a=0$ we have a standard pendulum, whereas for $a=1$ we have a new pendulum, that we call \emph{amended pendulum}. See the motivation for including this 
value of $a$ in Subsections \ref{subsubsec411}
and \ref{subsec42}.

\item
For $a=0$ there are theoretical results~\cite{Baldoma2006,BaldomaFGS2012} for Hamiltonian systems with $1+\frac12$ degrees of freedom that can be directly applied to the Toy CP problem~\eqref{eq:ToyCP}. According to these theoretical results, for $m\geq 4$ the Melnikov prediction with the standard pendulum gives the correct measure~\eqref{eq:homoL1Anal_compl} for the splitting:
\begin{equation}\label{eq:homoL1Anal_compl_prev}
  \Delta p
  \sim -\varepsilon^m\frac{8\pi}{3}\omega^{3}\left(1-\frac{2}{\omega^{2}}\right)
  \frac{\ee^{\pi\omega/2}}{1-\ee^{2\pi\omega}}.
\end{equation}
We have confirmed numerically these theoretical results in
Subsection \ref{subsubsec411}, but we have further found out that the Melnikov prediction gives also the right result for $m>1$ (see ~Subsection~\ref{41}, which opens the door to improve the quoted theoretical results).

Notice that the exponent $r$ of $\omega$ is $3$ in formula~\eqref{eq:homoL1Anal_compl_prev}. Therefore the above Melnikov prediction with $a=0$ clearly does not give the correct result~\eqref{eq:ResonantAsympt} for the CP problem~\eqref{eq:ResHam}, which takes place just for $m=1$. 

\item
However, for $a=1$, when we have an \emph{amended pendulum} in Hamiltonian~\eqref{eq:ToyCP}, the Melnikov prediction with the amended pendulum coincides with the value~\eqref{eq:ResonantAsympt} of the splitting computed numerically and gives also the right result for $m>1$, both for the exponent $\bar r$ 
 (of $\omega$) and the  leading constant $\bar A$. 
 It is worth noticing that in all formulas the exponential smallness
 $\ee^{\pi\omega/2}$ is the same. This is due to the fact the distance of the closest complex singularity of the separatrices to the real line is the same up to terms of order $1/\omega^4$, see Appendix~\ref{appendix:singularity}.  
 All these computations are detailed in Section \ref{PQP}
  and, in particular, illustrated in Figures 
  \ref{fig:intemel_mval1} and
   \ref{fig:CasCm1a6}. 
   This is, perhaps, the main novelty of this paper: the normal form or averaging procedure to transform a Hamiltonian system close to a resonance to a nearly integrable Hamiltonian must be carried out at a higher order than initially expected, and choosing uniquely some selected term. In passing, it is worth mentioning that the Melnikov integral is computed numerically, since we do not have an explicit expression for the separatrix of the amended pendulum. 
   \item From the numerical simulations done, we can conclude that the Toy CP problem with $a=m=1$ is a simplified model that already describes the splitting $\Delta p$ of the original CP problem in resonant coordinates, as far as asymptotic formulas for $\Delta p$, when $K$ tends to zero, are analyzed.  See Subsection \ref{subsec5p3} and, in particular, 
   Figure \ref{fig:difrealvsCasCm1}.
   Finally, we can also  conclude that the
fitting asymptotic formula obtained for the Melnikov integral provides
the same limit value of the exponent $\bar r=\mathring r=2.111\ldots=19/9$.
\end{enumerate}

\section{Some analytical computations}
\subsection{Polar coordinates}
The Hamiltonian of the CP problem
\[
   H=\frac12(\pxs^2+\pys^2)-\frac1{r}-(\xs\pys-\ys\pxs)+K \xs,\qquad r=\sqrt{\xs^2+\ys^2},
\]
takes the form in \emph{canonical polar coordinates} $(r,\theta,p_r,p_\theta)$
\[
H=\frac12\left(p_r^2+\frac{p_\theta ^2}{2}\right)
	-\frac{1}{r}-p_\theta+Kr\cos\theta :=H_{\mathrm{K}}
	-p_\theta+Kr\cos\theta.
\]
\subsection{Delaunay coordinates}
In the \emph{canonical Delaunay variables} $(\ell,g,L,G)$, which are the
action-angle variables for the Kepler Hamiltonian
$H_{\mathrm{K}}$ we have
\[
H=-\frac{1}{2L^2}-G+K\left(a\cos(\ell+g)-\frac32 a e\cos g
	+\frac{ae}{2}\cos(2\ell +g)+O(e^2)\right),
\]
where $a=L^2$ is the \emph{semi-major axis} and 
$\displaystyle e=\sqrt{1-\frac{G^2}{L^2}}$ the \emph{eccentricity}.

\subsection{Resonant coordinates}\label{subsecRC}
The equilibrium point $L_1$ satisfies $L=G=1$ (for $K=0$).
Changing to resonant variables $(x,y,\varphi,I)$,  and scaling
\begin{equation}\label{resvar1}
\ell=x+\pi+\varphi, \quad g=-\varphi,\quad L=1+\eps^2 y, 
\quad G=1+\eps^2y-\eps^2 I,
\end{equation}
\label{eq:change}
with $\displaystyle \eps=\left(\frac{K}{3}\right)^{1/4}$, as well as averaging (as in \cite{Simo1994,GelfreichS2008,GuardiaOS2010,Guardia13}),  we get the \emph{singular} \emph{resonant} Hamiltonian
\begin{equation}
\label{eq:ResHam}
H(x,y,\varphi,I;\eps)=\omega I +P_{\eps}(x,y)+\eps H_1(x,y,\varphi,I;\eps),
\end{equation}
where $\displaystyle \omega= -\frac{1}{3\eps^2}$ and
\begin{align*}
P_{\eps}(x,y)&=F(y,\eps)+ \cos x -1,\\
H_1(\varphi,x,I,y;\eps)&=\frac32\sqrt{2I}\cos\varphi-\frac12 \sqrt{2I}\cos(2x+\varphi) + O(\eps),
\end{align*}
with
\begin{equation}
\label{eq:F}
F(y,\eps)= \frac{y^2}{2}+\frac{f(\eps^2 y)}{3\eps^4}=\frac{y^2}{2}g(\eps^2 y),
\end{equation}
and
\begin{align*}
f(t)&=\frac{1}{2(1+t)^2}-\frac{1}{2}+t-\frac{3}{2}t^2=-2t^3+\frac{5}{2}t^4+\cdots+(-1)^n \frac{n+1}{2} t^n+\cdots=O(t^3),\\
g(t)&=1+\frac{2f(t)}{3t^2}=1-\frac{4}{3}t+\frac{5}{3}t^2+\cdots+(-1)^m\frac{m+3}{3}t^m+\cdots. 
\end{align*}
Taylor expanding in $\eps$ we have $F(y,\eps)=F_k(y,\eps)+O\left(\eps^{2k+2}y^{k+3}\right)$, where
\[
F_0(y)=\frac{y^2}{2}, \quad F_1(y)=\frac{y^2}{2}-\frac{2}{3}\eps^2y^3, \quad F_2(y)=\frac{y^2}{2}-\frac{2}{3}\eps^2y^3+\frac{5}{6}\eps^4y^4,
\]
and in general
\[
F_k(y)=\sum_{n=2}^{k+2} (-1)^n\frac{n+1}{6}\eps^{2n-4}y^n=y^2\sum_{m=0}^k(-1)^m\frac{m+3}{6}\left(\eps^2 y\right)^m, \quad k\geq 0.
\]

A remarkable property of this Hamiltonian, for our purposes, is that will allow to distinguish an integrable part with a separatrix plus a perturbation and therefore
to apply Melnikov theory. This will be done in the next two subsections.

Further, we introduce  new symplectic coordinates $(p,q)$ related with the action $I$ and the angle $\varphi$ by the change
\begin{equation}\label{eq:xyqp}
   q = \sqrt{2I}\cos{\varphi},\qquad\qquad p = -\sqrt{2I}\sin\varphi. 
\end{equation}
so in the coordinates~\eqref{eq:xyqp}, the resonant Hamiltonian~\eqref{eq:ResHam} casts
\begin{displaymath}
  H(q,p,x,y;\omega,\eps) = F(y) + V(x) + \omega\frac{q^{2} +
  p^{2}}{2} + \eps H_{1}(x,y,q,p;\eps),
\end{displaymath}
 with
\begin{equation}
V(x)=\cos x-1,\qquad 
H_{1}(x,y,q,p; \eps) = \frac{3}{2} q - \frac{q}{2} \cos(2x) - \frac{p}{2} \sin(2x) + O(\eps)   
\label{eq:H1plusO}
\end{equation}
(we keep $H$ and $H_{1}$ to denote, respectively, the Hamiltonian and perturbation term in these new~coordinates, which we also call resonant coordinates). 

{\bf Remark 1}. We notice that the simplest expression ---as a perturbation of
an integrable Hamiltonian---for the Hamiltonian \eqref{eq:HPML1} is
\begin{displaymath}
  H(q,p,x,y;\omega,\eps) = \frac{y^2}{2} + \cos x-1 + \omega\frac{q^{2} +
  p^{2}}{2} + \eps H_{1}(x,y,q,p;0),
\end{displaymath}

\subsection{Standard Melnikov prediction for 
 \texorpdfstring{$L_{1}$}{L1} for 
 \texorpdfstring{$\omega$}{om} independent of 
 \texorpdfstring{$\varepsilon$}{veps}}
Consider a Hamiltonian of type~\eqref{eq:HPML1}, 
\begin{equation}
H=  H(x,y,q,p;\omega,\varepsilon) = F(y) + V(x) + \omega\frac{q^{2} +
  p^{2}}{2} + \varepsilon H_{1}(x,y,q,p;\varepsilon),
  \label{eq:HPML1}
\end{equation}
where along this section we are going to assume that $\omega$ is independent of $\varepsilon$ and we are going to search the so-called ``Melnikov prediction''~\cite{DelshamsG2000,KoltsovaLDG05}.
Let us denote
\begin{equation}
  H_{0}(x,y,q,p) = F(y)+ V(x) + \omega\frac{q^{2} + p^{2}}{2},
  \label{eq:HPM0L1}
\end{equation}
the ``unperturbed'' Hamiltonian, where $V(x) =V(x+2\pi)$ has a non-degenerate global minimum at $x=0$: $V(x)
= -\lambda^{2}x^{2}/2 + O(x^{3})$, $\lambda > 0$ and $\omega\ne 0$, and $F(y)$ has a non-degenerate minimum at $y=0$: $F(0)=F'(0)=0$, $F''(0)=1$.
To distinguish variables it is convenient to introduce
\begin{displaymath}
  u = (x,y),\qquad v = (q,p)
\end{displaymath}
to denote the ``saddle'' coordinates and the
``elliptic'' or ``center'' coordinates in~\eqref{eq:HPML1}, respectively. In
particular, with this notation the perturbation term
can be written as 
\begin{displaymath}
    H_{1}(x,y,q,p;\eps) = H_{1}(u,v;\varepsilon)
\end{displaymath}
We are also going to assume that for the saddle part of
Hamiltonian~\eqref{eq:HPM0L1} 
\begin{displaymath}
  P(x,y) = F(y) + V(x)
\end{displaymath}
there exists a \emph{separatrix}
\begin{equation}
\label{eq:separatrix}
u_{0}(t) = (x_{0}(t), y_{0}(t))
\end{equation}
contained in $P=0$ satisfying
\begin{equation}
    \label{PropertiesSeparatrix}
x_0(0)=0, \quad x_0(t)\to 2\pi \text{ for }t\to\infty,\quad
x_0(t)\to 0 \text{ for }t\to-\infty, \qquad
y_{0}(t)\to 0 \text{ for } t\to\pm\infty.
\end{equation}

For instance, the \emph{standard pendulum} has potential energy is $V(x)=\cos x -1$ and kinetic energy $F(y)=\dfrac{{y}^2}{2}$, which give rise to two \emph{explicit} separatrices
\begin{equation}
\label{eq:StandardSeparatrices}
x_0(t)=x_0^{\pm}(t)=4\arctan\left(\mathrm{e}^{\pm t}\right), \quad y_0(t)=y_0^{\pm}(t)=\dot{x}_0^{\pm}(t)=\pm\frac{2}{\cosh t}.
\end{equation}

On the other hand, for the resonant Hamiltonian~\eqref{eq:ResHam}, $V(x) = \cos x - 1$ is the potential of a standard pendulum,
but the kinetic energy  $F(y)$ given in~\eqref{eq:F} depends essentially on $\varepsilon ^2 y$: 
\begin{displaymath}
F(y)=y^2\left(\frac{1}{2}-\frac{2}{3}\eps^2y+\frac{5}{6}\eps^4y^2+\cdots\right), 
\end{displaymath}
and we have an \emph{amended pendulum},
which coincides with the standard pendulum only for $\eps=0$, and which has also two separatrices.

For $H_1= 0$ the origin $(u,v)=(0,0)$ is a saddle-center equilibrium point with a separatrix $(u_0,v_0)$ satisfying  $P(u_{0}) =
0$, $v_0 = 0$. 
For $0 < \abs{\varepsilon} \ll 1$ there exist
stable and unstable invariant manifolds, $W^{\st,\un} =
\left\{\left(u^{\st,\un}(t), v^{\st,\un}(t)\right)\right\}$, associated to the new
saddle-center $\left(u^{\text{eq}},v^{\text{eq}}\right)$ which is a
critical point of $H$.  Writing 
\begin{displaymath}
  \left(u^{\text{eq}}, v^{\text{eq}}\right) =
  \left(\eps u_{1}^{\text{eq}}, \eps v_{1}^{\text{eq}}\right) +
  \Ord{\eps^{2}},
\end{displaymath}
we get
\begin{equation}
  u_{1}^{\text{eq}} = -\begin{pmatrix}
    -1/\lambda^{2} & 0\\
  0 & 1\end{pmatrix}\nabla_{u}H_{1}(0,0;0),\qquad
  v_{1}^{\text{eq}} = -\frac{1}{\omega}\nabla_{v}H_{1}(0,0;0). 
  \label{eq:u1v1eq}
\end{equation}
To measure the distance between $W^{\un}$ and $W^{\st}$ we can fix the
level of energy of $(u^{\text{eq}}, v^{\text{eq}})$ and choose a
suitable surface of section just to measure $\Delta v(t):= v^{\un}(t) - v^{\st}(t)$.
Notice, on the one hand, that $\left(u^{\st,\un},
v^{\st,\un}\right) = \left(u_{0}(t),0\right)$ for $H_1=0$, so we can write
$v^{\st,\un}(t) = \eps v_{1}^{\st,\un}(t) + \Ord{\eps^{2}}$ and
 we are going to obtain an explicit formula for $\Delta v_1(t)=v_1^u(t)-v_1^s(t)$.
 On the
other hand, observe that the differential equation satisfied, for example, by
$v^{\un}(t)$ is
\begin{displaymath}
   \dot{v}^{\un} = J_{2}\left(\nabla_{v}H_{0}\left(u^{\un},v^{\un}\right)
   +\eps\nabla_{v}H_{1}\left(u^{\un},v^{\un};\eps\right)\right),
\end{displaymath}
with $J_{2} = \begin{pmatrix}
\phantom{-}0 & 1\\ -1 & 0\end{pmatrix}$. Hence,
taking only the first order in $\eps$ yields
\begin{equation}
  \dot{v}_{1}^{\un} = J_{2} \begin{pmatrix}
  \omega & 0\\ 0 & \omega\end{pmatrix}v_{1}^{\un} 
  + J_{2}\nabla_{v}H_{1}(u_{0}(t),0;0) = \omega J_{2} v_{1}^{\un} 
  + J_{2} \nabla_{v}H_{1}(u_{0}(t),0;0).
  \label{eq:EDOv1}
\end{equation}
Besides, we define 
\begin{displaymath}
  b(t) := J_{2} \nabla_{v}H_{1}(u_{0}(t),0;0),
\end{displaymath}
notice that then $b(t) \to b_{\infty} := J_{2}\nabla_{v}H_{1}(0,0;0)$, when
$t\to\pm\infty$. If further, we introduce $\tilde{b}(t) :=
b(t)-b_{\infty}$, the system of linear differential
equations~\eqref{eq:EDOv1} casts 
\begin{displaymath}
  \dot{v}_{1}^{\un} = \omega J_{2} v_{1}^{\un} + b_{\infty} + \tilde{b}(t)
\end{displaymath}
and clearly we can find the equilibrium point of the autonomous part of
this last system by solving $w J_{2} v_{1}^{u_{0}} + b_{\infty}=0$ to get 
\begin{displaymath}
  v_{1}^{\text{eq}} = -\frac{1}{\omega}\nabla_{v}H_{1}(0,0;0)
\end{displaymath}
(as in \eqref{eq:u1v1eq}) and introducing $\tilde{v}_{1}^{\un} = v_{1}^{\un} - v_{1}^{\text{eq}}$
we get
\begin{equation}
\dot{\tilde{v}}_{1}^{\un} = \omega J_{2} \tilde{v}_{1}^{\un} +
\tilde{b}(t),
  \label{eq:v1TildeODE}
\end{equation}
with $\tilde{v}_{1}^{\un}(t)\to 0$ (or $v_{1}^{\un}(t)\to
v_{1}^{\text{eq}}$) as $t\to -\infty$. Next, solving the linear system of
ordinary differential equations~\eqref{eq:v1TildeODE} by variation of
constants one gets
\begin{displaymath}
  \tilde{v}_{1}^{\un}(t) = \ee^{\omega t J_{2}} \tilde{v}_{1}^{\un}(0)
  + \int_{0}^{t} \ee^{\omega (t-s) J_{2}} \tilde{b}(s)\df s
\end{displaymath}
where we stress that $\ee^{\omega t J_{2}} = \begin{pmatrix}
  \phantom{-}\cos\omega t & \sin\omega t\\
-\sin\omega t & \cos\omega t\end{pmatrix} = R(-\omega t)$ is a rotation of
angle $-\omega t$. Now, rewriting this solution as
\begin{displaymath}
  R(\omega t) \tilde{v}_{1}^{\un}(t) = \tilde{v}_{1}^{\un}(0)
  +\int_{0}^{t} R(\omega s)\,\tilde{b}(s)\df s
\end{displaymath}
and taking limits $t\to -\infty$ at both sides, the rhs clearly goes to
zero, so that $\displaystyle\tilde{v}_{1}^{\un}(0) = \int_{-\infty}^{0} 
R(\omega s)\tilde{b}(s)\df s$, and therefore  
\begin{displaymath}
  \tilde{v}_{1}^{\un}(t) = R(-\omega t) \int_{-\infty}^{0} R(\omega
  s)\tilde{b}(s)\df s + \int_{0}^{t}
  R\left(\omega(s-t)\right)\tilde{b}(s)\df s 
  = \int_{-\infty}^{t} R\left(\omega (s-t)\right)\tilde{b}(s)\df s.
\end{displaymath}
Hence $\displaystyle v_{1}^{\un}(t) = v_{1}^{\text{eq}} +
\int_{-\infty}^{t} R\left(\omega(s-t)\right)\tilde{b}(s)\df s$, or more
explicitly,
\begin{equation}
  v_{1}^{\un}(t) = -\frac{1}{\omega}\nabla_{v} H_{1}(0,0;0)
  + \int_{-\infty}^{t} R\left(\omega(s-t)\right)
  J_{2}\left(\nabla_{v}H_{1}(u_{0}(s),0;0)
  -\nabla_{v} H_{1}(0,0;0)\right)\df s
  \label{eq:intv1u}
\end{equation}
If should one have started from the differential equation satisfied by
$v^{\st}$, i.e.,
\begin{displaymath}
  \dot{v}_{1}^{\st} = J_{2}\left(\nabla_{v}H_{0}(u^{\st}, v^{\st})
  + \eps\nabla_{v}H_{1}(u^{\st},v^{\st};\eps)\right),
\end{displaymath}
having kept, as in~\eqref{eq:EDOv1}, only the first order in $\eps$, and
proceeded the same way we had done to get~\eqref{eq:intv1u}, one would have
arrived to the formula for $v_{1}^{\st}$, which turns out to be
\begin{equation}
  v_{1}^{\st}(t) =
  -\frac{1}{\omega}\nabla_{v}H_{1}(0,0;0) 
  -\int_{t}^{\infty} R\left(\omega(s-t)\right) J_{2}
  \left(\nabla_{v}H_{1}(u_{0}(s),0;0) - \nabla_{v}H_{1}(0,0;0)\right)\df s,
  \label{eq:intv1s}
\end{equation}
so that, subtracting~\eqref{eq:intv1s} to~\eqref{eq:intv1u}, yields

\begin{equation}\label{eq:v1uv1s}
  \Delta v_1(t):= v_{1}^{\un}(t) - v_{1}^{\st}(t) =
  \int_{-\infty}^{\infty} R\left(\omega(s-t)\right) J_{2}
  \left(\nabla_{v} H_{1}(u_{0}(s),0;0)-\nabla_{v}H_{1}(0,0;0)\right)\df s
  =R(-\omega t)\Delta v^{0}_{1},
\end{equation}
where we have introduced
\begin{displaymath}
  \Delta v^{0}_{1} := J_{2}
  \int_{-\infty}^{\infty} R(\omega s)
  \left(\nabla_{v} H_{1}(u_{0}(s),0;0)-\nabla_{v}H_{1}(0,0;0)\right)\df s.
\end{displaymath}
Notice that $\lVert\Delta v(t)\rVert=\lVert  v_{1}^{\un}(t) -
v_{1}^{\st}(t)\rVert=\lVert v^{0}_{1}\rVert$.

Let us compute this \emph{constant} separation between invariant manifolds for the
amended pendulum, that is, for $V(x)=\cos x -1$, and for a given separatrix
$u_{0}(t) = \left(x_{0}(t),y_{0}(t)\right)$
satisfying~\eqref{PropertiesSeparatrix}.

Thanks to the even character of the potential $V(x)$, we have the invariance of
the amended pendulum under the reversibility $(t,x,y)\mapsto (-t,2\pi-x,y)$
then, as for small $\varepsilon$ the points $\left(\pi, \pm 2 +
O\left(\varepsilon^{2}\right)\right)$ belong to the separatrix, we can take a
parameterisation of the separatrix~\eqref{eq:separatrix}
such that~$\left(x^{\pm}_{0}(0), y^{\pm}_{0}(0)\right) = \left(\pi, \pm 2 +
O\left(\varepsilon^{2}\right)\right)$, and hence satisfies
\begin{equation}
    (x_{0}(t),y_{0}(t))=(2\pi-x_{0}(-t),y_{0}(-t))  
    \label{eq:reversibility-of-the-separatrix}
\end{equation}
which, in particular, implies that $\sin\left(x_{0}(t)\right)$, $\sin\left(2
x_{0}(t)\right)$ are odd in $t$, and $\cos\left(x_{0}(t)\right)$, $\cos\left(2
x_{0}(t)\right)$ are even in $t$. Note that, of course, this
is true also for the explicit separatrix of the standard pendulum ($\varepsilon
= 0)$ given in~\eqref{eq:StandardSeparatrices}.

In particular, if we consider the perturbation term $H_{1}(q,p,x,y;\eps)$ given by~\eqref{eq:H1plusO}
\begin{equation}
  H_{1}(x,y,q,p;0) = \frac{3}{2}q - \frac{q}{2}\cos 2x -
  \frac{p}{2}\sin 2x,
  \label{eq:H1EPS0}
\end{equation}
one easily gets that

\begin{displaymath}
  R(\omega s)\left(\nabla_{v}H_{1}\left(u_{0}(s),0;0\right) -
  \nabla_{v}H_{1}\left(0,0;0\right)\right) 
  =\frac{1}{2}\begin{pmatrix}
      \cos(\omega s) - \cos\left(2x_{0}(s)+\omega s\right)\\
      \sin(\omega s) - \sin\left(2x_{0}(s)+\omega s\right)
    \end{pmatrix}.
\end{displaymath}     
Introducing
\begin{equation}
\mathcal{A}(\omega):=\int_{-\infty}^{\infty}\ee^{\iu\omega s}
  \left(1 - \ee^{2\iu x_{0}(s)}\right)\!\df s,
  \label{eq:Aomega-complex}
\end{equation}
we just get
\begin{displaymath}
  \Delta v^{0}_{1} = \frac{1}{2}J_{2}\begin{pmatrix}
    \text{Re} \mathcal{A}(\omega)\\
  \text{Im} \mathcal{A}(\omega)\end{pmatrix} 
  = \frac{1}{2}\begin{pmatrix}
     \text{Im} \mathcal{A}(\omega)\\- \text{Re} \mathcal{A}(\omega)\end{pmatrix} = 
     \frac{1}{2}
     \renewcommand{\arraystretch}{2.5}
     \begin{pmatrix}
     \displaystyle
     \int_{-\infty}^{\infty} \left(\sin\left(\omega s\right) - 
            \sin\left(2 x_{0}(s) + \omega s\right)\right)\! \df s\\
     \displaystyle
     -\int_{-\infty}^{\infty} \left(\cos\left(\omega s\right) - 
            \cos\left(2 x_{0}(s) + \omega s\right)\right)\! \df s\\
     \end{pmatrix}.
\end{displaymath}
and we define $\mathcal{A}^{\pm}(\omega)$ depending on the component $x_0^{\pm}(s)$ of the
separatrix chosen. 

Taking into account the reversion
property~\eqref{eq:reversibility-of-the-separatrix} of the parameterisation 
we use for the separatrix, it is straightforward to see that $\text{Im}
\mathcal{A}(\omega)=0$, so that $\text{Re} \mathcal{A}(\omega)=\mathcal{A}(\omega)$; that is,
\begin{equation}\label{FAw}
\mathcal{A}(\omega):=\int_{-\infty}^{\infty}\left(\cos(\omega s)
  -\cos\left({2 x_{0}(s)+\omega s}\right)\right)\!\df s,
\end{equation}
and the difference~\eqref{eq:v1uv1s} is just given by
\begin{equation}
  \Delta v_{1}(t)=v^{\un}_{1}(t) - v^{\st}_{1}(t)
  = -\frac{1}{2} R(-\omega t)\begin{pmatrix}
     0\\ \mathcal{A}(\omega)\end{pmatrix},\label{eq:homoL1Anal}
\end{equation}
that is,
\begin{equation}\label{eq:deltapmel}
  \Delta q_{1}(t):=q^{\un}_{1}(t) - q^{\st}_{1}(t) = -\frac{1}{2}\mathcal{A}(\omega)\sin\omega t,\\ 
\qquad  \Delta p_{1}(t):=p^{\un}_{1}(t) - p^{\st}_{1}(t) = -\frac{1}{2}\mathcal{A}(\omega)\cos\omega t.  
\end{equation}

Notice that for $t=0$ the linear approximation of $\Delta p(0)$ is
\begin{equation}\label{deltap}
 \varepsilon \Delta p_1^{\pm}(0), \qquad \text{with} \qquad      \Delta p_{1}^{\pm}(0)=-\frac{\mathcal{A}^{\pm}(\omega)}{2}.
\end{equation}
and $\Delta q_1^{\pm}(0)=0$.

We finish this section to exhibit some explicit computations that can be
obtained for the standard pendulum, that is for $F(y)=\dfrac{{y}^2}{2}$, where
the separatrices have the  closed-form expression given in
\eqref{eq:StandardSeparatrices}, so that $\mathcal{A}^{\pm}(\omega)$ is given by,
\begin{align}
  \mathcal{A}^{\pm}(\omega):=\int_{-\infty}^{\infty}
  \left(\cos(\omega s)
  -\cos\left({2 x^{\pm}_{0}(s)+\omega s}\right)\right)\! \df s
  &= 
  \frac{4}{3}\omega^{3}\left(1-\frac{2}{\omega^{2}}\right)
  \left[\frac{1}{\cosh\left(\dfrac{\pi\omega}{2}\right)}\mp
  \frac{1}{\sinh\left(\dfrac{\pi\omega}{2}\right)}\right]\nonumber\\
  &= \frac{16\pi}{3}\omega^{3}\left(1-\frac{2}{\omega^{2}}\right)
  \frac{\ee^{c^{\pm}\pi\omega/2}}{1-\ee^{2\pi\omega}},
  \label{eq:Aomega}
\end{align}
 being 
\begin{displaymath}
  c^{\pm} = \begin{cases}
    1, & \text{ for $u^{+}(t)$, the external branch $y_{0}(t)>0$, }\\
    3, & \text{ for $u^{-}(t)$, the internal branch $y_{0}(t)<0$. }
    \end{cases}
\end{displaymath}
See Appendix~\ref{appendix:intAw}  for details.

\subsection{Averaging and Melnikov method}
The method developed in the previous section, based on using the variational equations associated with a separatrix to measure the splitting that occurs when a perturbation is added, is commonly known as the \emph{Melnikov method}, or the \emph{Poincaré-Melnikov-Arnold method} in the Hamiltonian case (see \cite{Guck,DelshamsG2000,KoltsovaLDG05}). In few words, we approximate the system of equations associated to Hamiltonian~\eqref{eq:HPML1}
\begin{equation}
\label{eq:HamSys}
\begin{split}
\dot x&= F'(y)+\eps\frac{\partial H_1}{\partial y}(x,y,q,p,\eps),\\    
\dot y&= -V'(x)-\eps\frac{\partial H_1}{\partial x}(x,y,q,p,\eps),\\ 
\dot q&= \omega p+\eps\frac{\partial H_1}{\partial p}(x,y,q,p,\eps),\\    
\dot p&= -\omega q-\eps\frac{\partial H_1}{\partial q}(x,y,q,p,\eps),    
\end{split}
\end{equation}
by the simplified system
\begin{equation}
\label{eq:VarHamSys}
\begin{split}
\dot x_0&= F'(y_0),\\    
\dot y_0&= -V'(x_0),\\ 
\dot q_1&= \omega p_1+\eps\frac{\partial H_1}{\partial p}(x_0,y_0,q_1,p_1,0),\\    
\dot p_1&= -\omega q_1-\eps\frac{\partial H_1}{\partial q}(x_0,y_0,q_1,p_1,0),
\end{split}
\end{equation}
which is just the linear variational equation~\eqref{eq:EDOv1} satisfied by $v_1^{\st,\un}=(q_1^{\st,\un},p_1^{\st,\un})$ along the separatrix~\eqref{eq:separatrix}.
In the previous section we imposed that 
$v^{\st,\un}(t)\to 0$ for $t\to\pm\infty$ and we found formula~\eqref{eq:v1uv1s} for the separation
$\Delta v_1(t):= v_{1}^{\un}(t) - v_{1}^{\st}(t)$ between separatrices. It is worth remarking that the separatrix~\eqref{eq:separatrix} satisfies the first two equations of system~\eqref{eq:VarHamSys}, where $F(y)$ may have a whole expansion in the variable $\eps$, as happens in equation~\eqref{eq:F}.

To find up to what order in $\eps$ we need to expand $F(y)$, it is useful to perform some previous averaging or normalization steps of system~\eqref{eq:HamSys}.
Indeed, looking at the particular solutions of system~\eqref{eq:HamSys} for $\eps=0$ and $x=y=0$, which are simply
\begin{displaymath}
v_0(t)=(q_0(t),p_0(t))=\left(\sqrt{2I}\cos(\omega t +\varphi),-\sqrt{2I}\sin(\omega t +\varphi)\right),
\end{displaymath}
where 
\begin{displaymath}
  q_0(0)=\sqrt{2I}\cos\varphi,\qquad p_0(0)=-\sqrt{2I}\sin\varphi,\qquad
  I=\frac{q_0(t)^2+p_0(t)^2}{2},  
\end{displaymath}
we notice that the variables $(q_0,p_0)$ move like an harmonic oscillator, and
one can then think about averaging system~\eqref{eq:HamSys}. For instance, for
the perturbation $H_1$ given in equation~\eqref{eq:H1EPS0}, we observe that when
we substitute $(q,p)$ by $(q_0(t),p_0(t))$, then $H_1(x,y,q_0(t),p_0(t),0)$ has
zero average with respect to $t$, and the same happens to its derivatives with
respect to $x$ and $y$. Therefore, averaging system~\eqref{eq:HamSys} we get
system~\eqref{eq:VarHamSys}. 

To find the error of this averaging approximation, we can try to find the  change of variables $\Phi$ from the Hamiltonian $H=H_0+\eps H_1$, where $H_0=\omega I + P$ and $P(x,y)=F(y)+V(x)$, to the averaged Hamiltonian as the 1-time flow of a Hamiltonian $\eps W$:
\begin{align*}
H\circ \Phi&=H+\{H,\eps W\}+O(\eps^2 W^2)=
H_0+\eps H_1 +\{H_0,\eps W\}+O(\eps^2 W)\\
&=H_0+\eps H_1 +\{\omega I,\eps W\}+\{P,\eps W\} +O(\eps^2 W)=
H_0+\{P,\eps W\} +O(\eps^2 W),
\end{align*}
as long as we solve the so-called \emph{cohomological equation}
\[
\{\omega I,W\}+H_1=0.
\]
For the zero-average perturbation 
  $\displaystyle H_{1}(x,y,q,p;0) = \frac{3}{2}q - \frac{q}{2}\cos 2x -
  \frac{p}{2}\sin 2x$ given in equation~\eqref{eq:H1EPS0}, it is easy to check that
  \begin{equation}
 \label{eq:W} 
W(x,y,q,p) = \left(-\frac{3}{2}p + \frac{p}{2}\cos 2x -
  \frac{q}{2}\sin 2x\right)\frac{1}{\omega}
\end{equation}
solves the cohomological equation. Notice that $W=O(1/\omega)$, so that
\[
H\circ \Phi=H_0+\eps\{P, W\} +O(\eps^2/\omega).
\]

For $\omega=O(1/\eps^2)$ and $H=H_{0}+\eps H_{1}$, then $W=O(\eps^2)$ and
\begin{displaymath}
    H\circ \Phi=H_0+\eps\{P, W\} +O(\eps^4)=H_0+O(\eps^3),    
\end{displaymath}
so that the approximation
\begin{displaymath}
    F(y)=F_1(y)=\frac{y^2}{2}-\frac{2}{3}\eps^2y^3, 
\end{displaymath}
that is, the amended pendulum, is adequate.

As a last important remark, one can see that the perturbation $\eps\{P, W\}$ can be replaced just by $-\eps H_1$ in the computation of all the Melnikov integrals, as they are evaluated on the separatrix $u_0(t)$.

Indeed, we can write the first-order perturbation $H_1$ given in~\eqref{eq:H1EPS0} and $W$ given in~\eqref{eq:W} as
\[
H_1=q+ \widetilde{H_1}, \qquad W=-\frac{p}{\omega}+\widetilde{W},
\]
where $\widetilde{H_1}$ and $\widetilde{W}$ given by
\begin{displaymath}
 \widetilde{H}_1(x,y,q,p;0) = \frac{q}{2}(1-\cos 2x) -
  \frac{p}{2}\sin 2x,\quad
\widetilde{W}(x,y,q,p) = \left(-\frac{p}{2}(1-\cos 2x) -
  \frac{q}{2}\sin 2x\right)\frac{1}{\omega}
\end{displaymath}
satisfy 
\[
\lim_{t\to\pm\infty}\widetilde{H}_1(x_0(t),y_0(t),q,p;0)=0,\qquad
\lim_{t\to\pm\infty}\widetilde{W}(x_0(t),y_0(t),q,p)=0.
\]
Denoting by
$\Phi_t^{0}$ the flow of $H_0$ on the separatrix $\left(x_0(t),y_0(t),q_0(t),p_0(t)\right)$
whereas $\Phi_t^{00}$ denotes the flow of $H_0$ restricted to $x=y=0$, that is $\left(0,0,q_0(t),p_0(t)\right)$, from
\begin{displaymath}
    \frac{\mathrm{d}}{\mathrm{d}t}\left(\widetilde{W}\circ \Phi_t^0\right)=\{\widetilde{W}, H_0\}\circ \Phi_t^0
\end{displaymath}
one finally gets
\[
\int_{-\infty}^{\infty} \{P, W\}\circ \Phi_t^0 \,\mathrm{d}t=
\int_{-\infty}^{\infty} \{P, \widetilde{W}\}\circ \Phi_t^0 \,\mathrm{d}t=
\int_{-\infty}^{\infty} H_1\circ \Phi_t^0 - H_1\circ \Phi_t^{00}\,\mathrm{d}t.
\]
Similarly one can see that
\[
\int_{-\infty}^{\infty} R(\omega t)\,\{P, W\}\circ \Phi_t^0 \,\mathrm{d}t=
\int_{-\infty}^{\infty} R(\omega t)\left(H_1\circ \Phi_t^0 - H_1\circ \Phi_t^{00}\right)\!\df t,
\]
and an analogous result holds for $\displaystyle \frac{\partial W}{\partial q}$ and $\displaystyle \frac{\partial W}{\partial p}$  instead of $W$, for instance
\[
\int_{-\infty}^{\infty} R(\omega t)\,\nabla_v\{P, W\}\circ \Phi_t^0 \,\mathrm{d}t=
\int_{-\infty}^{\infty} R(\omega t)\left(\nabla_v H_1\circ \Phi_t^0 - \nabla_v H_1\circ \Phi_t^{00}\right)\! \df t.
\]
%
\section{Numerical results. 
 Splitting of the invariant manifolds of 
 \texorpdfstring{$L_{1}$}{L1}, 
 \texorpdfstring{$W^{\un,\st}(L_1)$}{WSL1}}
In this part of the paper we discuss the methodology used and  the results
obtained from numerical computations.
   
\subsection{Numerical computation of the manifolds 
\texorpdfstring{$W^{\un,\st}(L_1)$}{WUSL1} of the equilibrium point 
\texorpdfstring{$L_{1}$}{L1}}   
Our first goal is to compute $W^{u,s}(L_1)$ numerically. An effective way to so
so is applying the parameterization method (shortly reminded in Appendix A) in
order to obtain high order expansions of $W^{u,s}(L_1)$. Moreover, since we will
need to deal with very small values of $K$ to check exponentially small
estimates for the splitting between $W^{u}(L_1)$ and $W^{s}(L_1)$, a multiple
precision arithmetic will be required. Also regarding the numerical
computations, we observe that using the Hamiltonian for the CP problem in
synodical coordinates $(\xs,\ys,\dot \xs,\dot \ys)$ implies a term $1/r^3$ in
the ODE. To make the computations more efficient, we will change to Levi-Civita
coordinates $(u,v,u',v')$ and a new time $\tau$ (with $'=\df/\df\tau$) defined
by
\begin{equation}\label{LC}
\xs=u^2-v^2,\qquad \ys=2uv,\qquad \frac{\df t}{\df\tau}=4(u^2+v^2)
\end{equation}
so the system of ODE now becomes

\begin{displaymath}
    \left\{
    \begin{aligned}
        & u'' - 8(u^2+v^2)v' =  -4Cu -16Ku^3 + 12(u^2+v^2)^2u,\\
        & v'' + 8(u^2+v^2)u'  = -4Cv +16Kv^3 + 12(u^2+v^2)^2v,
    \end{aligned}
    \right.
\end{displaymath}
where $C$ is given by \eqref{eq:jacobi} and $'=\df/\df \tau$, which is simply polynomial, and the
implementation of the parameterization method becomes simpler and more efficient.

We recall that:

\noindent 1. The Levi-Civita tranformation \eqref{LC} duplicates the configuration plane.

\noindent 2. There exists a first integral expressed by

$$
u'^2+v'^2=8(u^2+v^2)\left(
\frac{1}{2}(u^2+v^2)^2+\frac{1}{u^2+v^2}-K(u^2-v^2)-\frac{C}{2}
\right)
$$
which is regular everywhere (including the collision $u=v=0$).

\subsection{Splitting in synodic coordinates}\label{subsec3.2}
The next step consists of computing the distance between the unstable and
stable invariant manifolds associated with $L_1$ at some Poincar\'e section $\Sigma$, which is taken as
$\ys=0$, $\xs>0$ (or correspondingly to $v=0$, $u>0$ in Levi-Civita coordinates).
We will focus on the external manifolds $W^{\un,\mathrm{e}}$ and
$W^{\st,\mathrm{e}}$ (see Figure~\ref{figWUSL1}), and from now on we simply
denote them by $\Cu$ and $\Cs$. Fixed a value of $K>0$, we want to compute
the distance --also called the splitting-- between $W^{u}$ and $W^{s}$ at the
first crossing with $\Sigma$, that is the distance between the two points
$P^u=(\xs^u,\ys^u,\pxs^u,\pys^u)$ and $P^s=(\xs^s,\ys^s,\pxs^s,\pys^s)$, with $\ys^u=\ys^s=0$
(we will provide the splitting for the variables $(\xs,\ys,\pxs,\pys)$ although the
numerical simulations will be done using Levi-Civita coordinates). However
notice that due to the symmetry \eqref{eq:reversibility}, $\xs^u=\xs^s$,
$\pys^u=\pys^s$ and $\pxs^s=-\pxs^u$, the stable manifold $\Cs$ does not need to be
computed and the distance between $P^u$ and $P^s$ is simply
$2\pxs^u=2\dot \xs^u:=\Delta \dot \xs$, since $\pxs=\dot \xs-\ys=\dot \xs$ at $\Sigma$ (of
course the corresponding  reversibility applies in Levi-Civita coordinates).

Next, for different values of $K$ we compute the corresponding value
$\Delta\dot \xs$. We have taken a range of values of $K$ decreasing from
$0.001$ to $10^{-8}$. 
 
\begin{figure}[ht!]
    \centering
    \includegraphics[trim={0mm 0mm 0mm 0mm},clip,width=0.48\textwidth]{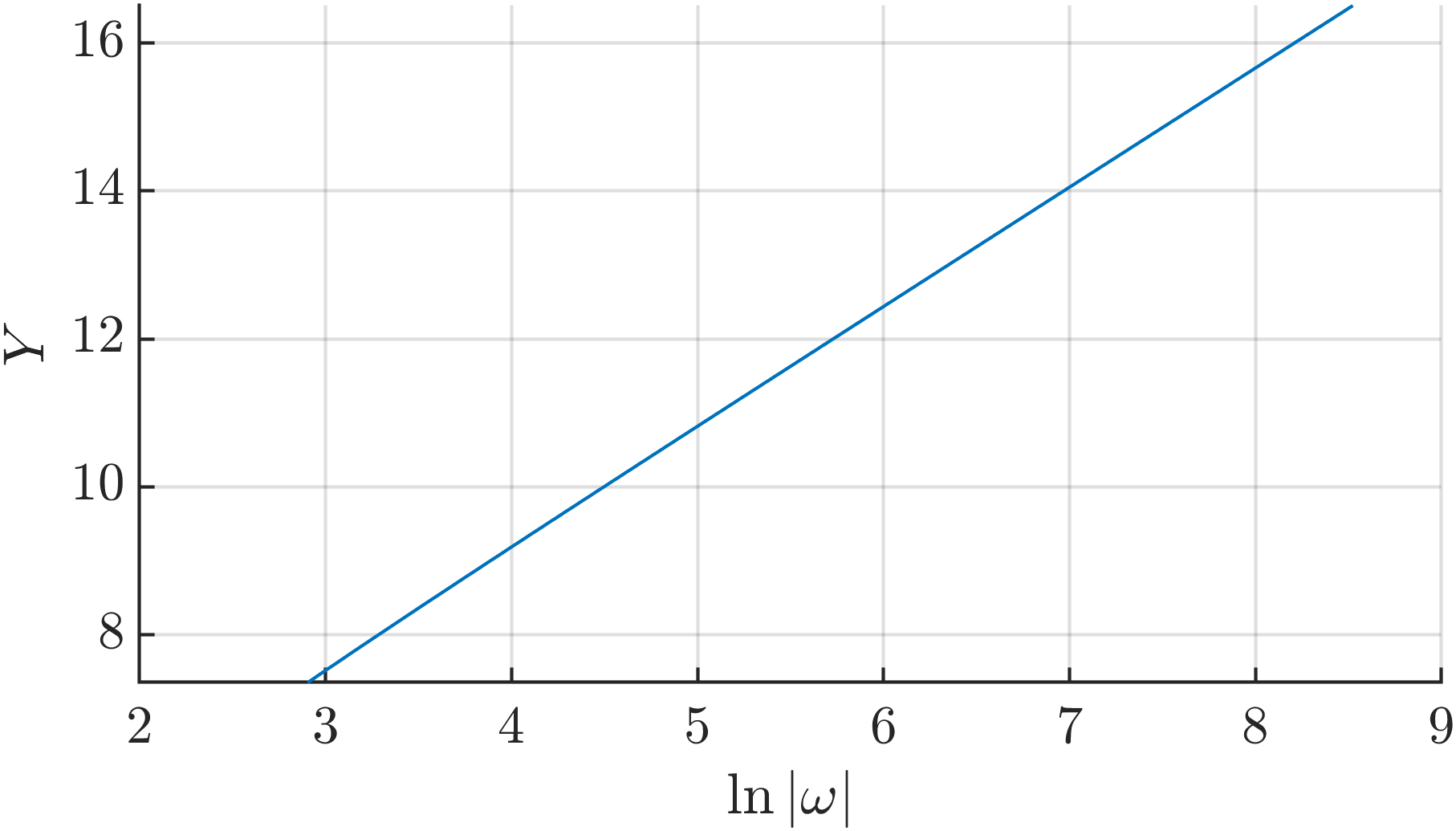}
    \caption{Points $(\ln |\omega|,Y_{\Delta\dot \xs})$ varying 
    (decreasing to zero) $K$ (or equivalently increasing $\ln |\omega|$). Notice the points {\sl lying} on the line 
    $Y=2.772+1.6115\ln |\omega|$. 
    }\label{fig:logw_Y_xp}
\end{figure}

\begin{figure}[ht!]
    \centering
    \includegraphics[trim={0mm 0mm 0mm 0mm},clip,width=0.48\textwidth]{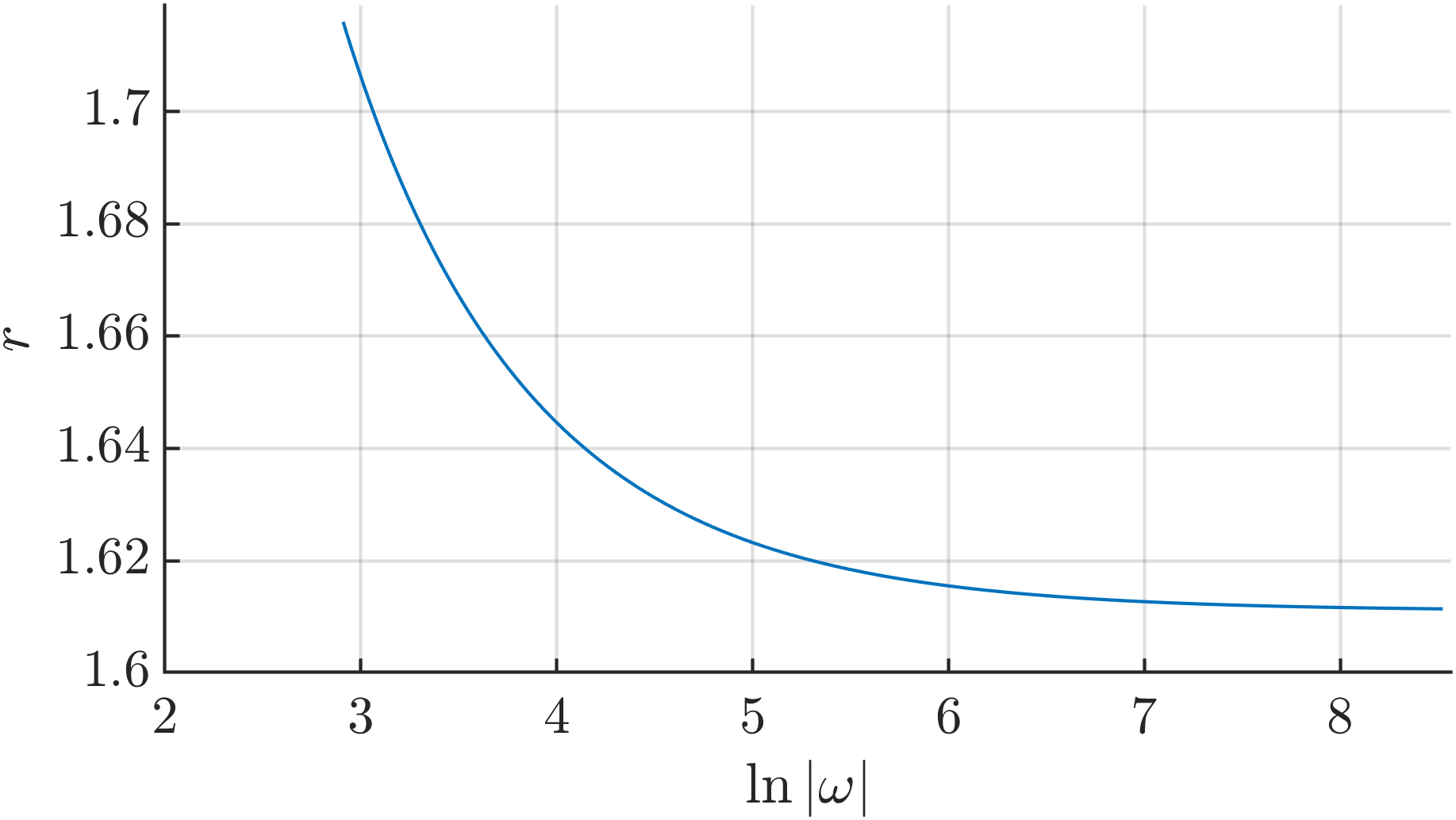}
    \includegraphics[trim={0mm 0mm 0mm 0mm},clip,width=0.48\textwidth]{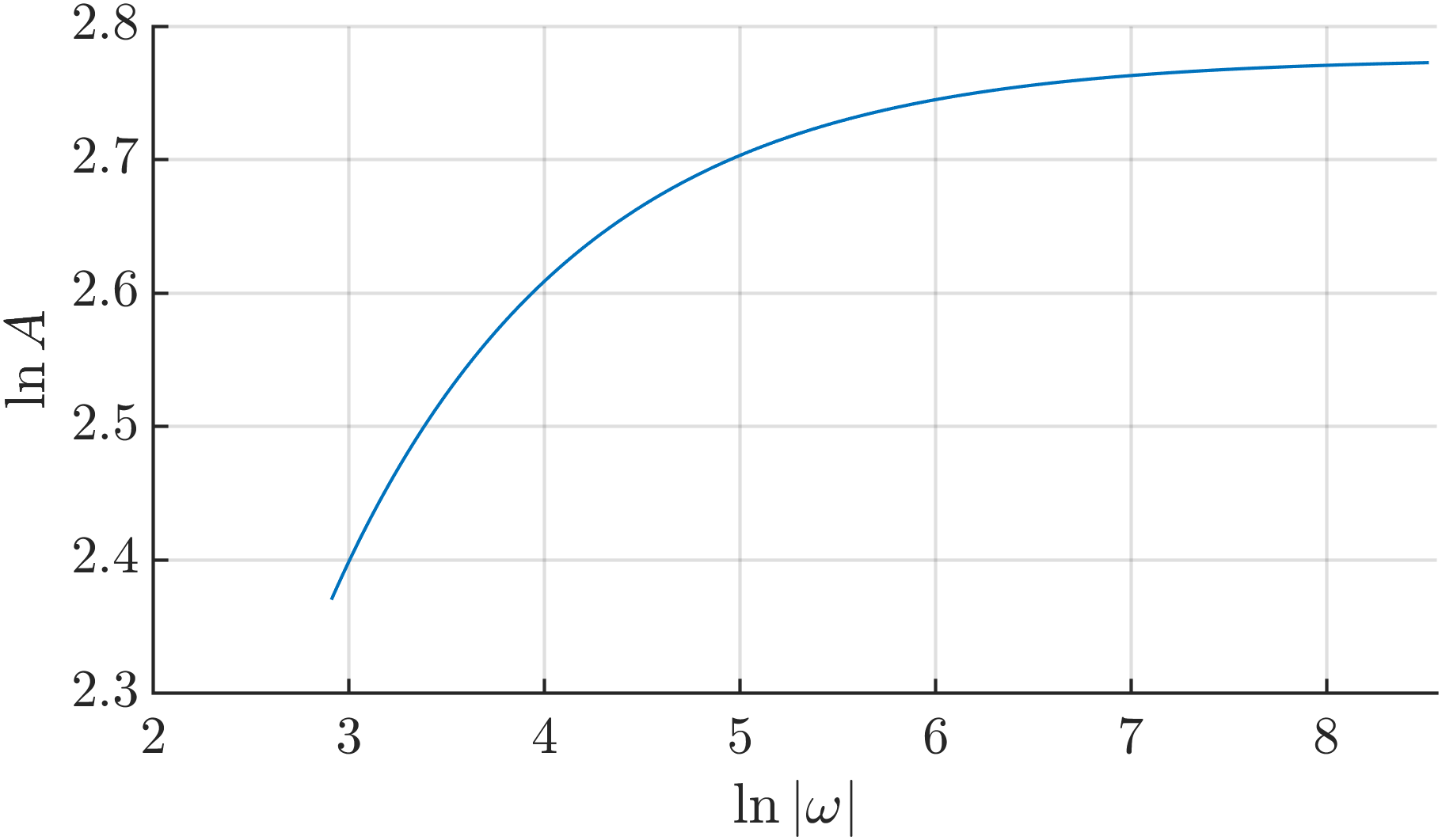}
    \caption{Values of $r$ (left) and $\ln A$ (right)
    obtained when taking pairs of points 
    $(K_i,\Delta \dot \xs_i)$ and $(K_{i+1},\Delta \dot \xs_{i+1})$
     varying (decreasing to zero) $K$ (or equivalently increasing $\ln |\omega|$).
    See the text for details.
    }
    \label{fig:real_logw_rA}
\end{figure}

Our next purpose is to fit the resulting data by an asymptotic formula. To do
so, we start with a {\sl naive} fit, that is a formula $\Delta \dot \xs\sim
c\cdot \exp (\omega d)$, or equivalently $\ln \Delta \dot \xs\sim \ln c+d\ln
|\omega|$, where $\omega=-\dfrac{1}{3\varepsilon ^2}$ and we look for the
constant $d$. Using divided differences we obtain that $d=\pi/2$. The next fit
asymptotic formula (inspired in formulas \eqref{deltap} and \eqref{eq:Aomega})
is
\begin{equation}\label{ajustxp}
\Delta \dot \xs\sim  \varepsilon A|\omega|^r\cdot \exp\left(
\frac{\omega\pi}{2}\right)
\end{equation}
or equivalently
\begin{equation}\label{ajustxpamblogw}
\ln |\Delta \dot \xs|-\ln  \varepsilon -\frac{\omega\pi}{2} \sim \ln A+  r\ln |\omega|.
\end{equation}
The crucial point is to determine the values of $\ln A$ and $r$. To do so, we
have proceeded following three  strategies: 

\emph{(i)} from the output data $(K,\Delta \dot \xs)$, we plot the points $(\ln
|w|,Y_{\Delta\dot \xs})$ where $Y_{\Delta\dot \xs}=\ln |\Delta \dot \xs|-\ln
\varepsilon -\dfrac{\omega\pi}{2}$  --see Figure~\ref{fig:logw_Y_xp}--. We
 clearly see that the points lie on a line. Actually, a
 linear regression approximation using formula \eqref{ajustxpamblogw}
provides
\begin{displaymath}
r= 1.6115,\qquad\ln A=2.772,\qquad A=15.99058.
\end{displaymath}
\emph{(ii)} Since formulas~\eqref{ajustxp} or~\eqref{ajustxpamblogw} are
asymptotic formulas, it seems quite reasonable to take for each pair of points
$(K_i,\Delta \dot \xs_i)$ and $(K_{i+1},\Delta \dot \xs_{i+1})$, the segment passing
through the points $(\ln |\omega_i|, Y_{\Delta \dot \xs_i})$ and $(\ln
|\omega_{i+1}|, Y_{\Delta \dot \xs_{i+1}})$ and look at the tendency of the slope
of such segments when $K$ decreases. We associate to this segment the
expression $ r\ln |\omega|+\ln  A$ and we plot the obtained values of $r$ in
Figure \ref{fig:real_logw_rA} left and the obtained values of $\ln  A$ in the
right figure. For the last pair of points we obtain 
\begin{displaymath}
r=1.6114670\dots,\qquad  \ln A=2.7726505\dots,\qquad A=16.000988\dots, 
\end{displaymath}
which are very similar values to those obtained in strategy \emph{(i)}.
 
A remark must be done at this point: notice that the decreasing values of $K$
considered range from $0.001$ to $10^{-8}$, or correspondingly the range of
increasing values in $\ln |w|$ is from $2.8$ to $8.6$. However the range of
values in $\varepsilon$ is from $0.14$ to $0.0076$, which turns out to be poor.
Nevertheless, we emphasize that the numerical computations have been done using
multiple precision arithmetics dealing with  up to five thousand digits and a
high order for the parametrization. So we reach a limit computing capacity for
taking smaller values of $K$, and therefore smaller values of $\varepsilon$.

\emph{(iii)} Actually, a more precise version of formula \eqref{ajustxp} can be
written as
\begin{displaymath}
\Delta \dot \xs\sim \varepsilon |\omega|^{r}\exp
    \left(\frac{w\pi}{2}\right)
    \left(A+A_1\varepsilon^{j_1}+A_2\varepsilon^{j_2}+\cdots
    \right)
\end{displaymath}
or equivalently as
\begin{equation}\label{ajustxpNOVA_despejada}
Z_{\Delta \dot \xs}:=
 \frac{1}{\varepsilon}
|\omega|^{-r}\Delta \dot \xs\exp \left(-\frac{\omega\pi}{2}\right)
    \sim A + A_{1}\varepsilon^{j_1} +
                A_{2}\varepsilon^{j_2} +
                A_{3}\varepsilon^{j_3} + \cdots
\end{equation}
and our goal is to obtain an (as much precise as possible) value of $A$, once
$r$ is given. Although we do not know the exact value of $r$, from the numerical
computations done in next Sections for the perturbed pendulum and for the
perturbed  amended pendulum, we will show (later on) that the \emph{guess} value for
$r$ is $r=1.611111\ldots=29/18$.

Using formula \eqref{ajustxpNOVA_despejada}  we have done several steps of extrapolation to obtain 
a robust value of $A$ (using the last 50 values of $K$ close
to $10^{-8}$) (see values of $(K,\varepsilon)$ in Table \ref{tablexpunt}, where only the smallest 10 values of  K are shown). 
That is, we have taken values of $(\varepsilon,Z_{\Delta \dot \xs})$
(see column 2 and 3 in the table) and
after some attempts, varying the values of $j_i$, we conclude that we get 
the best fit of $A$ with $j_i=2\cdot i$ and the value of
$A$ is $A=16.055307843$ (or $\ln A=2.776039450$). See columns  3, 4, 5 in the table.
 Columns 4 and 5 have been obtained taking the first two extrapolation steps.

\begin{table}[!ht]
    \centering
\begin{tabular}{|c|c|c|c|c|}
\hline
$K$ & $\varepsilon$ & $Z_{\Delta \dot \xs} $ & 1st extrapolation & 2nd extrapolation \\ \hline
   1.432E-8 &  8.312E-3 & {\bf 16.049}3664120828 & {\bf 16.055307}49436553 & {\bf 16.055307843}9006\\ \hline
   1.419E-8 &  8.293E-3 & {\bf 16.049}3937101860 & {\bf 16.055307}49756964 & {\bf 16.055307843}8344\\ \hline
   1.406E-8 &  8.274E-3 & {\bf 16.049}4208816007 & {\bf 16.055307}50074373 & {\bf 16.055307843}7819\\ \hline
   1.393E-8 &  8.255E-3 & {\bf 16.049}4479266667 & {\bf 16.055307}50388795 & {\bf 16.055307843}7170\\ \hline
   1.380E-8 &  8.236E-3 & {\bf 16.049}4748478517 & {\bf 16.055307}50700268 & {\bf 16.055307843}6392\\ \hline
   1.367E-8 &  8.217E-3 & {\bf 16.049}5016497992 & {\bf 16.055307}51008853 & {\bf 16.055307843}5483\\ \hline
   1.355E-8 &  8.198E-3 & {\bf 16.049}5283244532 & {\bf 16.055307}51314551 & {\bf 16.055307843}3730\\ \hline
   1.342E-8 &  8.179E-3 & {\bf 16.049}5548764688 & {\bf 16.055307}51617363 & {\bf 16.055307843}4478\\ \hline
   1.330E-8 &  8.160E-3 & {\bf 16.049}5813084435 & {\bf 16.055307}51917351 & {\bf 16.055307843}3717\\ \hline
   1.318E-8 &  8.141E-3 & {\bf 16.049}6076165450 & {\bf 16.055307}52214528 & {\bf 16.055307843}3357\\ \hline
\end{tabular}
\vskip 0.2cm
\caption{Approximate values of $A$ in formula \eqref{ajustxpNOVA_despejada} obtained by extrapolation.}\label{tablexpunt}
\end{table}

In summary, with the data at hand the specific values of $A$ and $r$ turn out to be
\begin{equation}\label{Ar}
 A=16.055307843\dots,\qquad r=1.6111111\ldots =29/18.
 \end{equation}


\subsection{Splitting in resonant coordinates}\label{subsec33}

Motivated, on one hand, by formulas~\eqref{eq:Aomega} and~\eqref{deltap} at
$t=0$ which provide explicit expressions of the splitting  for the perturbed pendulum, and on the 
other hand, by the numerical simulations done in the previous subsection, we
want to analyse the splitting (distance between the external manifolds $W^{u,e}$ and $W^{s,e}$ at the
Poincar\'e section $\Sigma$) when taking into account, not the original
variables $(\xs,\ys,\pxs,\pys)$ but the resonant ones $(x,y,p,q)$. So, for
each value of $K$, once we obtain the point $P^u$ from the computations
described in the previous subsection for the CP problem, we apply the change of
coordinates from Levi-Civita $(u,v,u',v')$ to synodical ones $(\xs,\ys,\pxs,\pys)$ and
to resonant coordinates $(x,y,p,q)$, so now $P^u$ becomes
$(x^u,y^u,q^u,p^u)$, with $x^u\approx\pi$. Taking the same range of values
of $K$ tending to zero, more precisely a set of values $K_i$, $i=1,\dots,N$ with
$K\in [10^{-8},10^{-3}]$ and using multiple precision
computations, we compute the splitting which now we define as $\Delta p=2p^u$. %
(notice that it is not $p^u-p^s$, since $x$ is not exactly equal to $\pi$).
Taking into account formulas \eqref{deltap} and \eqref{eq:Aomega}, where  we skip the $+$ notation since we only consider the external manifolds, we fit  the
splitting by the asymptotic formula 
\begin{equation}\label{deltapu}
\Delta p\sim \varepsilon \bar A|\omega|^{\bar r}\exp
\left(\frac{\omega\pi}{2}\right)
\end{equation}
or equivalently
\begin{equation}\label{logdeltapu}
\ln |\Delta p|-\ln  \varepsilon -\frac{\omega\pi}{2}
\sim \ln \bar A + \bar{r}\ln |\omega|.
\end{equation}
Our next goal is to determine $\bar A$ and $\bar r$. We proceed applying the
three strategies mentioned above.

\begin{figure}[ht!]
    \centering
    \includegraphics[trim={0mm 0mm 0mm 0mm},clip,width=0.48\textwidth]{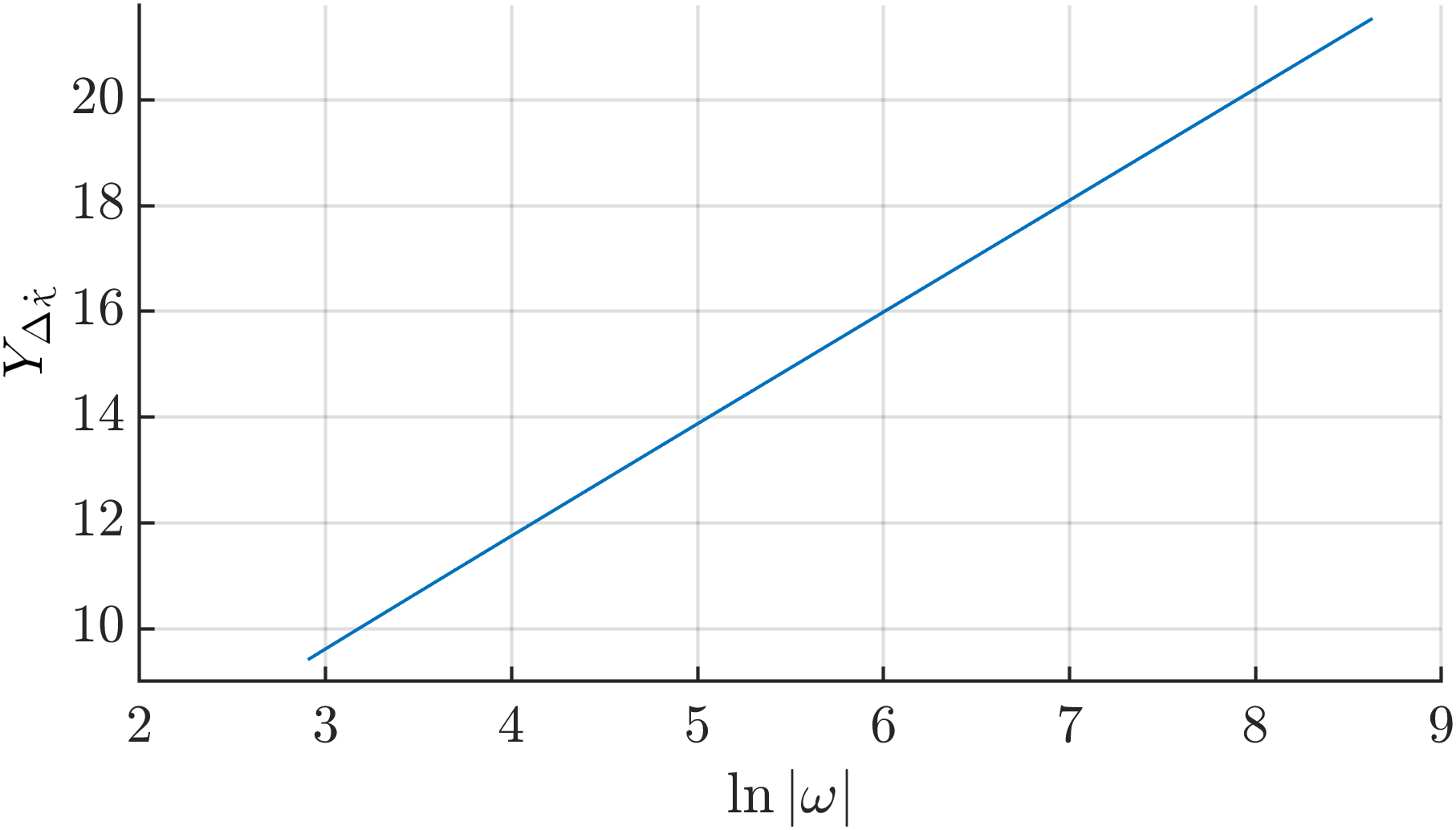}
    \caption{Points $(\ln |\omega|,Y_{\Delta p})$  varying (decreasing to zero)
    $K$ (or equivalently increasing $\ln |\omega|$). Notice the points {\sl lying} on the line $Y=3.325+2.111\ln|\omega|$.
    } \label{fig:logw_Y_Deltap}
\end{figure}

\emph{(i)} We plot the points $(\ln |\omega|,Y_{\Delta p})$, with $Y_{\Delta p}:=
\ln |\Delta p|-\ln  \varepsilon -\dfrac{\omega\pi}{2}$, see Figure
\ref{fig:logw_Y_Deltap}.
Again we observe the good fit of points
 $(\ln |\omega|,Y_{\Delta p})$ by 
 a line. A linear regression approximation using
 formula \eqref{logdeltapu} provides 
$\ln \bar A=3.325$ and $\bar r=2.111$.

\begin{figure}[ht!]
    \centering
    \includegraphics[trim={0mm 0mm 0mm 0mm},clip,width=0.48\textwidth]%
        {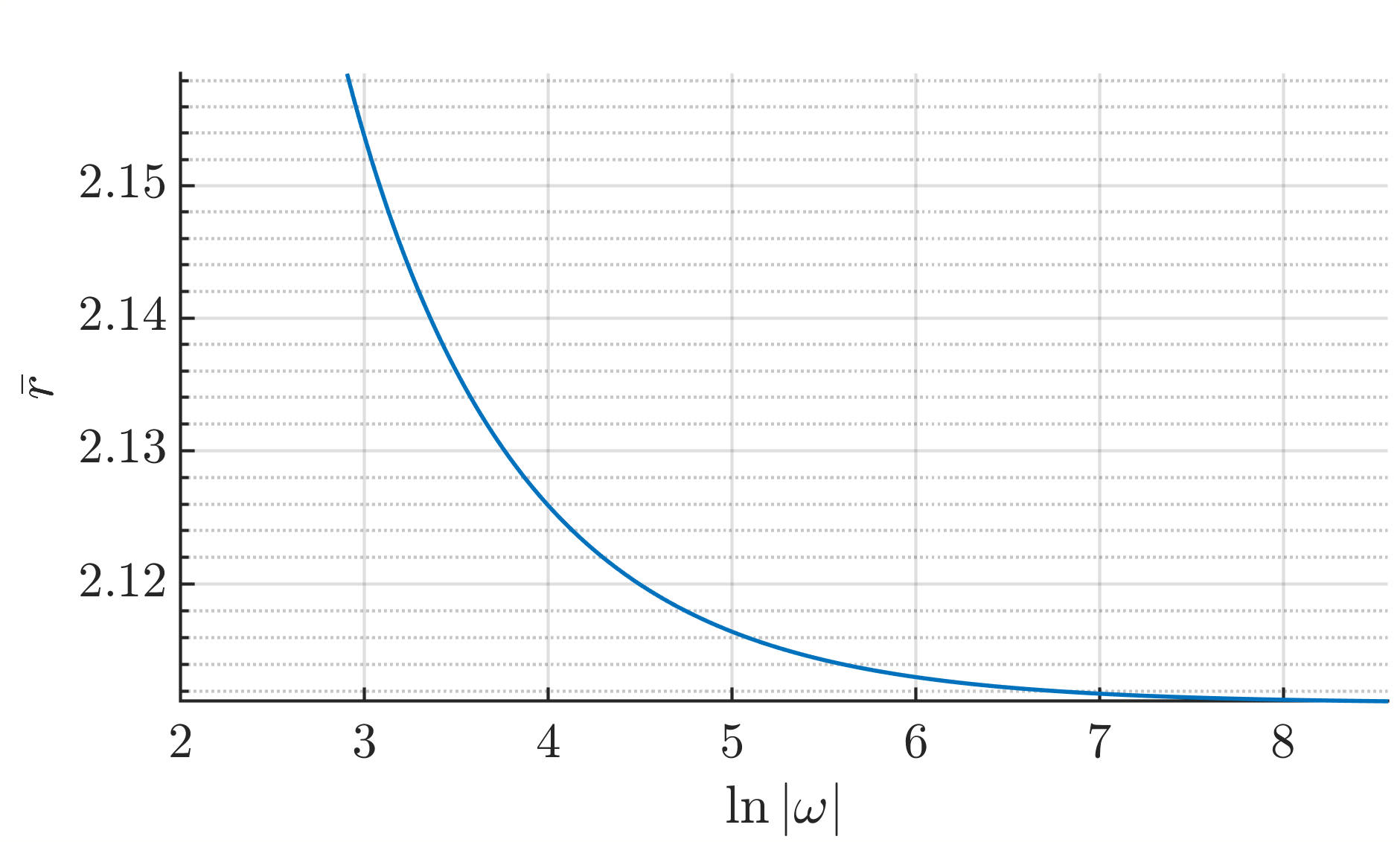}
    \includegraphics[trim={0mm 0mm 0mm 0mm},clip,width=0.48\textwidth]%
        {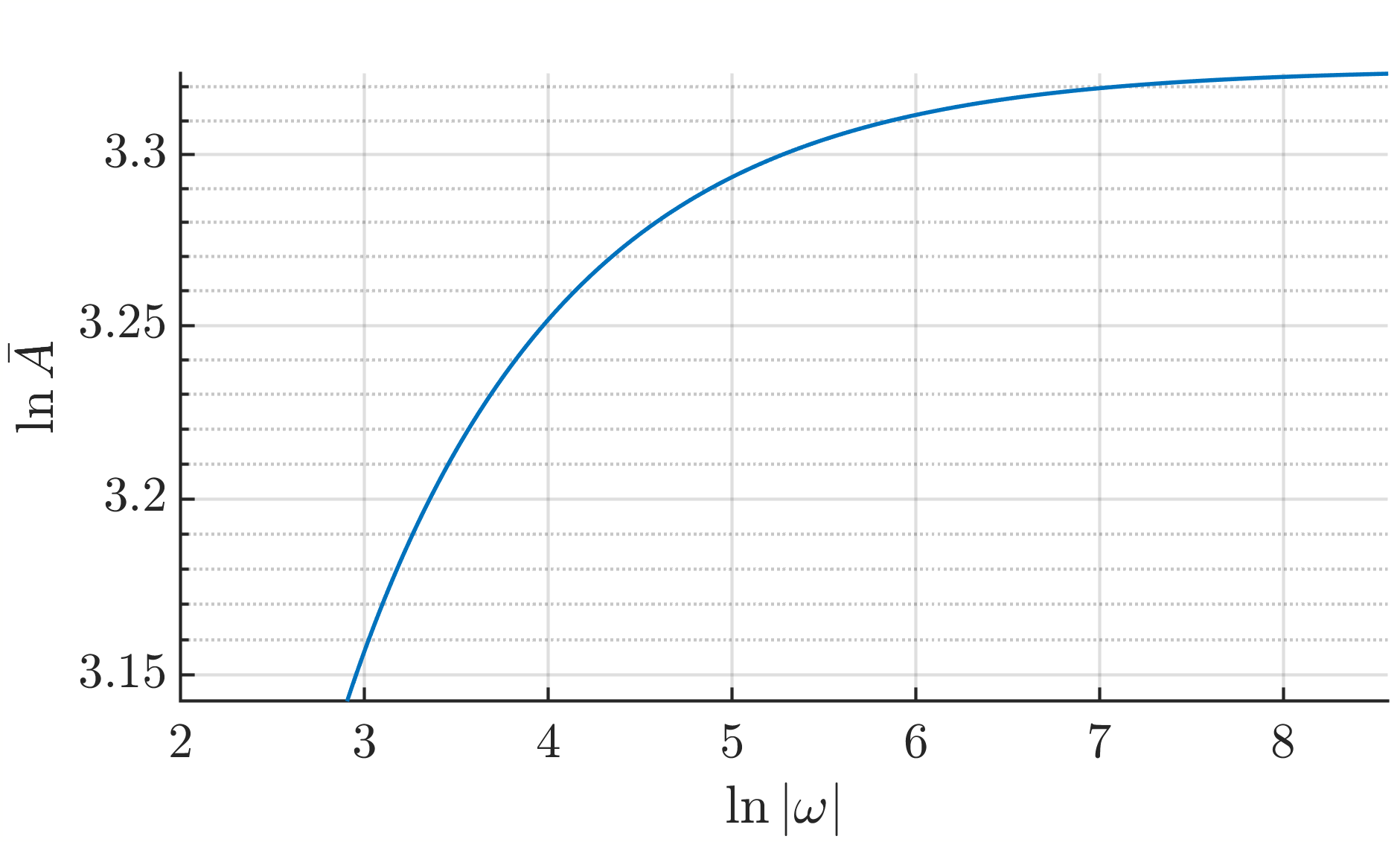}    
    \caption{Values of $\bar r$ (left) and $\ln\bar A$ (right) obtained when
    taking pairs of points $(K_i,\Delta p_i)$ and $(K_{i+1},\Delta 
    p_{i+1})$, in resonant coordinates, varying (decreasing to zero) $K$ (or
    equivalently increasing $\ln |\omega|$). See the text for details.}
    \label{fig:realRlogA}
\end{figure}

\emph{(ii)} Taking each successive pair of points $(K_i,\Delta p_i)$ and
$(K_{i+1},\Delta p_{i+1})$ and the segment passing through the points $(\ln
|\omega_i|, Y_{\Delta p_i})$ and $\left(\ln |\omega_{i+1}|,
Y_{\Delta {p_{i+1}}}\right)$. We associate to this segment the  expression $\bar
r\ln |\omega|+\ln \bar A$ and we plot the obtained values of $\bar r$ in Figure
\ref{fig:realRlogA} left and the obtained values of $\ln \bar A$ in the right
figure. With the data at hand, the last values obtained for  $\ln\bar A$ and  $\bar
r$ are $\ln\bar{A}=3.323854\dots$, $\bar{r} = 2.111267\dots$, which are very similar to
those values in strategy \emph{(i)}. 

\emph{(iii)} We now make use of the fit formula
\begin{displaymath}
    \Delta p\sim \varepsilon |\omega|^{\bar r} \exp 
    \left(\frac{\omega\pi}{2}\right)
    \left(\bar A +
          \bar A_1\varepsilon^{j_1} + 
          \bar A_2\varepsilon^{j_2} + 
          \dots
          \right),
\end{displaymath}
or equivalently as 
\begin{equation}\label{ajustpNOVA_despejada}
Z_{\Delta p}:=
 \frac{1}{\varepsilon}
|\omega|^{-\bar r}\Delta p\exp\left(-\frac{\omega\pi}{2}\right)\sim
           \bar A + 
           \bar A_{1}\varepsilon^{j_1} + 
           \bar A_{2}\varepsilon^{j_2} +
           \bar A_{3}\varepsilon^{j_3} +\dots
\end{equation}
and as we have done previously with $\Delta \dot x$, we want to apply
extrapolation to obtain an (as much precise as possible) value of $\bar A$, once
$\bar r$ is given. Again, we need a guess value which turns out to be $\bar{r}=2.11111\ldots=19/9$. This value will be justified in Section \ref{PQP}.

Using formula \eqref{ajustpNOVA_despejada}  we have done two steps of
extrapolation to obtain a robust value of $\bar A$ (again using the same last
$50$ values of $K$ close to $10^{-8}$, see Table \ref{tablep}) where we have taken
$j_i=2\cdot i$  (which provides the best fit).
 We list the values of $\bar
A$ obtained through the successive extrapolation steps (we only provide the last
10 values but we have used $50$):

\begin{table}[!ht]
    \centering
\begin{tabular}{|c|c|c|}
\hline
 $Z_{\Delta p} $ & 1st extrapolation & 2nd extrapolation  \\ \hline
   27.8040812616518 &   {\bf 27.808609}32094579 &  {\bf 27.80860891}67525   \\ \hline
   27.8041020671644 &   {\bf 27.808609}31723887 &  {\bf 27.80860891}66344   \\ \hline
   27.8041227760921 &   {\bf 27.808609}31356489 &  {\bf 27.80860891}65020   \\ \hline
   27.8041433886943 &   {\bf 27.808609}30992373 &  {\bf 27.80860891}63859   \\ \hline
   27.8041639068522 &   {\bf 27.808609}30631499 &  {\bf 27.80860891}62863   \\ \hline
   27.8041843341052 &   {\bf 27.808609}30273799 &  {\bf 27.80860891}62036   \\ \hline
   27.8042046643138 &   {\bf 27.808609}29919275 &  {\bf 27.80860891}60301   \\ \hline
   27.8042249010261 &   {\bf 27.808609}29567928 &  {\bf 27.80860891}59491   \\ \hline
   27.8042450462224 &   {\bf 27.808609}29219689 &  {\bf 27.80860891}58544   \\ \hline
   27.8042650969822 &   {\bf 27.808609}28874547 &  {\bf 27.80860891}57148   \\ \hline
\end{tabular}
\vskip 0.2cm
\caption{
Approximate values of $\bar A$ in formula \eqref{ajustpNOVA_despejada} obtained by extrapolation.
}\label{tablep}
\end{table}

In summary,  with the data at hand, the values of $\bar A$ and $\bar r$ turn out to be

\begin{equation}\label{Arbar}
\bar
A=27.80860891\dots \quad (\ln \bar A=3.325345645\dots), \qquad \bar r=2.11111\ldots=19/9,
\end{equation}

{\bf Remark}
We notice that the Poincar\'e section $y=0$ does not correspond to $x=\pi$ but
$x=\pi+\delta$. However, numerical computations reveal that  
\begin{displaymath}
 \delta=C\varepsilon^{2}|\omega|^{\bar r}\exp\left(\frac{\omega\pi}{2}\right)   
\end{displaymath}
so in terms of comparison of the fit formulas for the splitting, we can deal
with the Poincar\'e section $x=\pi$ or $x=\pi+\delta$, indistinctly (since
the contributionn of $\delta$ is of order $\varepsilon^2$, apart from the
exponentially small term).

At this point, two natural questions arise: \emph{(i)} Can we explain the
difference between $r=1.611\dots$ and $\bar{r} = 2.111\dots$? \emph{(ii)} Does the
Melnikov theory predict the asymptotic formula \eqref{deltapu} or
\eqref{logdeltapu} with these values of $\bar r$ and $\ln \bar A$?

Next section is devoted, in particular, to answer the second question.

Concerning the first question let us prove the following result:

At the Poincar\'e section $\Sigma:\quad \ys=0$, we have
\begin{displaymath}
\dot \xs\approx \varepsilon p
\end{displaymath}
or equivalently
\begin{displaymath}
\Delta\dot \xs\approx \varepsilon \Delta p.
\end{displaymath}
To prove this assertion, let us assume that $\xs=1+O(\delta)$, $\dot \ys=O(\delta)$,
and $\dot \xs >0$ being $\delta$ a small quantity. Then
\begin{align*}
a &=\frac{\xs}{2-\xs(\dot \xs^2+\dot \ys^2)-2\dot \ys\xs^2-\xs^3}=
\frac{1+O(\delta)}{2-(1+O(\delta))(\dot \xs^2+O(\delta ^2))-2O(\delta ^2)-1+O(\delta)} \\
&=\frac{1+O(\delta)}{1-\dot \xs^2(1+O(\delta))+O(\delta)}
 =\frac{1+O(\delta)}{1-\dot \xs^2+O(\delta)}\\
 &=\frac{1}{1-\dot \xs^2}+O(\delta),\\[8pt]
L&=\frac{1}{\sqrt{1-\dot \xs^2}}+O(\delta)\\[8pt]
G&=\xs(\dot \ys+\xs)=1+O(\delta)
\end{align*}
so
\begin{displaymath}
L-G=\frac{\dot \xs^2}{\sqrt{1-\dot \xs^2}+1-\dot \xs^2}+O(\delta)
\end{displaymath}
and
\begin{gather*}
\sqrt{2(L-G)}=\frac{\sqrt{2}\dot \xs}{\sqrt{\sqrt{1-\dot \xs^2}+1-\dot \xs^2}}+O(\delta)\\
%
e=\sqrt{1-\frac{G^2}{a}}=\dot \xs+O(\delta).
\end{gather*}
Moreover the vector ${\bf e}=(e^{(1)},e^{(2)})$, with $e^{(2)}=-\xs\dot \xs (\dot \ys+\xs)=
 -\dot \xs+O(\delta)$. So finally
\begin{displaymath}
p=\frac{\sqrt{2(L-G)}}{\varepsilon} \sin g =\frac{\sqrt{2(L-G)}}{\varepsilon} 
\frac{e^{(2)}}{e}\approx \frac{1}{\varepsilon}\dot \xs,
\end{displaymath}
that is 
\begin{equation}\label{relacio_p_dotxfinal}\dot \xs\approx \varepsilon p.
\end{equation}

Expression \eqref{relacio_p_dotxfinal} allows us to relate the exponent $r$ and $\bar r$ and the constants $A$ and $\bar A$ in the two fit formulas:
\begin{alignat*}{3}
& \dot \xs\sim A\varepsilon |\omega|^r\exp\left(\frac{w\pi}{2}\right),\quad 
&& \text{with} \quad && r=1.61111\ldots=\frac{28}{19},\\ 
& p\sim\bar A\varepsilon |\omega|^{\bar r}\exp\left(\frac{w\pi}{2}\right),\quad 
&&  \text{with}\quad && \bar r=2.1111\ldots =r+\frac{1}{2}=\frac{19}{9}. 
\end{alignat*}
Using $\dot \xs=\varepsilon p$, we obtain 
\begin{equation}\label{relarrbar}
\bar r=r+\frac 12,\qquad A\approx \frac{\bar A}{\sqrt{3}}
\end{equation}
which agrees with the relation $\bar r=2.11111\ldots=1.611111+0.5=r+0.5$ and if we compare the numerical values obtained
$$ A=16.0553078\dots,\qquad \bar A=27.8086089\dots
 $$
we precisely get $\frac{\bar A}{\sqrt{3}}=16.0553078\dots$
which coincides with the value of $A$ up to the first 9 digits.


\section{
    The perturbed pendulum. Other possible models
}\label{secpendol}
This Section is devoted to answer the following question: does the Melnikov
theory predict the asymptotic formula \eqref{deltapu}  with these values of
$\bar r$ and $\ln \bar A$?

To answer this question we take, as the most natural model, the first order in
$\varepsilon$ Hamiltonian. More precisely, in Section \ref{subsecRC} we have
written the Hamiltonian of the CP problem in resonant coordinates $(x,y,q,p)$
(see \eqref{eq:HPML1}). We now consider the Hamiltonian as an expansion in
powers of $\varepsilon = \left(K/3\right)^{1/4}$:
\begin{align}\label{eq:H0H1H2}
H &= \omega\frac{q^2+p^2}{2}+\frac{y^2}{2}+\cos x-1+
 \varepsilon\left[ \frac{3}{2}q-\frac{q}{2}\cos 2x-\frac{p}{2}\sin 2x\right]
 \nonumber\\[8pt]
 &\quad+ \varepsilon^2\left[ \left(-\frac{3}{8}q^2-\frac{5}{8}p^2\right)\cos x
  +2y\cos x +\frac{1}{4}qp\sin x + \frac{3}{8}\left(q^2-p^2\right)\cos 3x + 
\frac{3}{4}qp\sin 3x - \frac{2}{3}y^3\right] +\cdots\nonumber\\[8pt]
&:= H_0+\varepsilon H_1+\varepsilon ^2H_2+\cdots
\end{align}
and we take as the simplest model  the first order expansion,  

\begin{equation}\label{rotorpendulum} 
 \hat H(x,y,q,p)=H_0(x,y,q,p)+\varepsilon H_1(x,y,q,p)
 \end{equation}
 with
\begin{equation}\label{pend}
H_0(x,y,q,p)=\omega\frac{q^2+p^2}{2}+\frac{y^2}{2}+\cos x-1
\end{equation}
and
\begin{equation}\label{h1resso}
H_1(x,y,q,p)=\frac{3}{2}q - \frac{q}{2}\cos 2x-\frac{p}{2}\sin 2x, 
\end{equation}
and we call \eqref{rotorpendulum}  the perturbed pendulum (although $H_0$ is a rotor times a pendulum).

 We want to test the Melnikov theory (that is, to check if the Melnikov integral
 \eqref{FAw} and \eqref{deltap} using Hamiltonian \eqref{rotorpendulum} predicts
 the asymptotic behavior for the splitting provided by \eqref{deltapu} or
 \eqref{logdeltapu}). A clear advantage of this Hamiltonian is that we know the
 explicit expression of the separatrices (given by
\eqref{eq:StandardSeparatrices})  in $(q,p)=(0,0)$ of the integrable Hamiltonian
$H_0(x,y,q,p)$, which is a rotor times a pendulum. We observe that the
associated system of ODE to $\tilde H$ has two equilibrium points:
$L_-=(0,0,3\varepsilon^3,0)$ and $L_+=(2\pi,0,3\varepsilon^3,0)$. We will focus
our attention on the external unstable manifold of $L_-$ (that will start
describing a curve in the $(x,y)$ projection with $y>0$, $x>0$), and the
external stable manifold of $L_+$ (that will start describing a curve in the
$(x,y)$ projection with $y>0$, $x<2\pi $). We want to measure the distance
between them (the splitting) at a given Poincar\'e section. Due to the symmetry
$(x,y,q,p,t)\to (2\pi -x,y,q,-p,-t)$ satisfied by the associated system of ODE,
a Poincar\'e section that turns out to be convenient is $x=\pi$, $y>0$. From now
on we will omit the \emph{external} mention when discussing the manifolds and we
will simply refer to them as the manifolds.

\subsection{Melnikov integral for the perturbed pendulum}\label{41}
Formula \eqref{eq:Aomega} provides the Melnikov integral for the separatrix of
the pendulum and Formula \eqref{eq:deltapmel} provides (the linear approximation
of) the splitting. In particular at $t=0$ we have 

\begin{equation}\label{eq:homoL1Anal_compl}
  \Delta p
  \sim -\varepsilon^m\frac{8\pi}{3}\omega^{3}\left(1-\frac{2}{\omega^{2}}\right)
  \frac{\ee^{\pi\omega/2}}{1-\ee^{2\pi\omega}},
\end{equation}
with $m=1$.
So, comparing the splitting fit formulas \eqref{deltapu}  and \eqref{eq:homoL1Anal_compl} we obtain $\bar r=3$ from Melnikov formulation, instead of the expected value $\bar r=2.111\dots$
and $\bar A=\frac{8\pi}{3}=8.377\dots$
instead of $\bar A=27.808\dots$ (see \eqref{Arbar}).
 We can conclude that 
the numerical results for $\Delta p$ do not coincide with the theoretical prediction 
      given by the Melnikov formula.

\subsubsection{Analysis of the Toy CP problem with 
\texorpdfstring{$a=0$}{a=0}}\label{subsubsec411} 
Although the Hamiltonian $\hat H$ in \eqref{rotorpendulum} turns out not to be a
good approximation to get the asymptotic fit for $\Delta p$ for the CP problem,
it is interesting to analyze the following Hamiltonian---per se and to compare
our results with previous ones in the literature---, taking the same integrable
part plus a more general perturbation, that is the Toy CP problem with $a=0$,
\begin{equation}\label{pendol}
\tilde H(x,y,q,p)=H_0(x,y,q,p)+\varepsilon ^m H_1(x,y,q,p)
\end{equation}
where $m\ge 1$, $H_0$ and $H_1$ are given by \eqref{pend} and \eqref{h1resso}. 
 Let us note first that there are theoretical
results~\cite{Baldoma2006,BaldomaFGS2012} for Hamiltonian systems with
$1+\sfrac{1}{2}$ degrees of freedom that can be directly applied to the Toy CP
problem with $a=0$~\eqref{pendol}. According to these theoretical results, for
$m\geq 4$ the Melnikov prediction with the standard pendulum gives the correct
measure~\eqref{eq:homoL1Anal_compl} for the splitting.

Let us consider the equilibrium points of $\tilde H$,
$L_-=(0,0,3\varepsilon^{2+m},0)$
and
$L_+=(2\pi,0,3\varepsilon^{2+m},0)$. 
Our goal, in this subsection, is to check if the splitting between the unstable manifold of $L_-$
(in $y\ge 0$) and the stable manifold of $L_+$ (in $y\ge 0$)
measured at the Poincar\'e section $\tilde \Sigma$ defined by $x=\pi$, $y>0$  
can be predicted by the Melnikov formula \eqref{eq:homoL1Anal}
obtained analytically regardless the value of $m$.

To do so, we fix a value of $m$, and for each value of $K$ (or equivalently $\varepsilon$ or $\omega$):

1.~We compute the associated equilibrium point $L_-=(0,0,3\varepsilon^{2+m},0)$ of the corresponding Hamiltonian system.

2.~We compute, using the parameterization method and multiple precision arithmetics,
 the unstable manifold $\Cu(L_-)$ up to the Poincar\'e section $\widetilde \Sigma$.  
 We call $(x^{u},y^{u},q^{u},p^{u})$ the corresponding point
 with $x^u=\pi$.

 3.~We would proceed similarly with the stable manifold of $L_+$ and we would compute the intersection point 
  $(x^{s},y^{s},q^{s},p^{s})$
 with $\widetilde \Sigma$.
However, due to the symmetry
\begin{equation}\label{simetria}
(x,y,q,p,t) \to (2\pi -x,y,q,-p,-t)
\end{equation}
only the unstable manifold $\Cu(L_-)$ needs to be computed since
at the Poincar\'e section we have  $x^u=x^s=\pi$, $q^u=q^s$, $p^u=-p^s$ and $y^u=y^s$.
 We are interested in the splitting between 
 $\Cu(L_-)$ and $\Cs(L_+)$ at $\widetilde \Sigma$, that is,
  the value
\begin{displaymath}
   \Delta p=p^u-p^s=2p^u 
\end{displaymath}
4.~We vary the values of $K$ tending to 0 (or  $\varepsilon$ or $\omega$), and we obtain 
a set of points  $(K_i,\Delta p_i)$, for $i=1,\dots,N$. In the numerical simulations we have taken decreasing values of $K$
up to $3\cdot 10^{-8}$. Again order of thousands for the parametrization and multiple precision with thousands of digits are required.

\begin{figure}[ht!]
    \centering
    \includegraphics[trim={0mm 0mm 0mm 0mm},clip,width=0.48\textwidth]{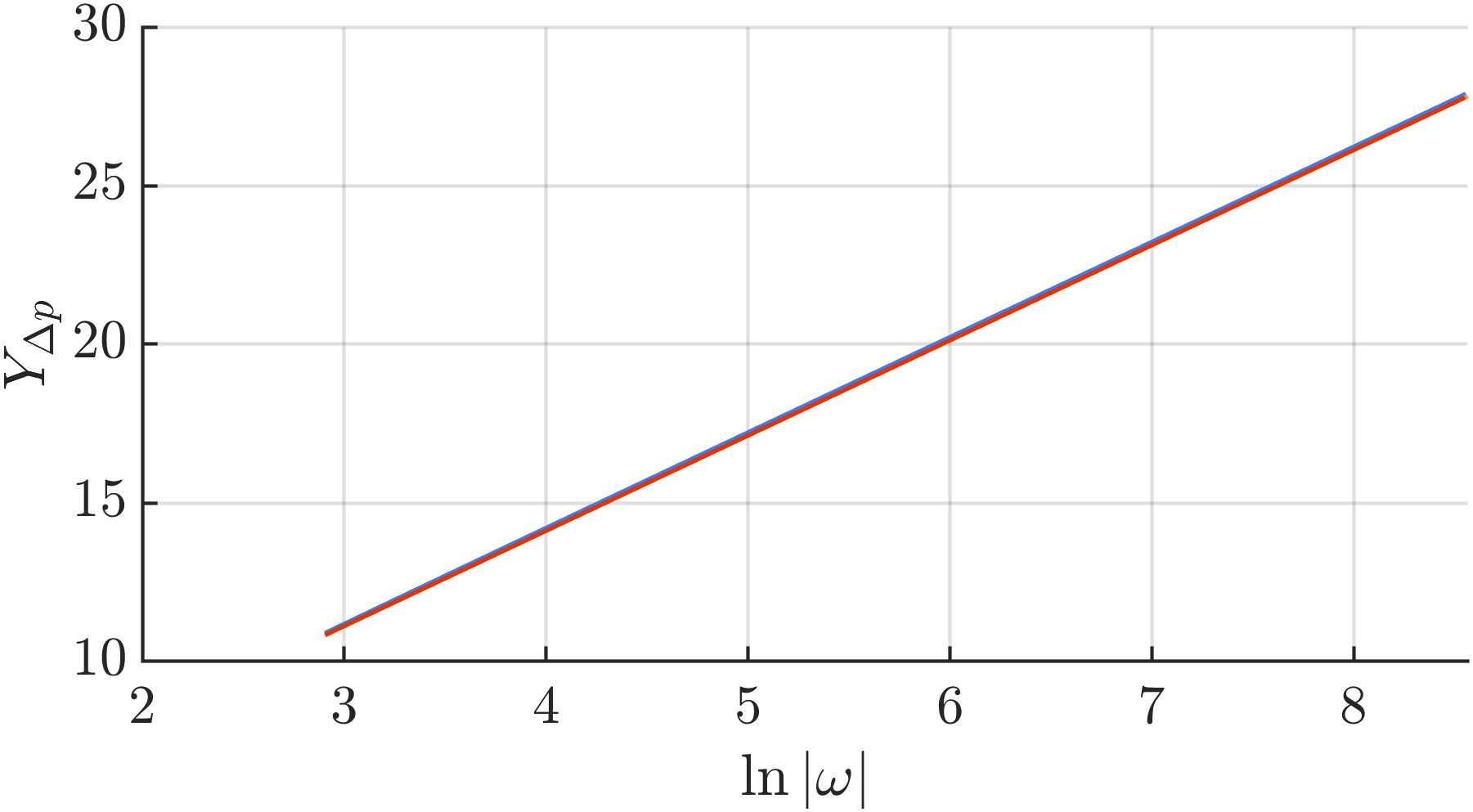}
    \caption{Points $(\ln|\omega_i|,Y_{\Delta p_i})$ for
     $m=1,2,3,4,5,6$. Notice they {\sl overlap} on a line.
}
    \label{fig:refer}
\end{figure}

5.~We fit the numerical computed values by a formula of type,
\begin{equation}\label{Deltap_m}
\Delta p\sim \varepsilon \tilde A|\omega|^{\tilde r}\cdot 
\exp\left(\frac{\omega\pi}{2}\right)
\end{equation}
or equivalently,
\begin{displaymath}
\ln |\Delta p|-\ln  \varepsilon -\frac{\omega\pi}{2} \sim \ln \tilde A+ \tilde r\ln |\omega|.
\end{displaymath}
Let us discuss the results obtained. Notice that we want to check two different
values: the exponent $\tilde r$ which should be equal to $3$ and the constant
$\tilde A$ which should be $8\pi/3$ or $\ln \tilde A=2.125591\dots$. In Figure
\ref{fig:refer} we show a first (rough) plot with  the points $(\ln
|\omega_i|,Y_{\Delta p_i})$ with $Y_{\Delta p_i}:=\ln |\Delta p_i|-\ln
\varepsilon_i -\omega_i\pi/2$ (and $\varepsilon _i=(K_i/3)^{1/4}$ and
$\omega_i=-1/(3\varepsilon _i^2)$), where we have taken some values of $K$
ranging from $0.005$ to $5\cdot 10^{-8}$. 
 Taking such a big set of values of $K$ for each value of $m$, we can see how the set of points for each case
 $m=1,2,3,4,5,6$ overlap on the line $Y=2.1255+3\ln|\omega|$.

\begin{figure}[ht!]
    \centering
    \includegraphics[trim={0mm 0mm 0mm 0mm},clip,width=0.48\textwidth]{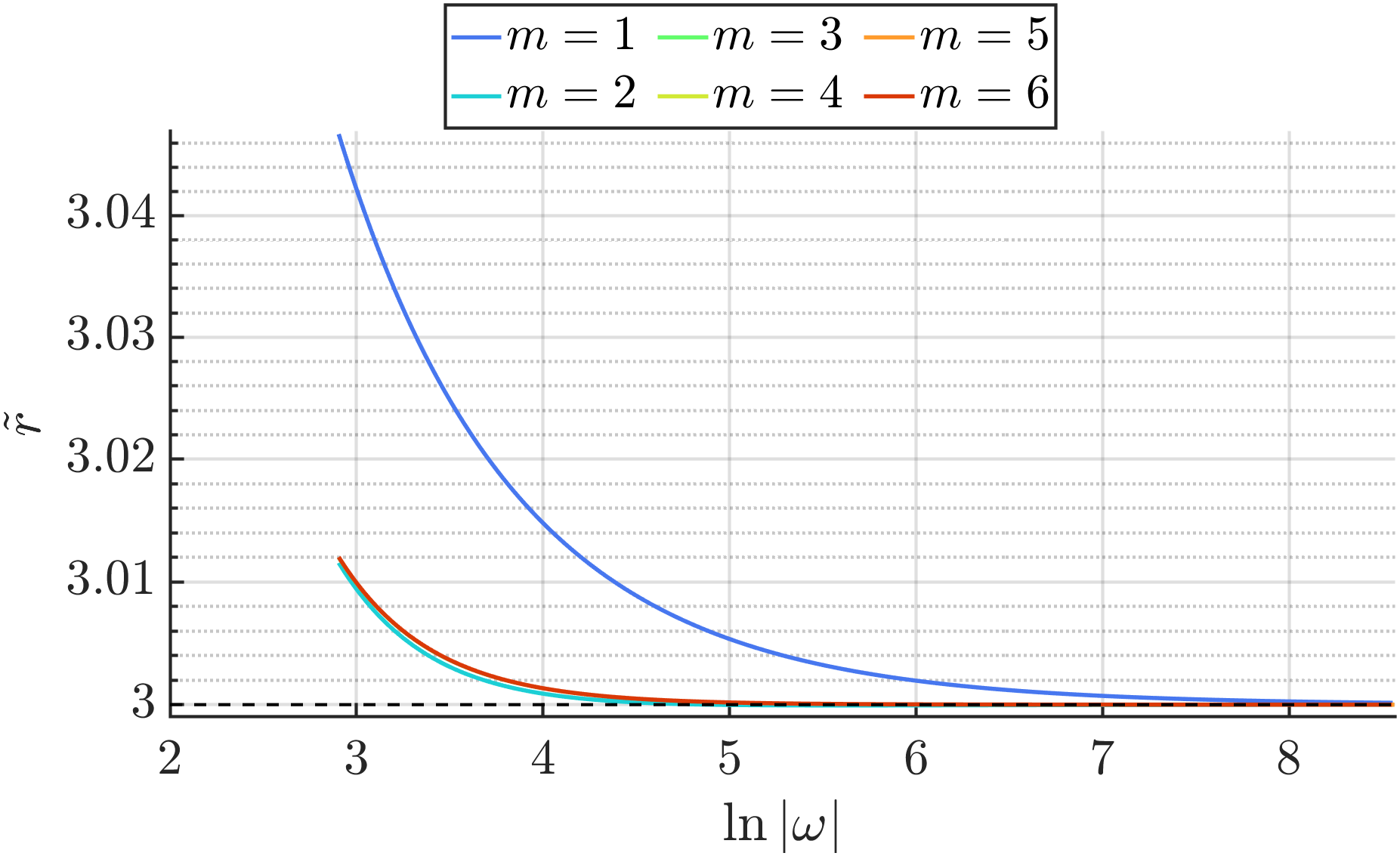}
    \includegraphics[trim={0mm 0mm 0mm 0mm},clip,width=0.48\textwidth]{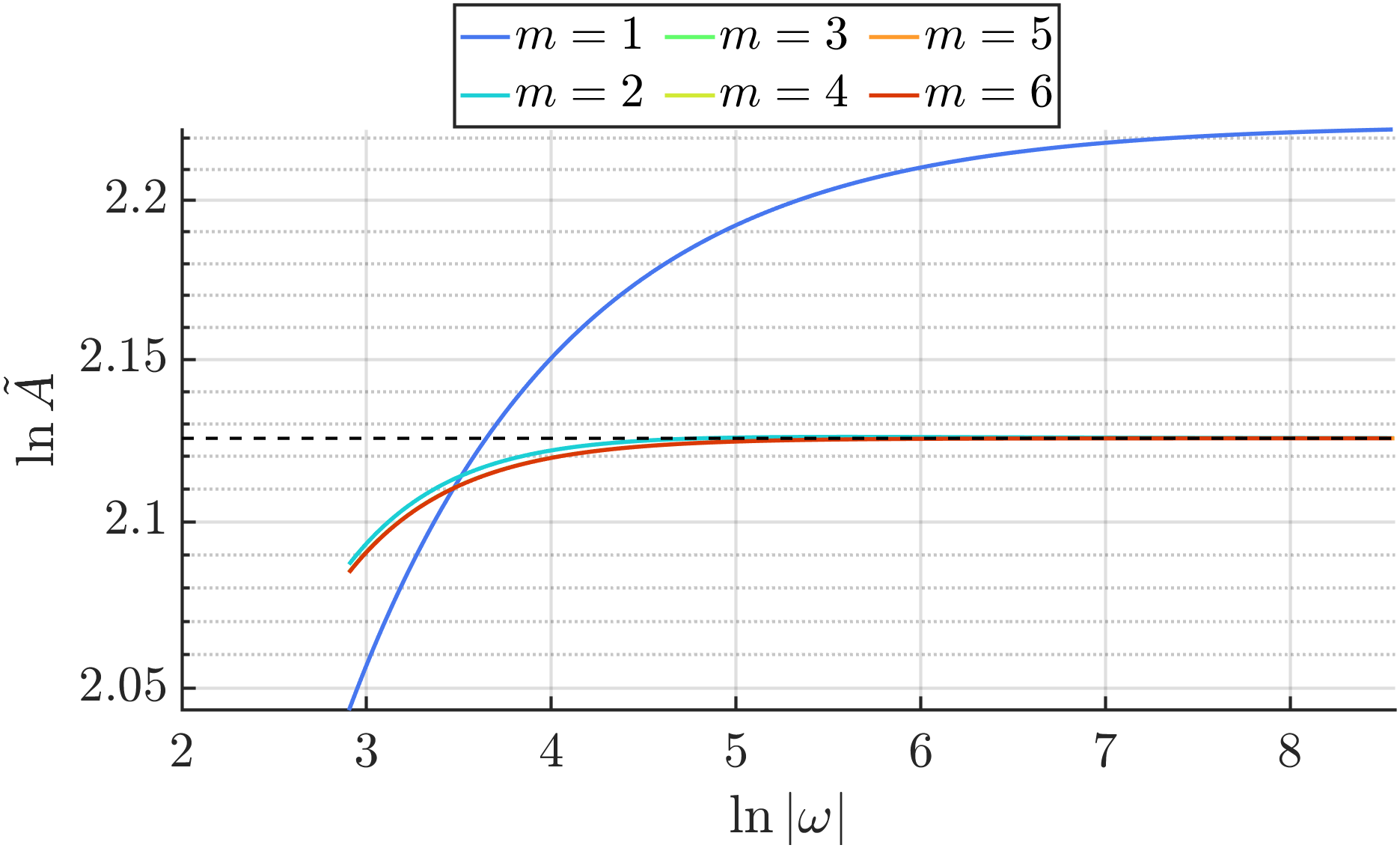}
    \caption{The perturbed pendulum. Points $(\ln |\omega|,\tilde r)$
     (left), and $(\ln |\omega|,\ln \tilde A)$ (right) 
    obtained for $m=1,\ldots,6$. The dashed line (on the right figure) corresponds to the value $\ln (8\pi/3)=2.1255\dots$.}
    \label{fig:CasAm1a6}
\end{figure}

\begin{figure}[ht!]
    \centering
    \includegraphics[trim={0mm 0mm 0mm 0mm},clip,width=0.48\textwidth]{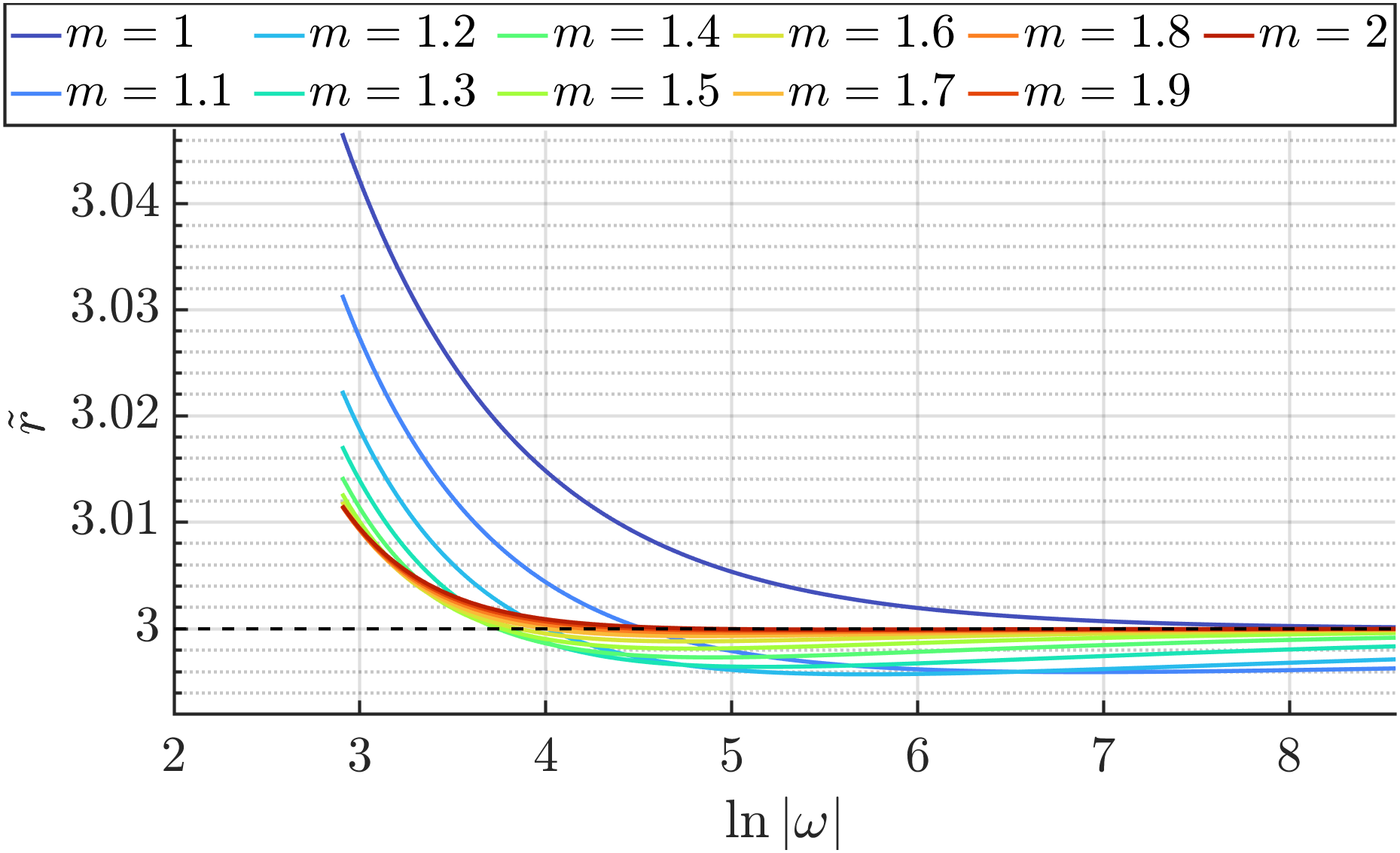}
    \includegraphics[trim={0mm 0mm 0mm 0mm},clip,width=0.48\textwidth]{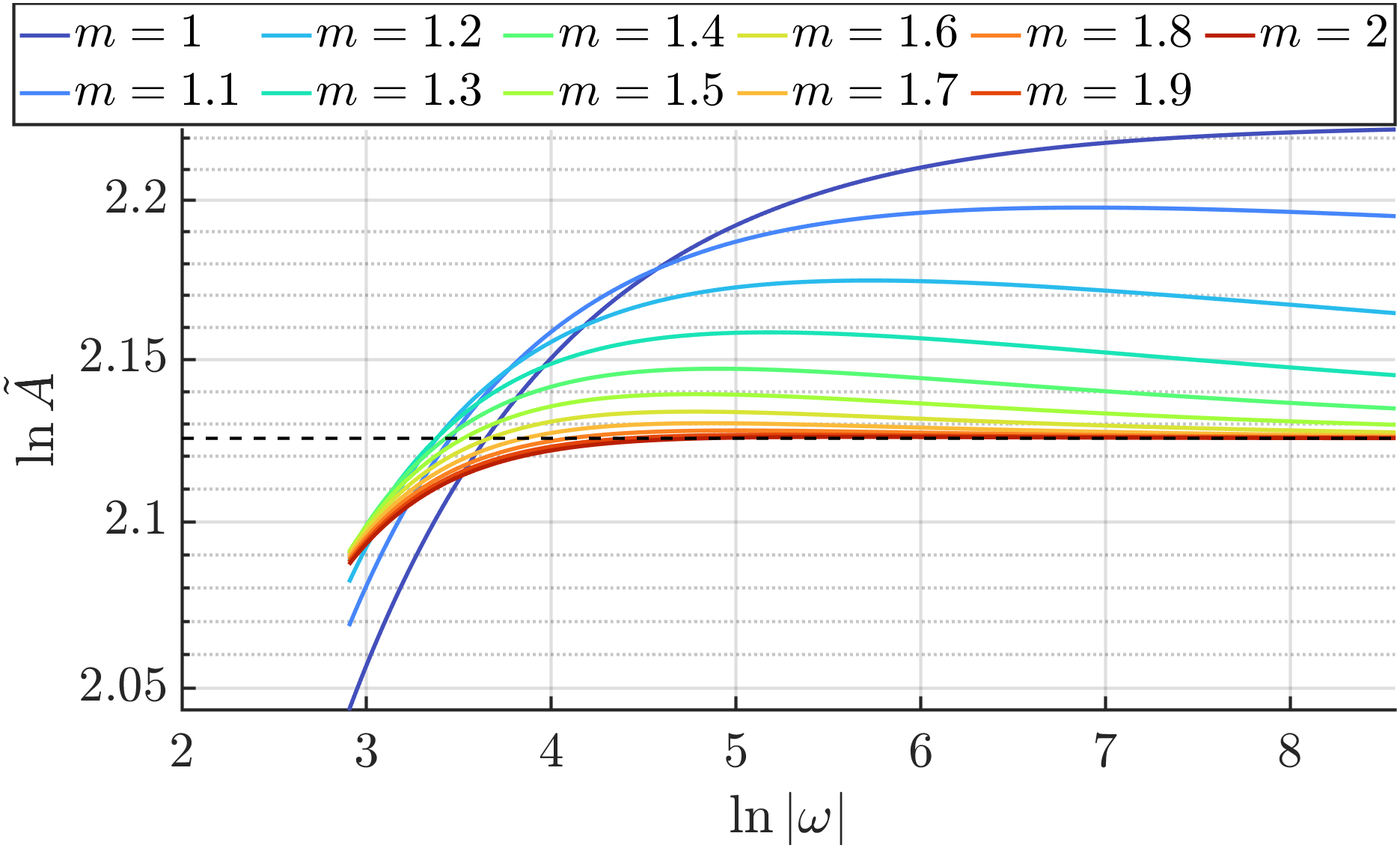}
    \caption{The Toy CP problem with $a=0$. Points $(\ln |\omega|,\tilde r)$
     (left), and $(\ln |\omega|,\ln \tilde A)$ (right) 
    obtained for  $m=1,1.1,\ldots,2$. 
 The dashed line (on the right figure) corresponds to the value $\ln (8\pi/3)=2.1255\dots$.}
    \label{fig:CasAm1a2}
\end{figure}

 However, formula \eqref{Deltap_m} 
  is an asymptotic formula, 
 so to find out the tendency of $\tilde r$ towards $3$ and $\ln A$
 towards $\ln (8\pi/3)=2.1255\dots$ (see the dashed horizontal line in Figures \ref{fig:CasAm1a6} and \ref{fig:CasAm1a2}), we proceed as above, that is, we compute, 
  for each two successive points $(\ln |w_i|,Y_{\Delta p_i})$ 
 and $(\ln |w_{i+1}|,Y_{\Delta p_{i+1}})$, 
 the segment passing through these two points, with an expression 
 $\tilde r_i\ln |w|+\ln \tilde A_i$. 
  Taking values of $K$ decreasing from $0.005$ to $5\cdot 10^{-8}$, we plot the set of points $(\ln |w_i|,\tilde r_i)$, $i=1,\dots,N-1$,  
 in Figure 
 \ref{fig:CasAm1a6} left 
and $(\ln |w_i|,\ln \tilde A_i)$, 
$i=1,\dots,N-1$, in  Figure 
 \ref{fig:CasAm1a6} right.
 These computations are shown
for  $m=1,2,3,4,5,6$.
Some remarks must be mentioned:

 (i) it is clear that the value of $\tilde r$
tends to $3$ for the values of $m$ considered so the exponent 3 in the theoretical Melnikov formula \eqref{eq:deltapmel}
fits with the numerical exponent obtained. 
Moreover we also observe that for $m=2,3,4,5,6$ (say $m\ge 2$) the tendency of $\tilde r$ is mainly the same and fast towards the value $3$, but for 
 $m=1$ it takes a bit longer to tend to $3$. But in all cases, $\tilde r$ tends to $3$ when $K$ tends to zero.
 So we may infer that for any $m\ge 1$, the exponent $\tilde r$ is 3. 
 
\begin{figure}[ht!]
    \centering
    \includegraphics[trim={0mm 0mm 0mm 0mm},clip,width=0.48\textwidth]{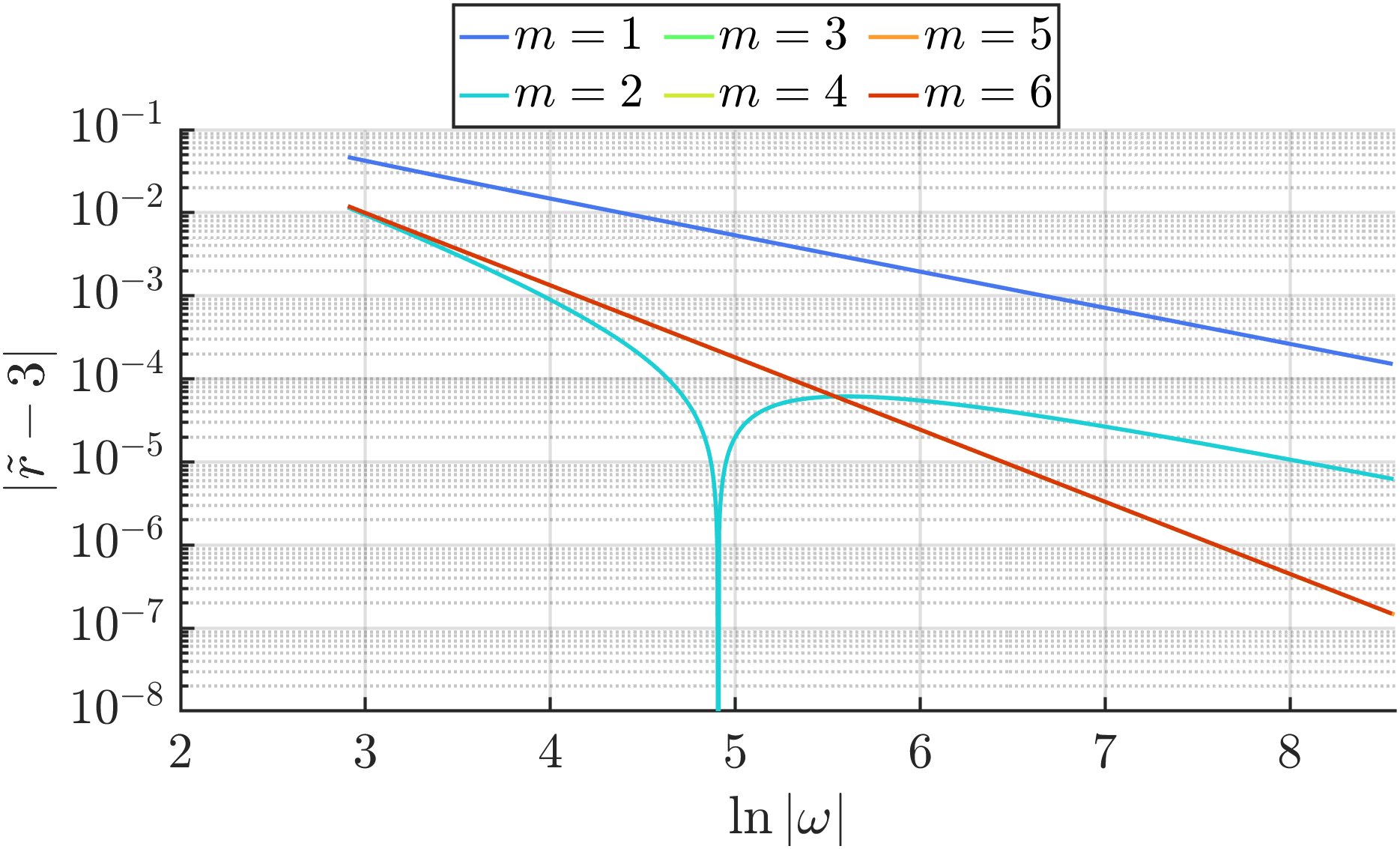}
    \includegraphics[trim={0mm 0mm 0mm 0mm},clip,width=0.48\textwidth]{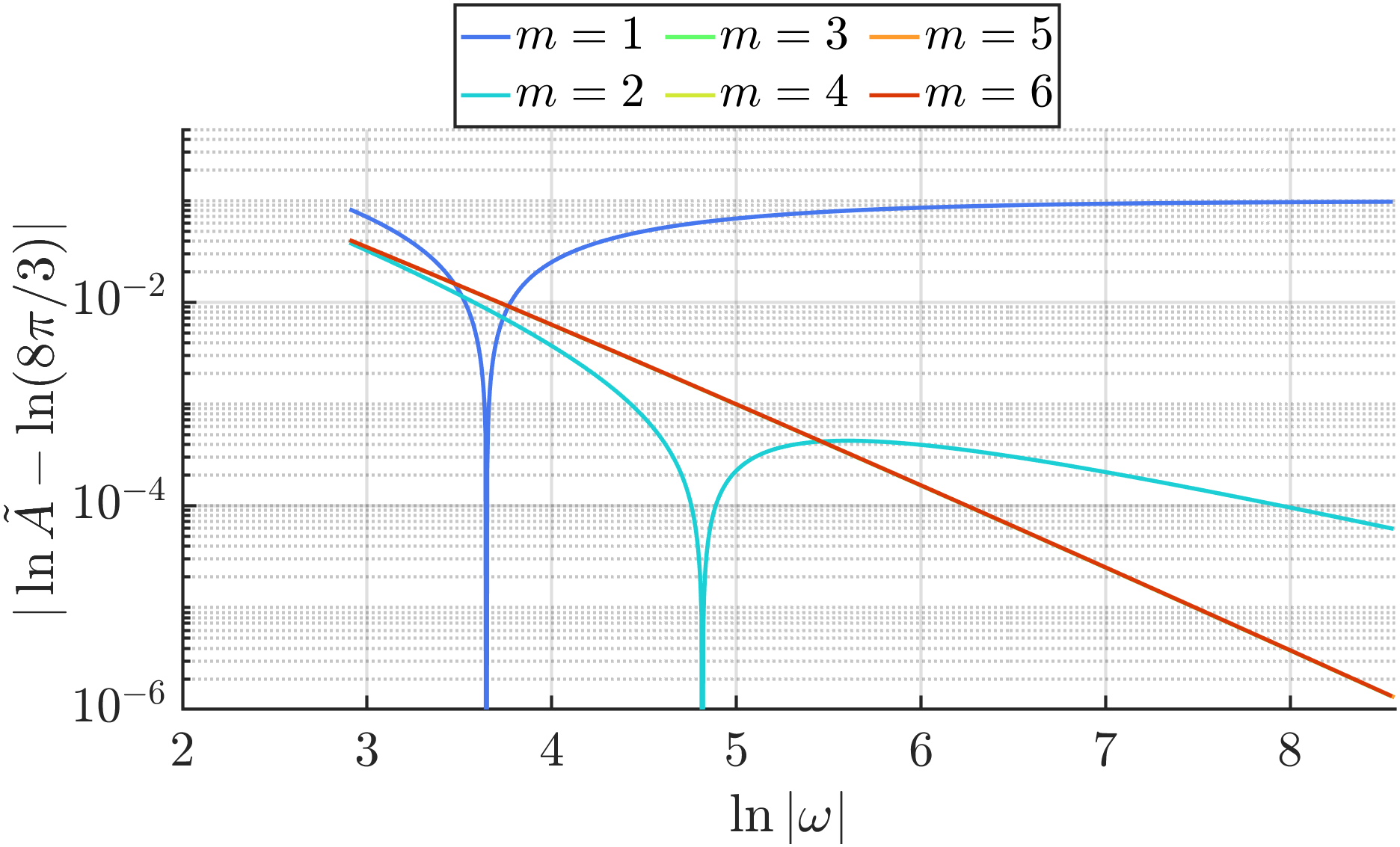}
    \caption{The Toy CP problem with $a=0$. Points $(\ln |\omega|,|\tilde r-3|)$
     (left), and $(\ln |\omega|,|\ln \tilde A-\ln (8\pi/3)|)$ (right) 
    obtained for $m=1,\ldots,6$.}
    \label{fig:CasAm1a6dif}
\end{figure}

(ii) Concerning the value of the constant $\ln \tilde A$, it is clear that, for
$m\ge 2$, $\ln \tilde A\to \ln  (8\pi/3)$. However, in principle it is not that
clear for $m=1$. Another way to show the tendency of $\tilde r_i$ and $\ln \tilde
A_i$ to $3$ and $\ln  (8\pi/3)=2.1255\dots$ respectively,  is based on the
computation of the points $(\ln |w_i|, |\tilde r_i-3|)$ in Figure
\ref{fig:CasAm1a6dif} left and the points $(\ln |w_i|, |\ln \tilde A_i-\ln
(8\pi/3)|)$ in Figure \ref{fig:CasAm1a6dif} right. We see in the left figure
that for $m=3,4,5,6$ the points {\sl overlap} on the same curve  and $|\tilde
r-3|$ tends clearly to zero. In particular, for $m=2$, the curve $\tilde r$
crosses the value $\tilde r=3$. See Figure \ref{fig:CasAm1a6} left (i.e.~the
points on the curve $(\ln |w|, |\tilde r-3|)$ describe a sharp minimum for a
value of $\ln |w|$ close to 5 and the values of $\tilde r$ tend to 3) and we
observe the decreasing tendency to zero in Figure \ref{fig:CasAm1a6dif} left).
We also see the decreasing tendency to zero for $m=1$, but not as fast as the
one for $m=2$ or $m\ge 3$.
   
Similarly, looking at the Figure \ref{fig:CasAm1a6dif} right, we observe for
$m=3,4,5,6$ the fast tendency of
\begin{displaymath}
    \left|\ln \tilde A-\ln\left(\frac{8\pi}{3}\right)\right|
\end{displaymath}
 to zero. However for $m=2$, the curve $|\ln \tilde A|$ (see Figure
 \ref{fig:CasAm1a6} right) crosses the value $\ln (8\pi/3)$ (it is clearly seen
 in the sharp minimum for $|\ln \tilde A-\ln (8\pi/3)|$ in Figure
 \ref{fig:CasAm1a6dif} right), has a local maximum  and goes on decreasing
 tending to $\ln (8\pi/3)$ (we observe a slower tendency to zero in Figure
 \ref{fig:CasAm1a6dif} right). Similarly we expect this to happen for $m=1$.
 Actually we see the crossing value $\ln (8\pi/3)$, and (almost) the maximum.
 However we cannot observe the decreasing tendency to $\ln (8\pi/3)$ as
 expected. We conjecture that this is the case.

To support this conjecture, the previous computation naturally leads to explore
the behavior of $\tilde r$ and $\ln \tilde A$ for intermediate values of $m$
between 1 and 2, to see the continuous evolution. This has been done for
$m=1.1,1.2,1.3,\dots,1.9$ and we show the results in Figure \ref{fig:CasAm1a2}
left for $\tilde r$ and right for $\ln \tilde A$. We notice the continuous
evolution of the curves for the intermediate different values of $m$. In
particular, on the left plot we see the tendency of $\tilde r$ to $3$ for any
value of $m$. On the right plot we see the maximum in each curve $\ln \tilde A$
for $m=1.1,\dots,1.9$ and the asymptotic tendency to $2.1255\dots$, as $\ln |w|$
increases, however for $m=1$ we see the maximum but not the expected decreasing
tendency to $2.1255\dots$. To do so, we should consider smaller values of $K$.
But this becomes numerically unfeasible (to compute the value of $\Delta p$ for
$K=5\cdot 10^{-8}$ we need thousands of digits and an order of thousands for the
parametrization of the manifold, so smaller values of $K$ become prohibitive).

6.~A double check concerning the value of $\tilde A$ is related to  the
extrapolation procedure done in Subsection~\ref{subsec33}. More precisely, for
every value of $m\ge 1$ fixed, we apply (finite) extrapolation steps to the
first fit formula: 
\begin{displaymath}
\frac{\Delta p|\omega| ^{-\tilde r}\exp\left(-\dfrac{\omega \pi}{2}\right)}%
{\varepsilon ^m}\sim\tilde A+\tilde A_2\varepsilon^2+\tilde A_4\varepsilon^4
 +\cdots
\end{displaymath}
or to the second fit formula
\begin{displaymath}
\frac{\Delta p|\omega| ^{-\tilde r}\exp\biggl(-\dfrac{\omega \pi}{2}\biggr)}%
{\varepsilon ^m\left(1-\dfrac{2}{\omega^2}\right)}\sim\tilde A+\bar {\tilde A}_2\varepsilon^2+\bar{\tilde A}_4\varepsilon^4+\cdots
\end{displaymath}
both formulas with $r=3$. We also want  to compare both formulas and to check if the second one
somewhat improves the first one.

Our simulations show that, as $m$ increases, the number of digits obtained for
the value of $\tilde A=8\pi/3$ also increases when taking the second fit formula
instead of the first one. For example, for $m=6$ we obtain, respectively,  6, 10
digits when taking a first step and a second step of extrapolation using the
first fit formula, whereas we obtain 20, 20 digits with two successive steps of
extrapolation applied to the second fit formula. However, as $m$ decreases, for
example $m=1$, neither the first nor the second fit formula, even using
extrapolation, are able to provide good digits of $\tilde A$. As mentioned
above, we should need smaller values of $K$, which  is not feasible.

So, from the numerical simulations, we can conclude that for the Toy CP problem
with $a=0$, Hamiltonian \eqref{pendol}, 

(i) formula \eqref{Deltap_m} provides a good fit formula for the splitting
$\Delta p$, and

(ii) the Melnikov formula predicts the splitting for any $m>1$ both in $\tilde
r$ and $\tilde A$. However for $m=1$, while the tendency of $\tilde r$ to $3$ is
clear, we conjecture this tendency to happen for $\ln \tilde A$ to $2.1255\dots$
but the numerical results are not conclusive enough.

\subsection{Other strategies/models to fit the numerical computation of 
\texorpdfstring{$\Delta p$}{Delta} for the CP problem}\label{subsec42}
Recall that our purpose is to find a simplified model that {\sl mimics} the
splitting $\Delta p$ of the CP problem in resonant coordinates. And we have just
shown that the perturbed pendulum, $H=H_0+\varepsilon H_1$, is not good enough.

So  we want  to find out another simple model that {\sl does} describe the
splitting of the CP problem.

A natural procedure is to take the second order approximation
\begin{displaymath}
    H=H_0+\varepsilon H_1+\varepsilon^2 H_2
\end{displaymath}
see \eqref{eq:H0H1H2}, and consider different truncated Hamiltonians adding to
$H_0+\varepsilon H_1$ different terms of $\varepsilon^2 H_2$. For each selected
Hamiltonian, we, first, take a value of $K$ and we compute the splitting $\Delta
p$ (obtained from the intersection of the unstable manifold of the corresponding
equilibrium point and the section $x=\pi$, $y>0$). Second, we repeat the
procedure for different decreasing values of $K$ and, third, we fit $\Delta p$
by an expression of the type
\begin{displaymath}
\Delta p\sim \varepsilon \hat A|\omega|^{\hat r}\cdot 
\exp \left(\frac{\omega\pi}{2}\right)
\end{displaymath} 
and we want to know which terms of the Hamiltonian are relevant that give
rise to a value of $\hat r$ equal to the one obtained for the CP problem, that is  $\bar r =2.111\ldots=19/9$.
We have done the numerical computations considering different truncated
Hamiltonians taking different terms of $H_2$. We have used just quadruple
precision and decreasing values of $K$ ranging from $0.001$ to $5\cdot 10^{-5}$ to have a first insight.
From these numerical simulations we can conclude that the simplest model that
reproduces the value of $\bar r=2.111\dots$~for the fitting splitting formula is
the Hamiltonian provided by $H_0+\varepsilon H_1-(2/3)\varepsilon^2y^3$, that is,
only the term $-(2/3)\varepsilon^2y^3$ in $H_2$ is responsible for changing the
exponent $\tilde r=3$ (obtained from the perturbed pendulum) to $\bar
r=2.1111\dots$.

So, in next Section we will consider this precise Hamiltonian and we will use
multiple precision computations in order to show that we really obtain the value
$\bar r=2.1111\ldots=19/9$.

\section{
The Toy CP problem with \texorpdfstring{$a=1$}{a=1}}\label{PQP} 
The previous discussion motivates to analyse the new Hamiltonian adding the term
$(-2/3)\varepsilon ^2y^3$ to the integrable part $H_0$ giving rise to a new
integrable part,
\begin{equation}\label{H0breve}
   \breve H_0(x,y,q,p)=H_0(x,y,q,p)-\frac{2}{3}\varepsilon ^2y^3,
   \end{equation} 
which we will call {\sl amended pendulum}, plus the perturbation, that is
\begin{displaymath}
\breve H_0(x,y,q,p)+\varepsilon H_1(x,y,q,p),
\end{displaymath}  
called from now on the {\sl perturbed amended pendulum}. Or, more generally, we
consider the Toy CP problem with $a=1$ and $m\ge 1$
\begin{displaymath}
\breve H(x,y,q,p)=\breve H_0(x,y,q,p)+\varepsilon^m H_1(x,y,q,p),
\end{displaymath}
In this Section we will carry out two kind of computations:
       
(i) the computation of the Melnikov integral for $\Delta p$ (see formulas
\eqref{FAw} and \eqref{deltap}). In this case, however, we cannot compute the
Melnikow integral for the amended pendulum analytically (as we did for the
pendulum) since no parametrisation of the homoclinic loop (the separatrices),
for the integrable part $\breve H_0$, is available. So the Melnikov integral
will be computed numerically. This will be done in Subsection \ref{subsec5p1}.
    
(ii) The computation of the splitting $\Delta p$ (of the unstable and stable
manifolds of the corresponding equilibrium points) at the Poincar\'e section
$\widetilde \Sigma$. This will be done in Subsection \ref{subsec5p2}.

(iii) Finally we will discuss the Melnikov prediction of the splitting for the
Toy CP problem, on the one hand (in Subsection \ref{subsec5p2}), and for the CP
problem, on the other hand (in Subsection \ref{subsec5p3}).

\subsection{Numerical computation of the Melnikov formula using the amended pendulum}\label{subsec5p1}
Recall formula \eqref{eq:v1uv1s} that provides the linear approximation of the
splitting $\Delta v_1(t)=(\Delta q_1(t), \Delta p_1(t))$. We are interested in
the splitting at $t=0$, so we use formula \eqref{deltap} for $\Delta
v_1(0)=(\Delta q_1(0),\Delta p_1(0))$. In particular, for $H_1(x,y,q,p)$ in
\eqref{eq:H1EPS0} we obtain 
\[
\Delta q_1(0)=0, \qquad 
 \Delta p_1(0)=-\frac12\int_{-\infty}^{\infty}\left(  \cos \left(2x(t)+\omega t\right)-\cos \omega t\right)\:\!\!\df t,
\]
where $(x(t),y(t))$ defines the separatrix of $\breve H_0$ and the right hand side term is precisely the {\sl Melnikov integral}.

So an strategy to compute the Melnikov integral consists of 
integrating the system of ODE 
\begin{equation}\label{eq:sistintemel}
    \left\{ 
    \begin{aligned}
        \dot{x} & = y -2\varepsilon^2y^2,\\
        \dot{y} & = \sin x,\\
        \dot{z} & = \cos(2x+\omega t)-\cos\omega t,
    \end{aligned}
    \right.
\end{equation}
and more particularly the Melnikov integral  
\[
\frac12\int _{-\infty}^{\infty} \left[\cos \left(2 x+\omega t\right)-\cos \omega t\right]\:\!\!\df t
=
\int _{-\infty}^0\left[\cos \left(2 x+\omega t\right)-\cos\omega t\right]
\:\!\!\df t
\]
due to the symmetry \eqref{simetria}. So if $\displaystyle z(t)=\int _{-\infty}^t \dot z(u) \mathrm{d}u$
we want to compute the value $z(0)$.

We proceed with the following steps: 

\textbf{Step 1:}
Applying the parameterization method (see the Appendix) we obtain a
local approximation  
of the unstable manifold (separatrix), $\bm{W}^u(\bm{0})$, associated with the equilibrium point $(0,0)$
of the system 

\begin{equation}\label{eq:separQP}
    \left\{
    \begin{aligned}
        \dot{x} & = y -2\varepsilon^2y^2\\
        \dot{y} & = \sin x
    \end{aligned}
    \right.
\end{equation}
that is, we have a (high order expansion of a) parametrization 
of the separatrix given by:
\begin{equation}\label{eq:param}
    \sum_{k=1}^N \bm{w}_k s^k + O(s^{N+1})
\end{equation}

\textbf{Step 2:} We compute a suitable value $\hat{s}$ small enough such that the error in this approximation is less than some given tolerance.

\textbf{Step 3:} We follow numerically the solution (the separatrix)
of system \eqref{eq:separQP}  from the point
\begin{displaymath}
  (x,y)= \sum_{k=1}^N \bm{w}_k \hat{s}^k  
\end{displaymath}
up to the section $\widetilde \Sigma$ ($x=\pi$) and compute the necessary time,
$T$, to reach $\widetilde \Sigma$.

\textbf{Step 4:} We define $s_0 = \cfrac{\hat{s}}{e^{-T}}$.

\textbf{Step 5:}
We compute the integral of the function $z(t)$ from $-\infty$ to $-T$. To do so,
we notice that from the Hamiltonian $\breve H_0$ in \eqref{H0breve} and taking
into account the value $H=h=0$ at $(0,0)$,
\[
    \cos x= 1-\frac{y^2}{2}+ \frac{2}{3}\varepsilon^2y^3
\]
and we have
\begin{equation}\label{eq:zpunt}
\dot{z}  = \cos(2x+\omega t) -\cos\omega t
=-2\dot{y}^2\cos\omega t +\dot{y}\left(-2+y^2
-\frac{4}{3}\varepsilon^2y^3\right)\sin\omega t.
\end{equation}

We notice that from \eqref{eq:param} we have, in particular, the
expansion of $y(t)$: 
\[    y(t) = y_1 s_0e^t+ y_2 s_0^2e^{2t}+y_3 s_0^3e^{3t}
           + \cdots= y_1s + 2y_2s^2 + 3y_3s^3 + \cdots
\]

Thus, the expression of $\dot z$ in \eqref{eq:zpunt}
becomes an expansion in $s$ (or equivalently in $s_0e^t$)
together with the terms $\sin(\omega t)$ and $\cos(\omega t)$.
We just need to integrate $\dot z(t)$ from $-\infty$ to $-T$.
 Let us provide a formula for the appearing terms in the integral:
\begin{align*}
    \int_{-\infty}^{-T}a_ks^k\cos\omega t\:\! \df t 
                & = a_k\hat{s}^k\frac{k\cos(-\omega T) + \omega \sin(-\omega T)}{k^2+\omega^2},\\[8pt]
    \int_{-\infty}^{-T}a_ks^k\sin\omega t\:\! \df t 
          & = a_k\hat{s}^k\frac{k\sin(-\omega T) - \omega \cos(-\omega T)}{k^2+\omega^2}.
\end{align*}

\textbf{Step 6:}
Once we get the value $z(-T)=\displaystyle\int _{-\infty}^{-T}\!\!\dot z(t)\df
t$ we integrate numerically system \eqref{eq:sistintemel} from $-T$ to $0$, and
we obtain $z(0)$.

\begin{figure}[ht!]
    \centering
    \includegraphics[trim={0mm 0mm 0mm 0mm},clip,width=0.48\textwidth]{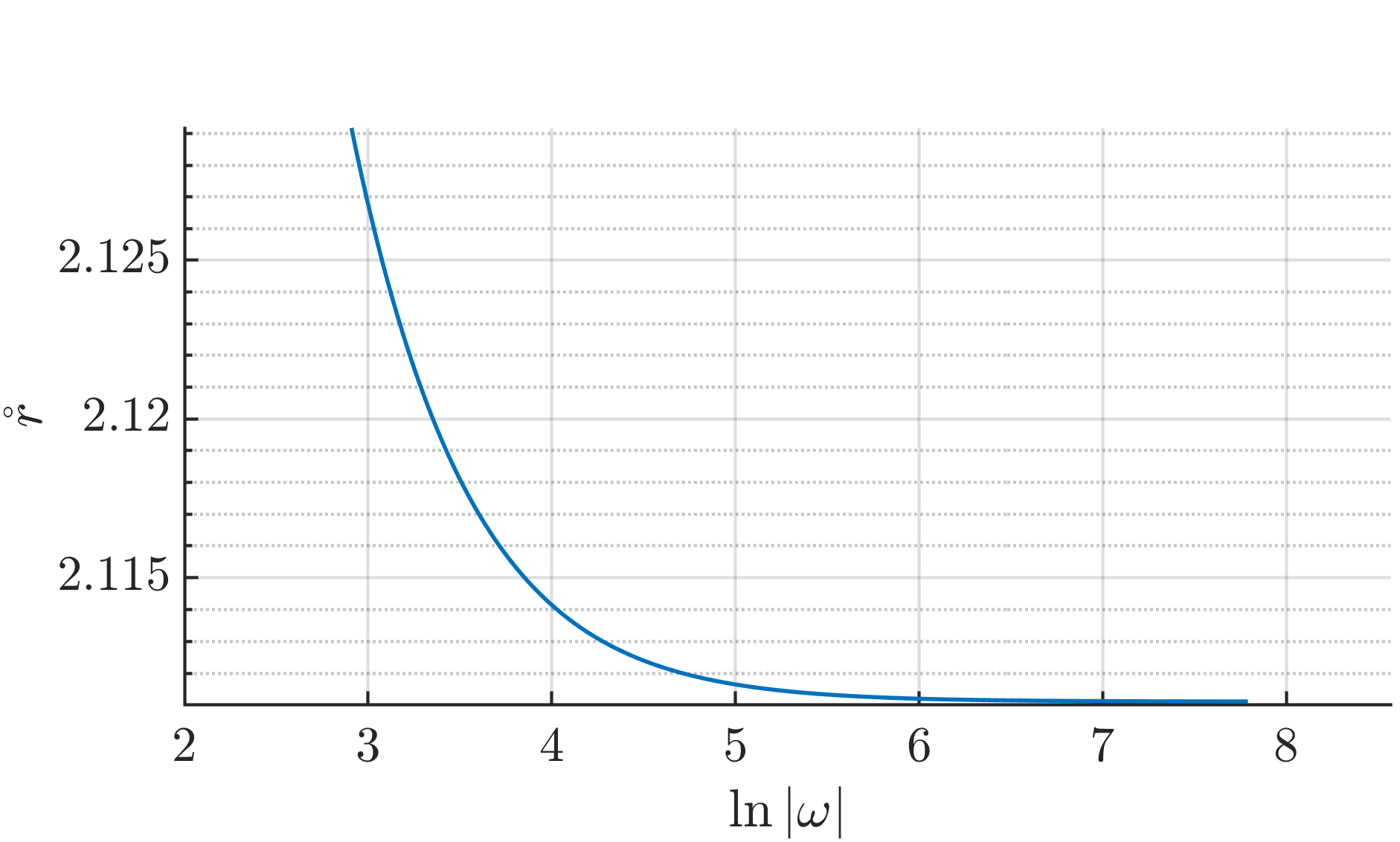}
    \includegraphics[trim={0mm 0mm 0mm 0mm},clip,width=0.48\textwidth]{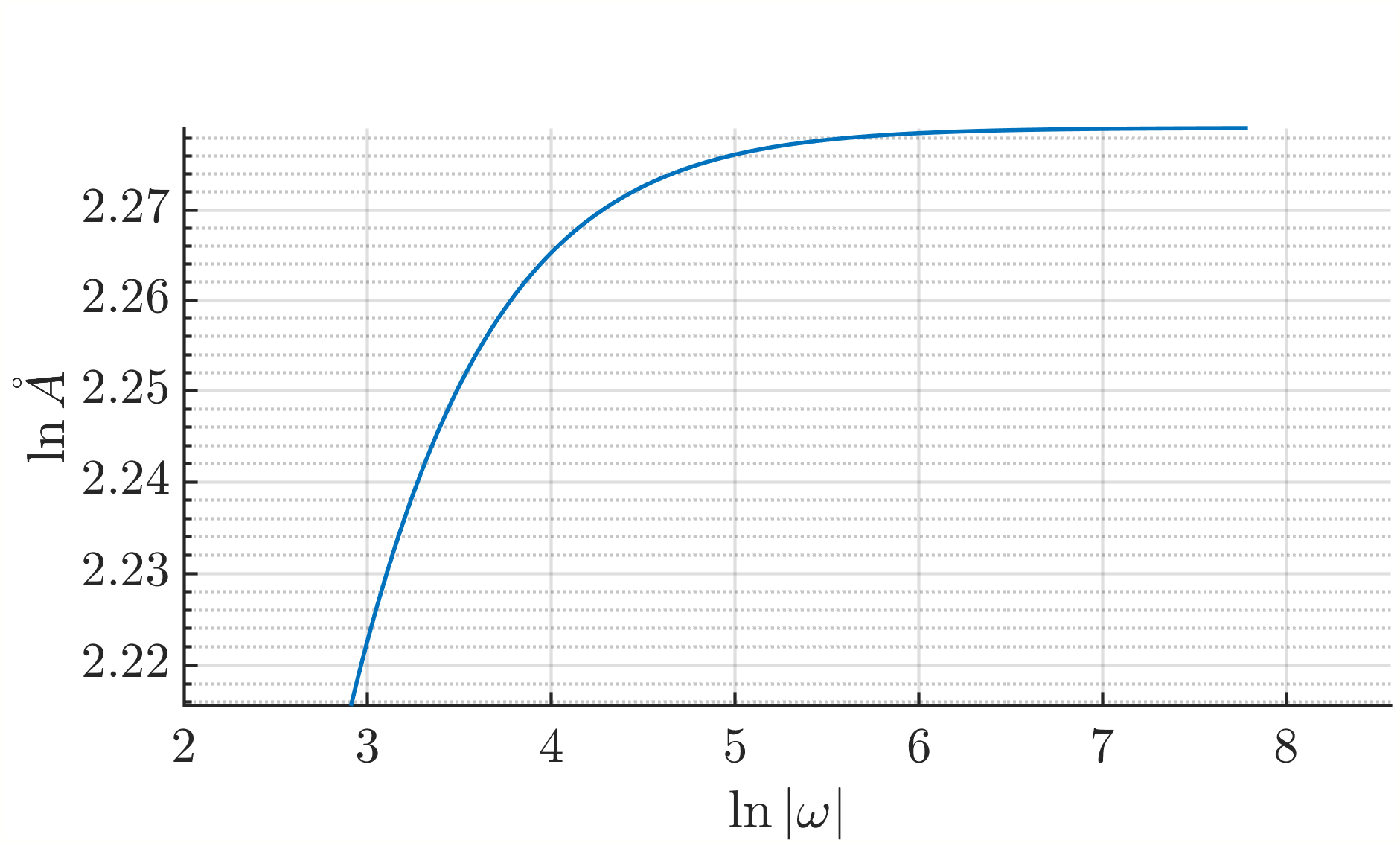} 
    \caption{Curves $(\ln |\omega|,\mathring r)$ and $(\ln |\omega|,\ln \mathring A)$ 
    obtained from the computation of the Melnikov integral.
    }
    \label{fig:intemel_mval1}
\end{figure}

Taking different decreasing values of $K$ from $K=0.001$ to $K=5.7\cdot 10^{-8}$, we fit the output values for the Melnikov integral
  by the formula
\begin{equation}\label{zde0}
z(0)\sim \mathring A|\omega|^{\mathring r}\cdot \exp\left(\frac{\omega\pi}{2}\right)
\end{equation}
or
\[
\ln z(0) -\frac{\omega\pi}{2}
\sim  \mathring r\ln |\omega|+\ln \mathring A ,   
\]
and we obtain (see Figure \ref{fig:intemel_mval1})
 \[
  \mathring r=2.111\dots,\qquad \qquad \ln \mathring A=2.279\dots .
  \]
We remark the good coincidence of the value of $\mathring r$ with that  obtained for
the CP problem in resonant coordinates $\bar r$.

Concerning the value of  $\ln \mathring A$, there is not such coincidence with
that obtained for the CP problem in resonant coordinates, $\ln \bar A$.

\subsection{Straight numerical computation of the splitting 
\texorpdfstring{$\Delta p$}{Delta(p)}}\label{subsec5p2}  
Similarly as we proceeded for the perturbed pendulum, we consider the Toy CP
problem with $a=1$, that is,
\begin{displaymath}
\breve H(x,y,q,p)=H_0(x,y,q,p)- \frac{2}{3}\varepsilon ^2y^3 +\varepsilon ^m H_1(x,y,q,p).
\end{displaymath}  
We fix a value of $m\ge 1$, we take different values of $K$ (decreasing up to
$3\cdot 10^{-8}$), and for each $K$ given we compute the unstable manifold of the
equilibrium point $L_-=(0,0,3\varepsilon ^{2+m},0)$  intersection with $\widetilde
\Sigma$ and we focus on the splitting $\Delta p=2p^u$, which we want to fit by 

\begin{equation}\label{fitdeltapbreve}
 \Delta p\sim \varepsilon |\omega|^{\breve r}\exp 
 \left(\frac{\omega\pi}{2}\right),
\end{equation}
or equivalently
\begin{equation}\label{fitdeltapQP_rodona}
\ln \Delta p-m\ln \varepsilon -\frac{\omega\pi}{2}
\sim \breve r\ln |\omega|+\ln \breve  A    
\end{equation}
In this case, unlike the perturbed pendulum problem, we do not know neither the value of $\breve r$ nor the value of $\ln \breve A$ in formula
\eqref{fitdeltapQP_rodona}.

\begin{figure}[ht!]
    \centering
    \includegraphics[trim={0mm 0mm 0mm 0mm},clip,width=0.48\textwidth]{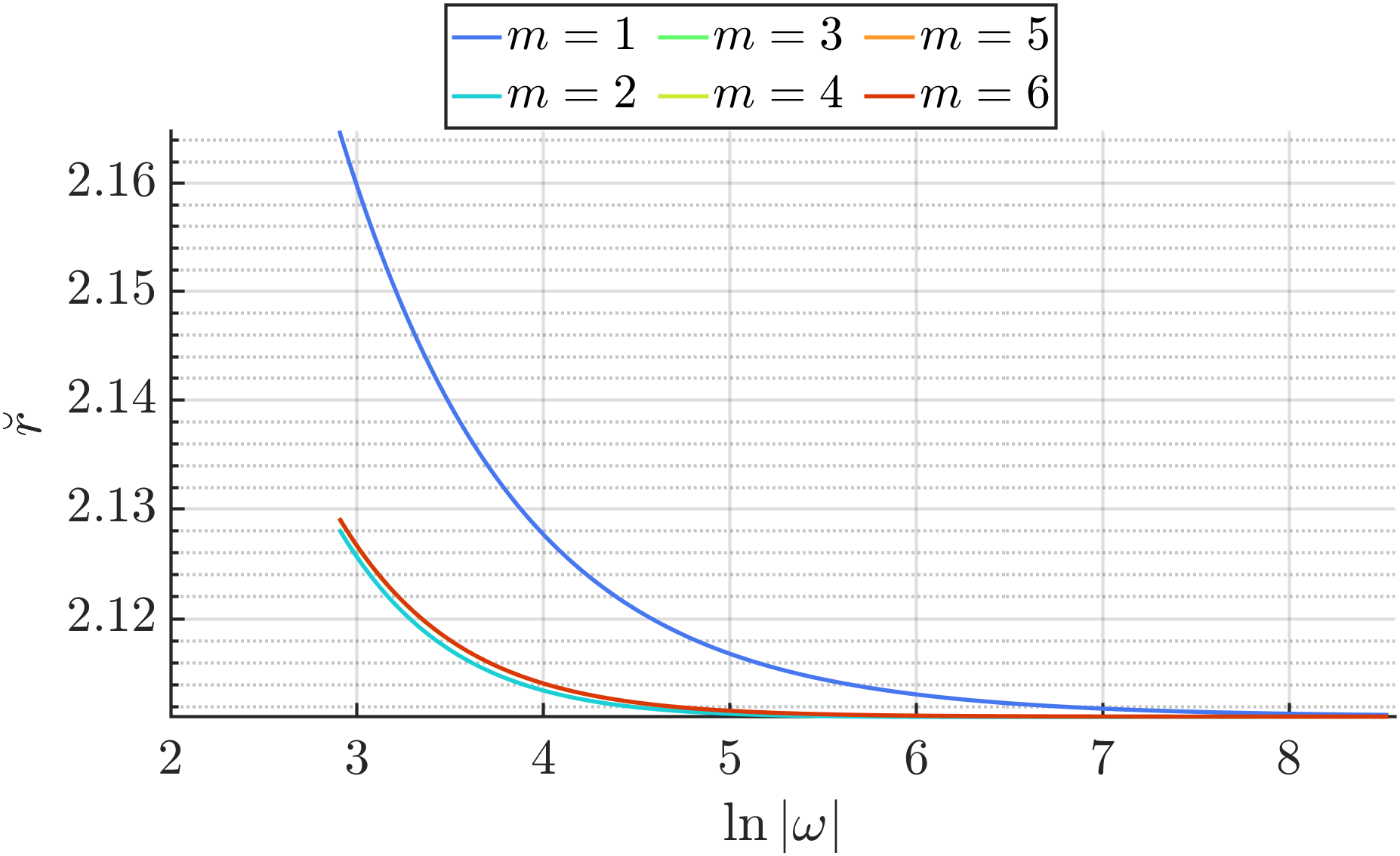}
    \includegraphics[trim={0mm 0mm 0mm 0mm},clip,width=0.48\textwidth]{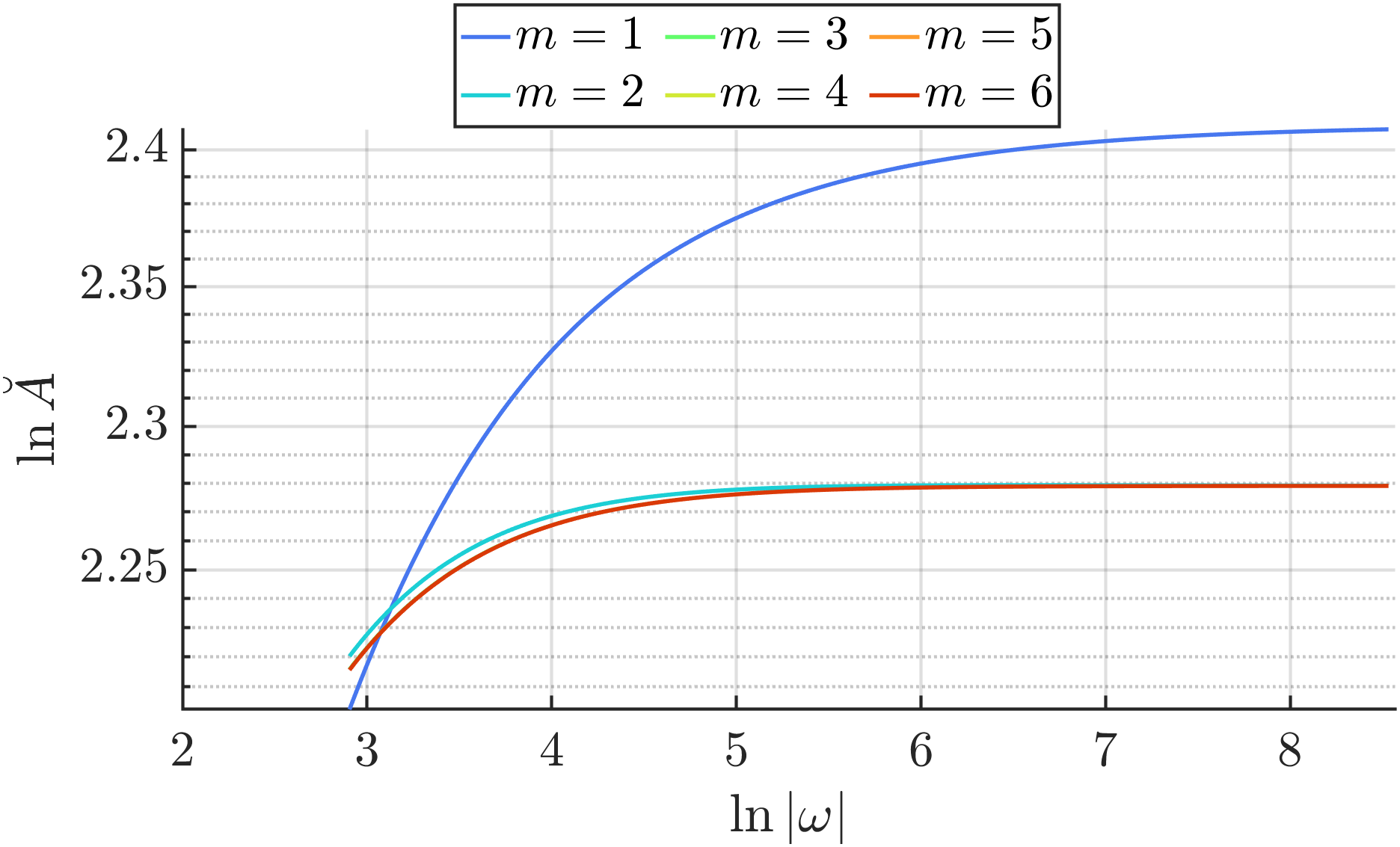}
    \caption{The Toy CP problem with $a=1$. Points $(\ln |\omega|,\breve r)$
     (left), and $(\ln |\omega|,\ln \breve A)$ (right) 
    obtained for $m=1,\ldots,6$.}
    \label{fig:CasCm1a6}
\end{figure}

We plot in Figure~\ref{fig:CasCm1a6} left the obtained points $(\ln
|\omega_i|,\breve r_i)$ (similarly as we did above, that is, taking a segment
between two successive points $(\ln |\omega_i|,Y_{\Delta p_i})$, $(\ln
|\omega_{i+1}|,Y_{\Delta p_{i+1}})$, with $Y_{\Delta p}:=\ln \Delta p-m\ln \varepsilon
-\omega\pi/2$). We observe the tendency of $\breve r$ towards the value
$2.111115\dots$, which coincides with the value 
of $\mathring r$ obtained from the Melnikov integral. 

\begin{figure}[ht!]
    \centering
    \includegraphics[trim={0mm 0mm 0mm 0mm},clip,width=0.48\textwidth]{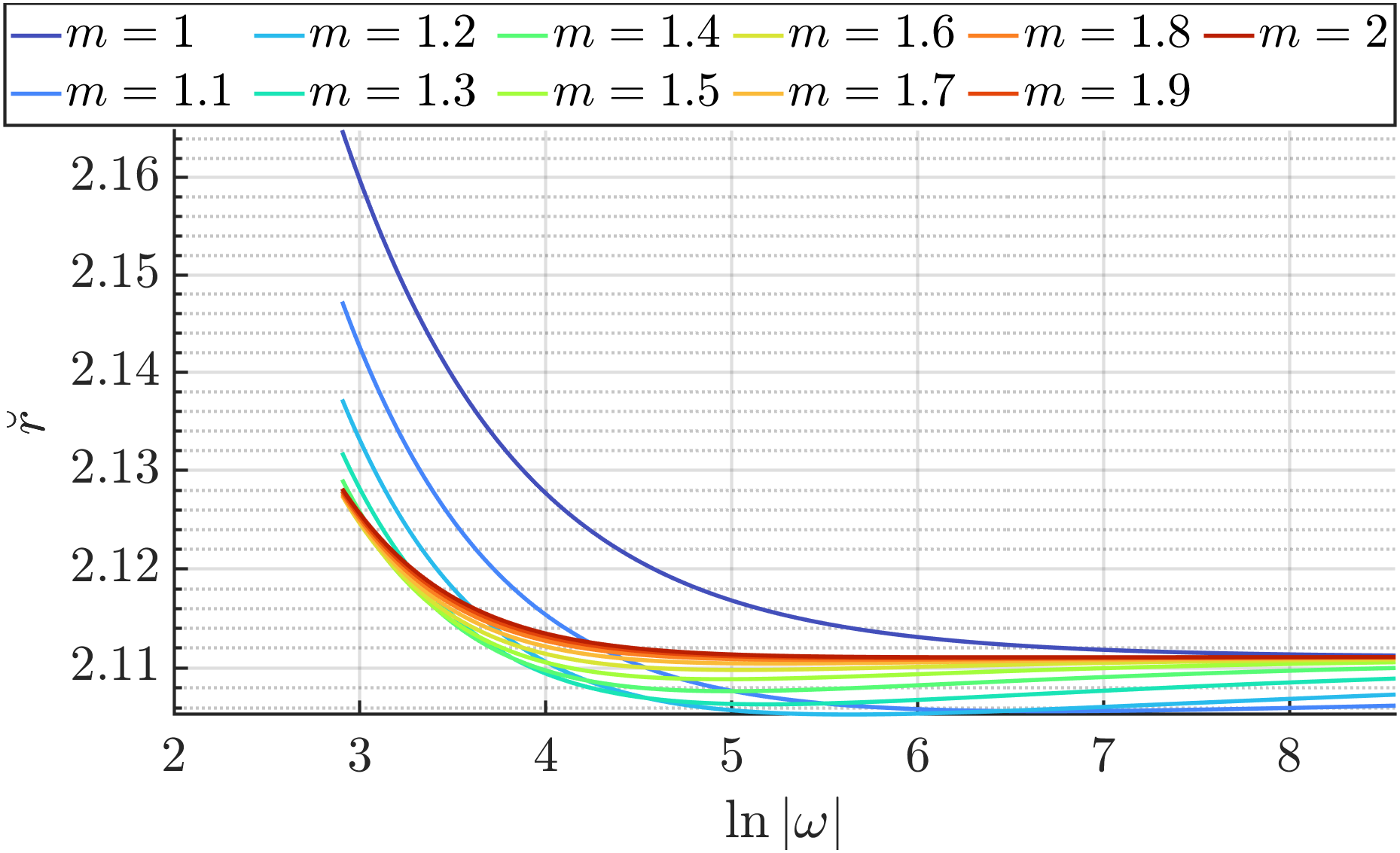}
    \includegraphics[trim={0mm 0mm 0mm 0mm},clip,width=0.48\textwidth]{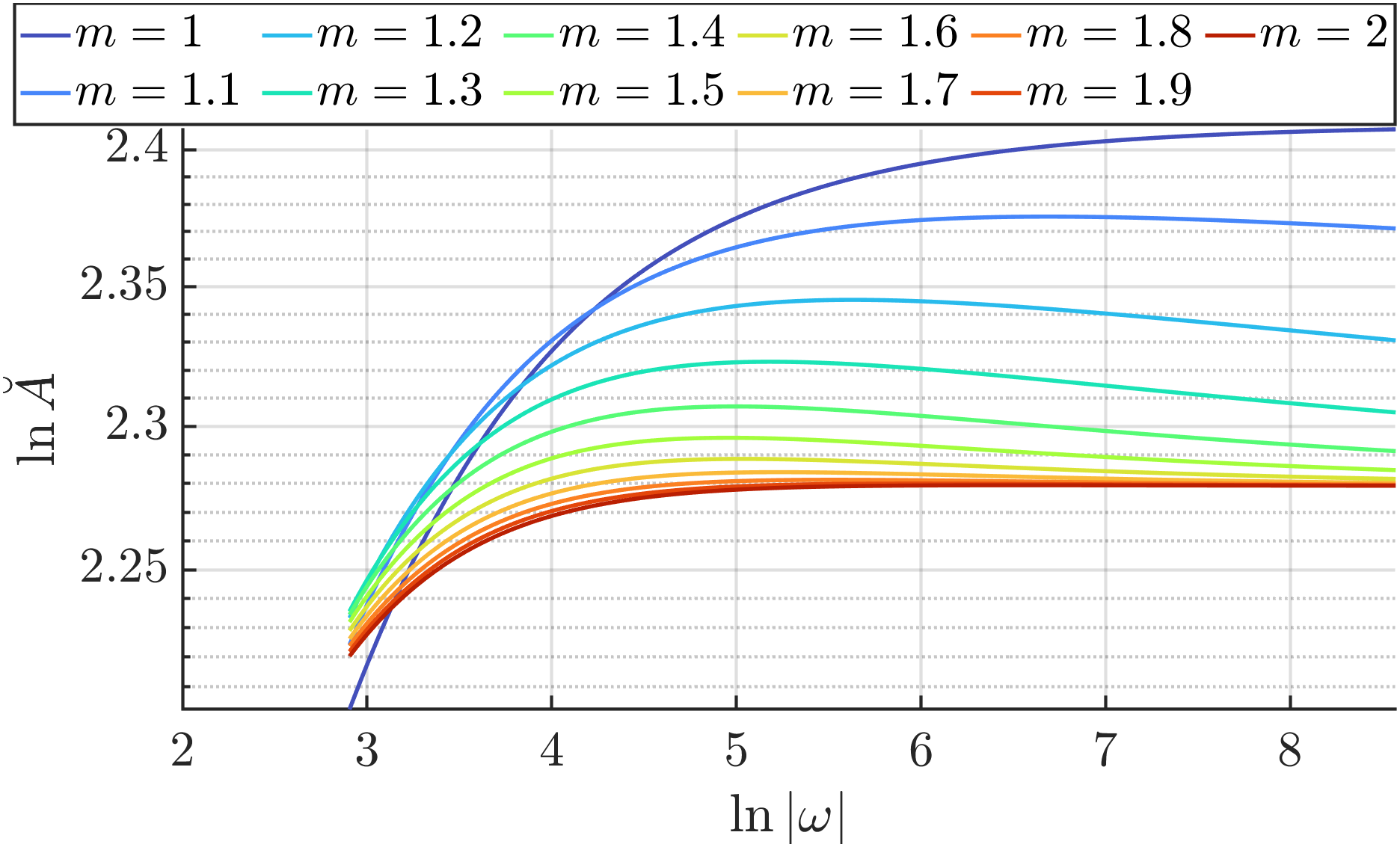}
    \caption{The Toy CP problem with $a=1$. Points $(\ln |\omega|,\breve r)$
     (left), and $(\ln |\omega|,\ln \breve A)$ (right) 
    obtained for $m=1,1.1,1.2,\ldots,2$}
    \label{fig:CasCm1a2}
\end{figure}

Regarding the  value of $\ln \breve A$, it is apparently clear
that for $m\ge 2$, $\ln \breve A$ tends to $2.279\dots$,  which coincides with the value 
of $\ln \mathring A$ obtained from the Melnikov integral. However it is not that
clear for $m=1$. Analogously as we did for the perturbed pendulum, we now
explore the intermediate values $m=1.1,1.2,1.3,1.4,\dots,1.9$ between $m=1$ and
$m=2$. The results are shown in Figure \ref{fig:CasCm1a2} left --for $\breve
r$-- and right --for $\ln \breve A$--. The continuous evolution is clear. So we
would expect,  for $m=1$, a tendency of $\ln \breve A$ to $2.279\dots$. But again,
as in the perturbed pendulum case, for $m=1$, we need smaller values of $K$
(which turns out to be prohibitive from a numerical point of view).

So, from the numerical computations done, we can conclude that for the Toy CP problem with $a=1$,

(i) the formula \eqref{fitdeltapQP_rodona}
provides a good fitting for $\Delta p$  
with $\breve r=2.111\ldots=19/9$, for $m\ge 1$, and $\ln \breve A=2.279\dots$ for $m>1$, and we expect/conjecture also the same asymptotic value of $\ln \breve A$ for $m=1$. 

(ii)  The Melnikov integral {\sl does} predict the value of the exponent $\mathring r=\breve r=19/9$
for $m\ge 1$ and the value of $\mathring A=\breve A=2.279\dots$ for $m>1$, and presumably also for
$m=1$.

\subsection{Agreement between the real CP problem and the Toy CP problem with
\texorpdfstring{$a=m=1$}{a = m = 1}}\label{subsec5p3}  

(i) Concerning the computation of the splitting  $\Delta p$ (from the unstable/stable manifolds   at the intersection with $x=\pi$, $y>0$), we now consider
 two different problems:
the original CP problem in resonant coordinates and the perturbed amended pendulum
 with $m=1$, that is, the Toy CP problem with $a=m=1$. 
 We want to compare the respective values of $\bar r$
and $\breve r$ in the corresponding fit formulas
  \eqref{deltapu} and \eqref{fitdeltapbreve}.
 
\begin{figure}[ht!]
    \centering
    \includegraphics[trim={0mm 0mm 0mm 0mm},clip,width=0.48\textwidth]{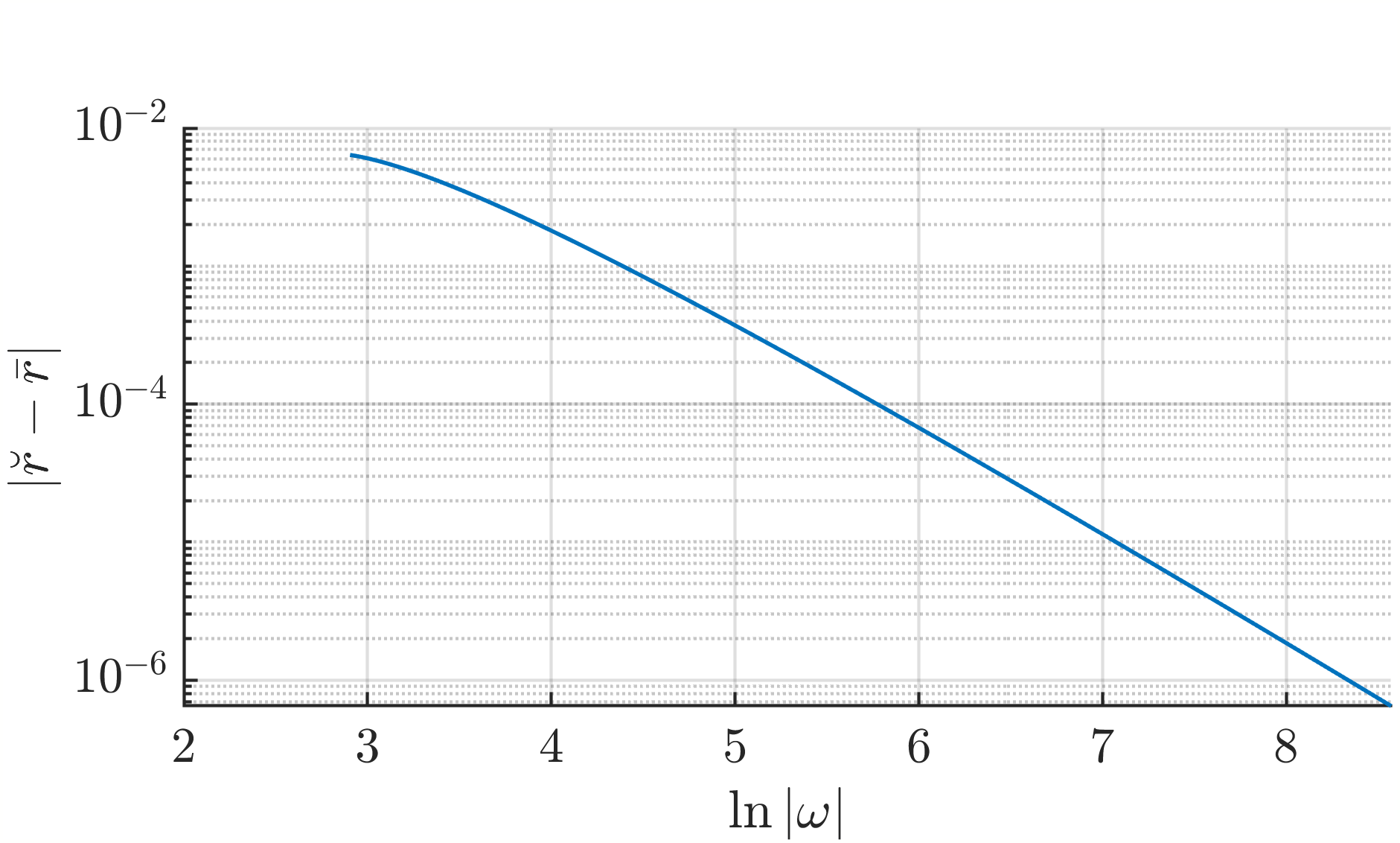}
    \caption{Curve $(\ln |\omega|,|\breve r-\bar r|)$  when computing the straight splitting from the CP problem and the perturbed amended pendulum (with $m=1$). See the text for details.}
    \label{fig:difrealvsCasCm1}
\end{figure}

More precisely, concerning the values of $\bar r$ and $\breve r$,   
we want to compare two curves, the numerical obtained curve,
 $(\ln |\omega|,\bar r)$,
for the splitting in resonant coordinates
of the original CP Hamiltonian, and the 
numerical obtained curve, $(\ln |\omega|,\breve r)$,
for the splitting of the perturbed amended pendulum Hamiltonian (with $m=1$). In Figure
 \ref{fig:difrealvsCasCm1} left we plot the curve
 $(\ln |\omega|,|\breve r-\bar r|)$ to find out the
 differences between both values. We remark the 
 very good agreement of both values as far as
 $K$ decreases (or equivalently $\ln |\omega|$  increases).

So we can conclude that the perturbed amended pendulum is a simplified model that already describes the splitting of the CP problem, as
far as asymptotic formulas (when $K$ tends to zero) are concerned. 

(ii) Finally, how does the Melnikov integral 
 (taking into account the CP problem with $a=m=1$) predict the splitting for the
  CP problem in resonant coordinates?
  That is, we want to compare the fit formulas for $\Delta p$ (of the CP problem, that is \eqref{deltapu})
  and for the Melnikov integral (of the perturbed amended pendulum problem with $m=1$, that is \eqref{zde0}),  more precisely, the values of the exponents $\bar r$ and
$\mathring r$. 

\begin{figure}[ht!]
    \centering
    \includegraphics[trim={0mm 0mm 0mm 0mm},clip,width=0.48\textwidth]{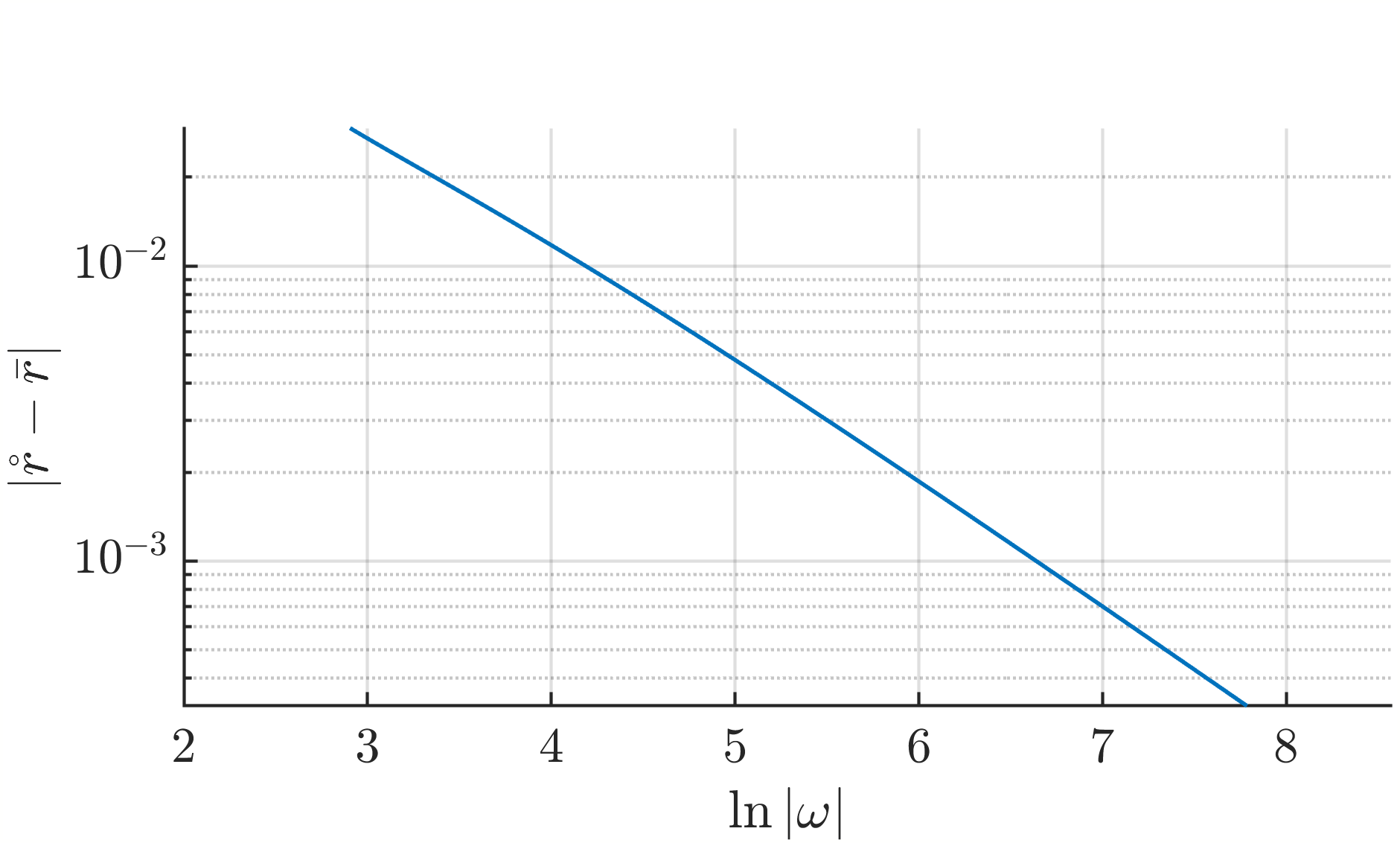}
    \caption{
Curve 
$(\ln |\omega|,|\mathring r - \bar r|)$  
when computing the straight splitting from the CP problem and the splitting using the Melnikov integral of the perturbed amended pendulum (with $m=1$). See the text for details. 
}
    \label{comparativaCPiMEL}
\end{figure}

In Figure \ref{comparativaCPiMEL} 
 we observe the difference $|\mathring r - \bar r|$
 tending to zero.
 So we can conclude that the
fitting asymptotic formula obtained for the Melnikov integral provides
the same limit value of the exponent $\bar r=\mathring r=2.111\ldots=19/9$.

{\bf Remark}. In Section \ref{secpendol}, where we considered the perturbed pendulum, we discussed the known results (in the literature) compared with the ones obtained in this paper. However, for the Toy CP problem we are not aware of any
published results.

\section*{Acknowledgments}
\addcontentsline{toc}{section}{\numberline{}Acknowledgements}%
This work is supported by the Spanish State Research Agency, through the Severo
Ochoa and María de Maeztu Program for Centers and Units of Excellence in R\&D
(CEX2020-001084-M). The authors were also supported by the Spanish grant
PID2021-123968NB-I00 (MICIU/AEI/10.13039/501100011033/FE\-DER/UE).
\noindent

\bibliographystyle{plain}
\bibliography{ChaosCPBiblio}

\def\cprime{$'$} \def \registered {$^{{\ooalign{\hfil\raise.07ex \hbox{\footnotesize R}\hfil\crcr \mathhexbox20D}}}$}\def \trademark {$^{\hbox{\sc tm}}$}
\begin{thebibliography}{10}

\bibitem{Baldoma2006}
I.~Baldom\'{a}.
\newblock The inner equation for one and a half degrees of freedom rapidly forced {H}amiltonian systems.
\newblock {\em Nonlinearity}, 19(6):1415--1445, 2006.

\bibitem{BaldomaGG2021Arxiv-I}
I.~Baldom\'{a}, M.~Giralt, and M.~Guardia.
\newblock Breakdown of homoclinic orbits to ${L}_{3}$ in the {RPC3BP} ({I}). {C}omplex singularities and the inner equation.
\newblock {\em Adv. Math.}, 408(6), 2022.

\bibitem{BaldomaGG2021Arxiv-II}
I.~Baldom\'{a}, M.~Giralt, and M.~Guardia.
\newblock Breakdown of homoclinic orbits to ${L}_{3}$ in the {RPC3BP} ({II}). {A}n asymptotic formula.
\newblock {\em Adv. Math.}, 430(109218), 2023.

\bibitem{BaldomaM2012}
I.~Baldom\'a and P.~Mart\'in.
\newblock The inner equation for generalized standard maps.
\newblock {\em SIAM J. Appl. Dyn. Syst.}, 11(3):1062--1097, 2012.

\bibitem{BaldomaFGS2012}
Inmaculada Baldom\'{a}, Ernest Fontich, Marcel Guardia, and Tere~M. Seara.
\newblock Exponentially small splitting of separatrices beyond {M}elnikov analysis: rigorous results.
\newblock {\em J. Diff. Eqs.}, 253(12):3304--3439, 2012.

\bibitem{BarrabesOBFM2012}
E.~Barrab\'{e}s, M.~Oll\'{e}, F.~Borondo, D.~Farrelly, and J.M. Mondelo.
\newblock Phase space structure of the hydrogen atom in a circularly polarized microwave field.
\newblock {\em Phys. D}, 241(4):333--349, 2012.

\bibitem{BOR2024}
E.~Barrab\'{e}s, M.~Oll\'{e}, and \'O. Rodr\'{\i}guez.
\newblock On the collision dynamics in a molecular model.
\newblock {\em Phys. D}, 460:1--31, 2024.

\bibitem{PhysRevA.90.063402}
Jaros\l{}aw~H. Bauer, Francisca Mota-Furtado, Patrick~F. O'Mahony, Bernard Piraux, and Krzysztof Warda.
\newblock Ionization and excitation of the excited hydrogen atom in strong circularly polarized laser fields.
\newblock {\em Phys. Rev. A}, 90:063402, December 2014.

\bibitem{Bialynicki-BirulaKE1994}
I.~Bialynicki-Birula, M.~Kalinski, and J.H. Eberly.
\newblock Lagrange equilibrium points in celestial mechanics and nonspreading wave packets for strongly driven {R}ydberg electrons.
\newblock {\em Phys. Rev. Lett.}, 73(13):1777--1780, 1994.

\bibitem{BrunelloUF1996}
A.F. Brunello, T.~Uzer, and D.~Farrelly.
\newblock Hydrogen atom in circularly polarized microwaves: {C}haotic ionization via core scattering.
\newblock {\em Phys. Rev. A}, 55(5):3730--3745, 1997.

\bibitem{BuchleitnerDG1995}
A.~Buchleitner, D.~Dominique, and J.-C. Gay.
\newblock Microwave ionization of three-dimensional hydrogen atoms in a realistic numerical experiment.
\newblock {\em Journal of the Optical Society of America. B, Optical Physics}, 12(4):505--519, 1995.

\bibitem{param1}
X.~Cabr{é}, E.~Fontich, and R.~de~la Llave.
\newblock The parameterization method for invariant manifolds {I}: Manifolds associated to non-resonant subspaces.
\newblock {\em Indiana University Mathematics Journal}, 52, 2003.

\bibitem{param2}
X.~Cabr{é}, E.~Fontich, and R.~de~la Llave.
\newblock The parameterization method for invariant manifolds {I}{I}: regularity with respect to parameters.
\newblock {\em Indiana University Mathematics Journal}, 52, 2003.

\bibitem{param3}
X.~Cabr{é}, E.~Fontich, and R.~de~la Llave.
\newblock The parameterization method for invariant manifolds {I}{I}{I}: overview and applications.
\newblock {\em J. Diff. Eqs.}, 218, 2005.

\bibitem{DelshamsG2000}
A.~Delshams and P.~Guti{\'e}rrez.
\newblock Splitting potential and the {P}oincar\'e-{M}elnikov me\-thod for whis\-ke\-red to\-ri in {H}amiltonian systems.
\newblock {\em J. Nonlinear Sci.}, 10(4):433--476, 2000.

\bibitem{PhysRevLett.74.1720}
David Farrelly and T.~Uzer.
\newblock Ionization mechanism of {R}ydberg atoms in a circularly polarized microwave field.
\newblock {\em Phys. Rev. Lett.}, 74:1720--1723, March 1995.

\bibitem{10.1145/1236463.1236468}
Laurent Fousse, Guillaume Hanrot, Vincent Lef\`{e}vre, Patrick P\'{e}lissier, and Paul Zimmermann.
\newblock Mpfr: A multiple-precision binary floating-point library with correct rounding.
\newblock {\em ACM Trans. Math. Softw.}, 33(2):13–es, June 2007.

\bibitem{PhysRevLett.64.511}
Panmimg Fu, T.~J. Scholz, J.~M. Hettema, and T.~F. Gallagher.
\newblock Ionization of rydberg atoms by a circularly polarized microwave field.
\newblock {\em Phys. Rev. Lett.}, 64:511--514, January 1990.

\bibitem{Ganesan}
K.~Ganesan and R.~GiBarowski.
\newblock Chaos in the hydrogen atom interacting with external fields.
\newblock {\em Pramana J. of Physics}, 48(2):379--410, 1997.

\bibitem{PhysRevA.51.1508}
Robert Ge\ifmmode~\mbox{\c{}}\else \c{}\fi{}barowski and Jakub Zakrzewski.
\newblock Ionization of hydrogen atoms by circularly polarized microwaves.
\newblock {\em Phys. Rev. A}, 51:1508--1519, February 1995.

\bibitem{GelfreichS2008}
Vassili Gelfreich and Carles Sim\'{o}.
\newblock High-precision computations of divergent asymptotic series and homoclinic phenomena.
\newblock {\em Discr. Contin. Dyn. Syst. Ser. B}, 10(2-3):511--536, 2008.

\bibitem{AutoDif}
Andreas Griewank and Andrea Walther.
\newblock {\em Evaluating Derivatives}.
\newblock Society for Industrial and Applied Mathematics, second edition, 2008.

\bibitem{PhysRevA.45.R2678}
Jennifer~A. Griffiths and David Farrelly.
\newblock Ionization of rydberg atoms by circularly and elliptically polarized microwave fields.
\newblock {\em Phys. Rev. A}, 45:R2678--R2681, Mar 1992.

\bibitem{Guardia13}
Marcel Guardia.
\newblock Splitting of separatrices in the resonances of nearly integrable {H}amiltonian systems of one and a half degrees of freedom.
\newblock {\em Discr. Contin. Dyn. Syst.}, 33(7):2829--2859, 2013.

\bibitem{GuardiaOS2010}
Marcel Guardia, Carme Oliv\'{e}, and Tere~M. Seara.
\newblock Exponentially small splitting for the pendulum: a classical problem revisited.
\newblock {\em J. Nonlinear Sci.}, 20(5):595--685, 2010.

\bibitem{Guck}
J.~Guckenheimer and P.~Holmes.
\newblock {\em Nonlinear oscillations, dynamical systems and bifurcations of vector fields}.
\newblock Springer-Verlag, 1983.

\bibitem{paramLlibre}
{À}. Haro, M.~Canadell, JL. Figueras, A.~Luque, and J.M. Mondelo.
\newblock {\em The Parameterization Method for Invariant Manifolds: From Rigorous Results to Effective Computations}.
\newblock Applied Mathematical Sciences 195. Springer International Publishing, 1 edition, 2016.

\bibitem{PhysRevLett.84.610}
Charles Jaff\'e, D.~Farrelly, and T.~Uzer.
\newblock Transition state theory without time-reversal symmetry: Chaotic ionization of the hydrogen atom.
\newblock {\em Phys. Rev. Lett.}, 84:610--613, Jan 2000.

\bibitem{JorbaZ2005}
\`A. Jorba and M.~Zou.
\newblock A software package for the numerical integration of {ODE}s by means of high-order {T}aylor methods.
\newblock {\em Experiment. Math.}, 14(1):99--117, 2005.

\bibitem{KoltsovaLDG05}
Oksana Koltsova, Lev Lerman, Amadeu Delshams, and Pere Guti\'errez.
\newblock Homoclinic orbits to invariant tori near a homoclinic orbit to center-center-saddle equilibrium.
\newblock {\em Phys. D}, 201(3-4):268--290, 2005.

\bibitem{O2018}
Merc\`e Oll\'e.
\newblock To and fro motion for the hydrogen atom in a circularly polarized microwave field.
\newblock {\em Comm. Nonlinear Sci. Numer. Simulat.}, 54:286--301, 2018.

\bibitem{OP2018}
Merc{\`{e}} Oll{\'{e}} and Juan~Ram{\'{o}}n Pacha.
\newblock Hopf bifurcation for the hydrogen atom in a circularly polarized microwave field.
\newblock {\em Comm. Nonlinear Sci. Numer. Simulat.}, 62:27--60, 2018.

\bibitem{PhysRevA.50.5077}
M.~J. Rakovi\ifmmode~\acute{c}\else \'{c}\fi{} and Shih-I Chu.
\newblock Approximate dynamical symmetry of hydrogen atoms in circularly polarized microwave fields.
\newblock {\em Phys. Rev. A}, 50:5077--5080, December 1994.

\bibitem{PhysRevA.47.R1612}
Kazimierz Rza$\mbox{\c{}}\dot{\text{z}}$ewski and Bernard Piraux.
\newblock Circular {R}ydberg orbits in circularly polarized microwave radiation.
\newblock {\em Phys. Rev. A}, 47:R1612--R1615, March 1993.

\bibitem{Simo1994}
Carles Sim\'{o}.
\newblock Averaging under fast quasiperiodic forcing.
\newblock In {\em Hamiltonian mechanics ({T}oru\'{n}, 1993)}, volume 331 of {\em NATO Adv. Sci. Inst. Ser. B: Phys.}, pages 13--34. Plenum, New York, 1994.

\bibitem{nature-569-75-77-2019}
Satya Undurti, Han Xu, Xiaoshan Wang, Atia Noor, William Wallace, Nicolas Douguet, Alexander Bray, I.~Ivanov, Klaus Bartschat, A.~Kheifets, Robert Sang, and Igor Litvinyuk.
\newblock Attosecond angular streaking and tunnelling time in atomic hydrogen.
\newblock {\em Nature}, 568:75--77, 04 2019.

\bibitem{ZakrzewskiDGR1993}
J.~Zakrzewski, D.~Delande, J.-C. Gay, and K.~Rzazewski.
\newblock Ionization of highly excited hydrogen atoms by a circularly polarized microwave field.
\newblock {\em Phys. Rev. A, Atomic, Molecular, and Optical Physics}, 47(4A):R2468--R2471, 1993.

\end{thebibliography}

\appendix

\section{Computation of 
\texorpdfstring{$\mathcal{A}^{\pm}\left(\omega\right)$}{A(om)}}
\label{appendix:intAw}
In this Appendix we explain how to derive the formula~\eqref{eq:Aomega} that
gives $\mathcal{A}^{\pm}(\omega)$ when $x_{0}(t) = x^{\pm}_{0}(t)$ is that of
the standard pendulum given in~\eqref{eq:StandardSeparatrices}. In particular we
take $x_{0}(t) = x_{0}^{+}(t)$, so we consider the external branch ($y_{0}(t) =
y_{0}^{+}(t) > 0$). For the internal branch ($y_{0}(t) = y_{0}^{-}(t) < 0$) the
process that leads to $\mathcal{A}^{-}(\omega)$ can be carried out in a similar
way, so we will not give the details here. 

Let $(x^{\pm}_{0}(s), y^{\pm}(s))$ be the parameterization of the separatrices
of the standard pendulum given in~\eqref{eq:StandardSeparatrices}, and let
$(x(s), y(s))$ denote specifically the external one, i.e.,  
\begin{displaymath}
    x_{0}(s) = x^{+}_{0}(s) = 4\arctan\left(\rme^{s}\right),
    \qquad
    y_{0}(s) = \dot{x}_{0}(s) = y^{+}_{0}(s) 
             =  \frac{2}{\cosh\:\!\! s}.
\end{displaymath}
Following~\cite{DelshamsG2000} we define
\begin{equation}
    \mathcal{J}(a,b) := \int_{-\infty}^{\infty} \rme^{\rmi a s} 
    \rme^{\rmi b \tfrac{x_{0}(s)-\pi}{2}} y_{0}(s) \df  s,\label{eq:Jab}
\end{equation}
that, as it is pointed out there, is real whenever $a$ and $b$ are real. With
the above definiton, integrating by parts the complex
expression~\eqref{eq:Aomega-complex} of the real formula~\eqref{eq:Aomega}, it
follows at once that 
\begin{equation}
    \mathcal{A}(\omega) = 
    \int_{-\infty}^{\infty} \rme^{\rmi\:\! \omega s}
    \left( 1 - \rme^{2\rmi x_{0}(s)} \right)\! \df s
    =
    \left[\mathcal{Q}_{\omega}(r)\right]_{-\infty}^{\infty} 
    + \frac{2}{\omega} \mathcal{J}(\omega, 4), \label{eq:Aomega-parts} 
\end{equation}
where 
\begin{displaymath}
    \mathcal{Q}_{\omega}(r)  = \frac{\rme^{\rmi \omega r}}{\rmi\:\!\omega}
    \left(1 - \rme^{\rmi 4 \tfrac{x_{0}(r)-\pi}{2}}\right) 
    = \frac{\rme^{\rmi \omega r}}{\rmi\:\!\omega}\left[ 1 - \left(
    \frac{1 + \rmi \sinh r}{\cosh r}\right)^{4}\right],
\end{displaymath}
so the ``integrated'' part at the rhs of~\eqref{eq:Aomega-parts} vanishes, since 
\begin{displaymath}
   \left[\mathcal{Q}_{\omega}(r)\right]_{-\infty}^{\infty} 
    = \lim\limits_{r\to\infty} \mathcal{Q}_{\omega}(r)
      - \lim\limits_{r\to-\infty} \mathcal{Q}_{\omega}(r)
\end{displaymath}    
and clearly $\mathcal{Q}_{\omega}(r)\to 0$ as $r\to\pm\infty$. Hence,
from~\eqref{eq:Aomega-parts}, $\mathcal{A}(\omega) = 2 \mathcal{J}(\omega,4) /
\omega$ and finally, we can use the formulas given in~\cite{DelshamsG2000},
Subsection 3.2, to compute~\eqref{eq:Jab} for $a = \omega$, $b = 4$. Thus, after
some algebra, it is seen that
\begin{displaymath}
    \mathcal{A}(\omega) = \frac{2}{\omega} \mathcal{J}(\omega, 4)  
        = \frac{4}{3}\omega^{3}\left(1 - \frac{2}{\omega^{2}}\right)\left[ 
            \frac{1}{\cosh\left(\dfrac{\pi\omega}{2}\right)} 
            -\frac{1}{\sinh\left(\dfrac{\pi\omega}{2}\right)} 
        \right]
        = \frac{16\pi}{3}\omega^{3}\left(1 - \frac{2}{\omega^{2}}\right)
        \frac{\rme^{\pi\omega / 2}}{1 - \rme^{2\pi\omega}},
\end{displaymath}
which is the equation~\eqref{eq:Aomega} corresponding to $c^{\pm} = 1$.

\section{Parameterization method}
\label{appendix:parmeth}
Let us assume that $z^{\ast}$ is a fixed point of the system $\dot{z} = F(z)$
with $F:\mathbb{R}^n\rightarrow\mathbb{R}^n$ regular enough. A simple strategy
to compute an invariant manifold ($W$) associated with the point $z^{\ast}$
involves first calculating the tangent linear space to $W$ at $z^{\ast}$, i.e.,
$T_{z^{\ast}} W$. Points on $T_{z^{\ast}} W$ that are sufficiently close to
$z^{\ast}$ can then be used as initial conditions for approximating $W$.

In this section, we focus specifically on computing the unstable manifold
$W^u(z^{\ast})$, which we assume to be 1-dimensional. The computation for
$W^s(z^{\ast})$ is analogous. Thus, we consider that $DF(z^{\ast})$ has a unique
eigenvalue $\lambda > 0$ with geometric multiplicity 1, and we denote its
associated eigenvector by~$v$.

A common strategy to compute a numerical approximation of the unstable invariant
manifold of $z^{\ast}$, $W^u(z^{\ast})$, is to take initial points $z^{\ast} +
sv$ and $z^{\ast} - sv$ with $s > 0$ sufficiently small, and integrate forward
in time. This approach provides a straightforward method to approximate the
desired invariant manifold. However, this strategy does not always yield
sufficiently accurate results. In such cases, the parameterization method
\cite{param1,param2,param3} can be used to obtain very efficient approximations
of these invariant manifolds to any desired order.

To find this approximation, we need to determine a parameterization $z = W(s)$,
with $s \in \mathbb{R}$, of the invariant manifold $W^u(z^{\ast})$ such that
$W(0) = z^{\ast}$. This way, the internal dynamics of the manifold are described
by the field $\dot{s} = f(s)$ with $f(0) = 0$, and the invariance equation is
given by:
\begin{equation} \label{eq:invariancia}
    F(W(s)) = DW(s)f(s).
\end{equation}
To solve \eqref{eq:invariancia}, we will find expansions for $W$ and $f$. in
particular, we consider:
\begin{displaymath} 
\left\{
\begin{aligned}
    &W(s) = z^* + \sum_{k\geq1}W_k(s) = z^* + \sum_{k\geq1}w_ks^k, \phantom{blabla} \text{with } w_1 = v,\\
    &f(s) = \sum_{k\geq1}f_k(s) = \sum_{k\geq1}\mathpzc{f}_k s^k, \phantom{blabla} \text{with } \mathpzc{f}_1 = \lambda,\\
\end{aligned}
\right.
\end{displaymath}
where $w_k\in\mathbb{R}^n$ and $\mathpzc{f}_k\in\R$ $\forall k$.

The goal is to compute $w_k$ and $\mathpzc{f}_k$ for $k > 1$ once the preceding terms $W_{<k}(s)$ and $f_{<k}(s)$ are known (i.e., given $w_1, \ldots, w_{k-1}$ and $\mathpzc{f}_1, \ldots, \mathpzc{f}_{k-1}$).

First, let us observe that we know the expression on the left side of the invariance equation \eqref{eq:invariancia} up to order $k-1$, i.e.,
\begin{equation*}
\left[F(W_{<k}(s))\right]_{<k} = \left[F\left(z^* + \sum_{i=1}^{k-1} W_i(s)\right)\right]_{<k}.
\end{equation*}
Thus, we will first calculate the homogeneous terms of degree $k$ for the composition $F(W_{k}(s))$ and the matrix product $DW_{<k}(s) f_{<k}(s)$, which we denote by:
\begin{equation*}
\begin{aligned}
    &\left[F(W_{<k}(s))\right]_k = \left[F\left(z^* + \sum_{i=1}^{k-1} W_i(s)\right)\right]_k, \\
    &\left[DW_{<k}(s) f_{<k}(s)\right]_k = \sum_{i=2}^{k-1} DW_i(s) f_{k+1-i}(s) = s^k \sum_{i=2}^{k-1} i\,w_i\, \mathpzc{f}_{k+1-i}.
\end{aligned}
\end{equation*}
In this way, we have:
\[
\begin{aligned}
    \left[F(W_{\leq k}(s))\right]_k &=  \left[F(W_{<k}(s))\right]_k + DF(z^*) W_k(s)\\
    & = \left[F(W_{<k}(s))\right]_k + DF(z^*) w_k s^k,
\end{aligned}
\]
and
\[
\begin{aligned}
    \left[DW_{\leq k}(s) f_{\leq k}(s)\right]_k 
    &= \left[DW_{<k}(s) f_{<k}(s)\right]_k + DW_k(s) f_1(s) + DW_1(s) f_k(s)\\
    & = \left[ \sum_{i=2}^{k-1} \left(i\,w_i\, \mathpzc{f}_{k+1-i}\right)
    +\lambda k w_k + v \mathpzc{f}_k
    \right] s^k.
\end{aligned}
\]

This leads to the cohomological equation of order $k$ for $w_k$ and $\mathpzc{f}_k$:
\begin{equation} \label{eq:coho}
    \big[\big(DF(z^*) - \lambda k\,Id\big) w_k  - v \mathpzc{f}_k\big] s^k = -E_k(s),
\end{equation}
where $Id$ denotes the identity matrix and $E_k(s) = \left[F(W_{<k}(s))\right]_k - s^k\sum_{i=2}^{k-1} \left(i\,w_i\, \mathpzc{f}_{k+1-i}\right)$ is the error of order $k$.

Note that system \eqref{eq:coho} is a linear system with $n$ equations and $n+1$ unknowns, giving us one degree of freedom. Therefore, we will assume all the terms $\mathpzc{f}_i = 0$ $\forall i > 1$, which is known as a normal form style of parameterization (see \cite{paramLlibre} for details), resulting in the system:
\begin{equation} \label{eq:coho2}
    \big[\big(DF(z^*) - \lambda k\,Id\big) w_k  \big] s^k = -E_k(s)= -\left[F(W_{<k}(s))\right]_k.
\end{equation}
Thus, to obtain successive orders of the invariant manifold, we will simply compute the error of order $k$ using automatic differentiation \cite{AutoDif}, and then solve the linear system \eqref{eq:coho2}.

\section{Singularities of the separatrices of the amended  pendulum}\label{appendix:singularity}
Consider the Hamiltonian of the amended pendulum 
\begin{equation}\label{eq:HQP}
P(x,y,\varepsilon)=\frac{y^2}{2}-\frac{2}3\varepsilon^2y^3+\cos x -1
\end{equation}
and consider the separatrices of the hyperbolic equilibrium point $(0,0)$
defined by  $P(x,y,\varepsilon)=0$. We know that when $\varepsilon =0$, the
separatrices coincide and become singular at $t=\pm \rmi\sfrac{\pi}{2}$. So we
expect, by continuity, the singularities of the separatrices
$P(x,y,\varepsilon)=0$ of the amended pendulum to be close to them. We will
focus on the separatrix such that at $t=0$ it passes through the point
$(x,y)=(\pi,y_0(\varepsilon))$ (when $\varepsilon =0$, $y_0=2$).

So, let us consider the change of time $t=-\rmi s$ (we will see below that the
singularity takes place for $t^*=\rmi(-s^*)=\sfrac{\pi}{2}+\delta
(\varepsilon)$), then the system of ODE associated with \eqref{eq:HQP}
\begin{align*}
\frac{\df x}{\rmd t} &= y-2\varepsilon^2y^2,\\
\frac{\rmd y}{\rmd t} &= \sin x,
\end{align*}
becomes
\begin{align*}
\frac{\rmd X}{\rmd s} &= -(Y-2\varepsilon^2Y^2),\\
\frac{\rmd Y}{\rmd s} &= \sinh\:\! X,
\end{align*}
which is a Hamiltonian system with Hamiltonian
\begin{displaymath}
Q(x,y,\varepsilon)=-\frac{Y^2}{2}+\frac{2}3\varepsilon ^2Y^3-\cosh\:\!\! x +1
\end{displaymath}
\begin{figure}[t!]
    \centering
    \includegraphics[trim={0mm 0mm 0mm 0mm},
                     clip,width=0.48\textwidth]
                     {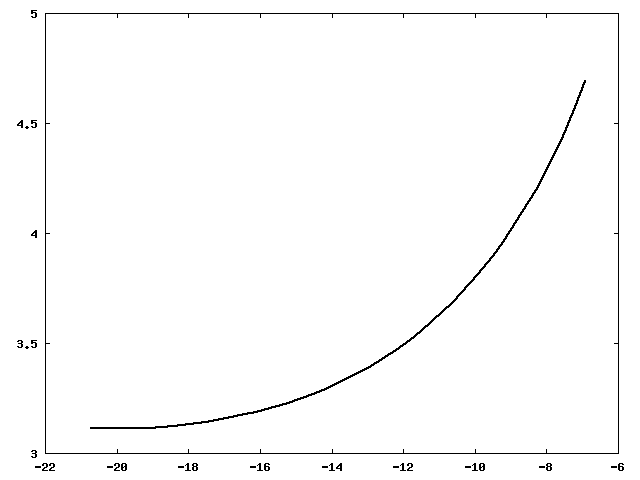}
\caption{%
    Curve $\left(\ln K,\tfrac{\delta (K)}{(K|\ln K|)^{2.051}}\right)$ where $K$
    ranges decreasing from $0.001$ to $10^{-9}$.
}
\label{fig:singtat185d}
\end{figure}
The following properties can be obtained easily:
\begin{align*}
   (P1) & & P(\pi,y,\varepsilon) &= \frac{y^2}{2}-\frac{2}3\varepsilon ^2y^3-2, &
   \nonumber\\
   (P2) & & \cos x &= 1-\frac{y^2}{2}+\frac{2}{3}\varepsilon^2y^3, &
   \nonumber\\
   (P3) & & \sin^2 x &= \left(\frac{2}{3}\varepsilon ^2y-\frac 12\right)
            P\left(\pi,y,\varepsilon\right) y^2 & \nonumber\\
        & & &= -\frac {2}{3}\varepsilon^2\biggl(
            \frac{2}{3}\varepsilon ^2y-\frac 12\biggr)
            \biggl(y-y_0(\varepsilon)\biggr)
            \biggl(y-y_1(\varepsilon)\biggr)
            \biggl(y-y_2(\varepsilon)\biggr)y^2,\\ 
   (P4) & & y(-\rmi s) &= Y(s)\\    
   (P5) & & \sinh\:\! X(s) &= -\rmi\sin x(-\rmi s) &\nonumber     
 \end{align*}  
with, 
 \begin{gather*}
          y_0(\varepsilon) = \dfrac{\sqrt{3}}{\cos\left( 
              \dfrac{\pi-\arccos(4\sqrt{3}\varepsilon^2)}{3}\right)},\qquad 
          y_1(\varepsilon) = \dfrac{\sqrt{3}}{\cos\left( 
              \dfrac{\pi+\arccos(4\sqrt{3}\varepsilon^2)}{3}\right)},\\
          y_2(\varepsilon) = -\dfrac{\sqrt{3}}{\cos\left(
              \dfrac{\arccos(4\sqrt{3}\varepsilon^2)}{3}\right)}.
\end{gather*}              
So, we obtain
\begin{displaymath}
\frac{\df Y}{\df s} (s) = \sinh\:\! X(s) 
        = -\sqrt{-\bigl(\sfrac{1}{2}-\sfrac{2\varepsilon^2Y(s)}{3}\bigr)
                    \bigl(\Scale[0.85]{Y(s)-y_0(\varepsilon)}\bigr)
                    \bigl(\sfrac{2\varepsilon^2Y(s)}{3} -
                        \sfrac{2\varepsilon^2 y_1(\varepsilon)}{3}\bigr)
                    \bigl(\Scale[0.75]{Y(s)-y_2(\varepsilon)}\bigr)\Scale[0.85]{Y^2(s)} 
      }, 
\end{displaymath}
that is
\begin{displaymath}
-s^{\ast}=\mathlarger{\mathlarger{\int}}_{_{y_{0}(\varepsilon)}}^{
                      \scriptscriptstyle\infty}
\dfrac{\df Y}
{\sqrt{-\bigl(\sfrac{1}{2}-\sfrac{2\varepsilon^2Y}{3}\bigr)
                    \bigl(\Scale[0.85]{Y-y_0(\varepsilon)}\bigr)
                    \bigl(\sfrac{2\varepsilon^2Y}{3} -
                        \sfrac{2\varepsilon^2 y_1(\varepsilon)}{3}\bigr)
                    \bigl(\Scale[0.75]{Y-y_2(\varepsilon)}\bigr)\Scale[0.85]{Y^2} 
      }, 
}
\end{displaymath}
which is a convergent improper integral or, equivalently, through the
change of variable $Y=2/v$, we~have
\begin{displaymath}
-s^{\ast} = \mathlarger{\mathlarger{\int}}_{_0}^{\sfrac{2}{y_{0}(\varepsilon)}}
\!\!\!\!
\dfrac{v \df v}
{\sqrt{\bigl(v-\sfrac{8\varepsilon^2}{3}\bigr)
       \bigl(v-\sfrac{2}{y_0(\varepsilon)}\bigr)
       \bigl(v-\sfrac{2}{y_1(\varepsilon)}\bigr)
       \bigl(v-\sfrac{2}{y_2(\varepsilon)}\bigr)}}.     
\end{displaymath}
For each given value of $\varepsilon>0$, we have computed, independently, both
improper integrals to double check that both  values of $s^*$ coincide. From the
numerical computations done, we conclude that  $s^*$ has real and imaginary
parts, that is, $-s^*=-(\Re s^*+\rmi \Im s^*)=\sfrac{\pi}{2}+\delta
(\varepsilon)$, where $\Re (-s^*)=\sfrac{\pi}{2}+\delta (\varepsilon)$ and $\Im
\delta (\varepsilon)=O(\varepsilon ^2)$. 

Our next purpose is to have a fit estimate of $\Re (-s^*)=\sfrac{\pi}{2}+\delta
(\varepsilon)$ when varying $\varepsilon$, or, more precisely, of
$\delta(\varepsilon)=|\Re (s^*)|-\sfrac{\pi}{2}$. The following fit is obtained
\begin{displaymath}
\delta(K)=A(K|\ln K|)^{2.051}
\end{displaymath} 
as can be seen in Figure \ref{fig:singtat185d}, where the points $\left(\ln
K,\tfrac{\delta (K)}{(K|\ln K|)^{2.051}}\right)$ are plotted and we see the
tendency towards a constant $A$. So we can infer, taking into account the
relation $K=3\varepsilon ^4$,  that $\delta(\varepsilon)$ is of order
$\varepsilon ^8\ln ^2 \varepsilon$.
\end{document}